\documentclass[11pt,a4paper]{report}

\usepackage{amsmath}
\usepackage{amsthm}
\usepackage{amsfonts}
\usepackage{amscd}
\usepackage{amssymb}
\usepackage{mathrsfs}
\usepackage[dvips]{geometry}
\usepackage[all,cmtip]{xy}
\newtheorem{thm}{Theorem}[section]
\newtheorem*{thm*}{Theorem}
\newtheorem{cor}[thm]{Corollary}
\newtheorem*{cor*}{Corollary}
\newtheorem{lem}[thm]{Lemma}
\newtheorem*{lem*}{Lemma}
\newtheorem{prop}[thm]{Proposition}

\theoremstyle{definition}
\newtheorem{defn}{Definition}[section]
\newtheorem*{defn*}{Definition}
\newtheorem*{conjecture}{Conjecture}
\newtheorem*{condition}{Condition}

\theoremstyle{remark}
\newtheorem{rem}{Remark}[section]
\newtheorem{example}{Example}[section]

\newcommand{\To}{\longrightarrow}
\newcommand{\ra}{\rightarrow}
\newcommand{\calQ}{\mathcal Q}
\newcommand{\U}{\mathcal U}
\newcommand{\Q}{\mathbb Q}
\newcommand{\R}{\mathbb R}
\newcommand{\C}{\mathbb C}
\newcommand{\cc}{\mathfrak c}

\newcommand{\XX}{\mathfrak X}

\newcommand{\Z}{\mathbb Z}
\newcommand{\PP}{\mathbb P}
\newcommand{\F}{\mathcal F}
\newcommand{\SC}{\mathscr H}
\newcommand{\E}{\mathscr E}
\newcommand{\HC}{\mathbb H}
\newcommand{\HH}{\mathfrak H}
\newcommand{\A}{\mathscr A}
\newcommand{\D}{\mathscr D}
\newcommand{\DD}{\mathbb D}
\newcommand{\SB}{\mathbb S}
\newcommand{\wt}{\widetilde}
\newcommand{\ol}{\overline}

\newcommand{\veps}{\varepsilon}

\newcommand{\z}{z}
\newcommand{\zc}{\bar z}
\newcommand{\I}{\mathrm i}

\newcommand{\G}{\mathfrak G}

\newcommand{\Hom}{\mathrm{Hom}}
\newcommand{\HHom}{\mathbf{Hom}}

\newcommand{\Fl}{\mathrm{Fl}}
\newcommand{\fl}{\mathrm{fl}}
\newcommand{\Int}{\mathrm{int}}
\newcommand{\cycle}{\mathrm{cycle}}
\newcommand{\const}{\mathrm{const}}
\newcommand{\Stab}{\mathrm{Stab}}

\newcommand{\vol}{\mathrm{vol}}
\newcommand{\alg}{\mathrm{alg}}
\newcommand{\analyt}{\mathrm{an}}
\newcommand{\eq}{\; = \;}
\newcommand{\parab}{\mathrm{par}}
\newcommand{\sym}{\mathrm{sym}}

\DeclareMathOperator{\res}{res}
\DeclareMathOperator{\Tr}{Tr}
\DeclareMathOperator{\Der}{Der}

\DeclareMathOperator{\id}{id}

\DeclareMathOperator{\Sym}{Sym}
\DeclareMathOperator{\ord}{ord}

\DeclareMathOperator{\smooth}{smooth}

\DeclareMathOperator{\spec}{Spec}
\DeclareMathOperator{\Div}{div}
\DeclareMathOperator{\cl}{cl}

\DeclareMathOperator{\cone}{cone}
\DeclareMathOperator{\image}{Im}
\DeclareMathOperator{\kernel}{Ker}
\DeclareMathOperator{\cokernel}{Coker}

\DeclareMathOperator{\dist}{dist}

\renewcommand{\Im}{\mathrm{Im}\;}

\author{Anton Mellit \footnote{mellit@gmail.com}}

\title{Higher Green's functions for modular forms}

\begin{document}

\bibliographystyle{alpha}
\maketitle

\begin{abstract}
Higher Green functions are real-valued functions of two variables on the
upper half plane which are bi-invariant under the action of a congruence
subgroup, have logarithmic singularity along the diagonal, but instead of
the usual equation $\Delta f=0$ we have equation $\Delta f = k(1-k) f$.
Here $k$ is a positive integer. Properties of these functions are related
to the space of modular forms of weight $2k$. In the case when there are no
cusp forms of weight $2k$ it was conjectured that the values of the
Green function at points of complex multiplication are algebraic multiples
of logarithms of algebraic numbers. We show that this conjecture can be
proved in any particular case if one constructs a family of elements of
certain higher Chow groups on the power of a family of elliptic curves. These families have to satisfy certain properties. A different family of elements of Higher Chow groups is
needed for a different point of complex multiplication. We give an example
of such families, thereby proving the conjecture for the case when the group is $PSL_2(\Z)$, $k=2$ and one of the arguments is $\I$.
\end{abstract}

\chapter*{Introduction}

The subject of the present thesis is higher Green's functions. For any integer $k>1$ and subgroup $\Gamma\subset PSL_2(\Z)$ of finite index there is a unique function $G_k^{\HH/\Gamma}$ on the product of the upper half plane $\HH$ by itself which satisfies the following conditions:
\begin{enumerate}
\item $G_k^{\HH/\Gamma}$ is a smooth function on $\HH\times\HH\setminus \{(\tau,\gamma \tau)\,|\, \gamma\in\Gamma,\, \tau\in \HH\}$ with values in $\R$.
\item $G_k^{\HH/\Gamma}(\gamma_1 \tau_1, \gamma_2 \tau_2) = G_k^{\HH/\Gamma}(\tau_1, \tau_2)$ for all $\gamma_1, \gamma_2\in \Gamma$.
\item $\Delta_i G_k^{\HH/\Gamma} = k(1-k) G_k^{\HH/\Gamma}$, where $\Delta_i$ is the hyperbolic Laplacian with respect to the $i$-th variable, $i=1,2$.
\item $G_k^{\HH/\Gamma}(\tau_1, \tau_2) = m\log|\tau_1-\tau_2|^2 + O(1)$ when $\tau_1$ tends to $\tau_2$ ($m$ is the order of the stabilizer of $\tau_2$, which is almost always $1$).
\item $G_k^{\HH/\Gamma}(\tau_1, \tau_2)$ tends to $0$ when $\tau_1$ tends to a cusp.
\end{enumerate}

This function is called the Green function. It is necessarily symmetric, 
\[
G_k^{\HH/\Gamma}(\tau_1, \tau_2) = G_k^{\HH/\Gamma}(\tau_2, \tau_1).
\]

Such functions were introduced in paper \cite{GZ}. Also it was conjectured in \cite{GZ} and \cite{GKZ} that these functions have ``algebraic'' values at CM points. A particularly simple formulation of the conjecture is in the case when there are no cusp forms of weight $2k$ for the group $\Gamma$:

\begin{conjecture}[1]
Suppose there are no cusp forms of weight $2k$ for $\Gamma$. Then for any two CM points $\tau_1$, $\tau_2$ of discriminants $D_1$, $D_2$ there is an algebraic number $\alpha$ such that
\[
G_k^{\HH/\Gamma}(\tau_1, \tau_2) \eq (D_1 D_2)^{\frac{1-k}2} \log \alpha.
\]
\end{conjecture}

The main result of this thesis is a general approach for proving this conjecture and an actual proof for the case $\Gamma = PSL_2(\Z)$, $k=2$, $\tau_2=\I$\footnote{We use ``$\I$'' to denote $\sqrt{-1}$ and ``$i$'' for other purposes.}. The number $\alpha$ in the latter case is represented as the intersection number of a certain higher Chow cycle on the elliptic curve corresponding to the point $\tau_1$ and an ordinary algebraic cycle which has a certain prescribed cohomology class.

The general approach can be formulated as follows. For an elliptic curve $E$ we consider the higher Chow group $CH^k(E^{2k-2}, 1)$ on the product of $E$ by itself $2k-2$ times. The elements of this group are represented by so-called ``higher cycles'', which are formal linear combinations
\[
\sum_i (W_i, f_i),
\]
where $W_i$ is a subvariety of $E^{2k-2}$ of codimension $k-1$ and $f_i$ is a non-zero rational function on $W_i$ such that the following condition holds in the group of cycles of codimension $k$:
\[
\sum_i \Div f_i = 0.
\]
Denote the abelian group of higher cycles by $Z^k(E^{2k-2}, 1)$, so that $CH^k(E^{2k-2}, 1)$ is its quotient (the relations defining $CH^k(E^{2k-2}, 1)$ will be explained in Section \ref{chow_groups}).

We have the Abel-Jacobi map
\[
AJ^{k,1}:CH^k(E^{2k-2},1) \To \frac{H^{2k-2}(E^{2k-2}, \C)}{F^k H^{2k-2}(E^{2k-2},\C)+2\pi\I H^{2k-2}(E^{2k-2},\Z)},
\]
and there is a canonical cohomology class $[\theta]\in H^{2}(E^{2}, \C)$, namely the one represented by the form
\[
\theta = \frac{\omega \otimes \bar\omega + \bar\omega \otimes \omega}{\int_{E} \omega\wedge\bar{\omega}},
\]
where $\omega$ is a holomorphic differential $1$-form on $E$. On may notice that for any element $[x]\in CH^k(E^{2k-2},1)$ there is a perfectly defined number
\[
\Re \left( \I^{k-1}(AJ^{k,1}[x], [\theta]^{k-1}) \right).
\]
So we hope that for those cases for which the conjecture above is formulated, for a fixed CM point $\tau_2$, there exists a family of higher cycles $\{x_s\}_{s\in S}$, $x_s\in CH^k(E_s^{2k-2},1)$ for a family of elliptic curves $\{E_s\}_{s\in S}$ ($S$ is an algebraic variety over $\C$) which ``computes'' the values of the Green function in the way described above. This approach appeared from an attempt to understand the discussion of the algebraicity conjecture in \cite{ZH}.

\vspace{2mm}
To be more specific I introduce the following $3$ classes of functions.

Let $A$ be a subgroup of $\C$ (we will usually take $A=\Z$, $A=\Q$ or $A=\I^{k-1} \R$). Consider a holomorphic multi-valued function $f$ on the upper half plane which is allowed to have isolated singularities and is defined up to addition of polynomials of $\tau$ with coefficients in $2\pi\I A$ of degree not greater than $2k-2$. This means that for a small disk $U$ which does not contain the singularities of $f$ one has an element of $O_{\mathrm{an}}(U)/2 \pi\I V_{2k-2}^A$, where $V_{2k-2}^A$ is the abelian group of polynomials with coefficients in $A$ of degree not greater than $2k-2$. Suppose $f$ transforms like a modular form of weight $2-2k$ with respect to $\Gamma$, i.e. $f\equiv f|_{2-2k} \gamma \mod 2\pi\I V_{2k-2}^A$ for any $\gamma\in\Gamma$. Denote the abelian group of all such functions by $M^{A}_{2-2k}(\Gamma)$. Note that if we put some growth condition on $f$ near the singularities, it will imply that the singularity of $f$ at a point $\tau_0$ will be of the form $F(\tau) \log (\tau-\tau_0)$ for $F\in V_{2k-2}^A$.

Similarly, let $M^{A}_{0, k}(\Gamma)$ be the space of multi-valued functions $f$ with isolated singularities satisfying $\Delta f = k(1-k) f$, almost holomorphic (expressible as polynomials of $\frac{1}{\tau-\bar\tau}$ with holomorphic coefficients), invariant (in weight $0$) with respect to $\Gamma$, defined up to addition of functions of the form $\sum_{i=0}^{2k-2} \alpha_i r_i(\tau)$, where $\alpha_i\in 2\pi\I \frac{i! (2k-2-i)!}{(k-1)!} A $ and the $r_i$ are defined by the equation
\[
\sum_{i=0}^{2k-2} r_i(\tau) X^i \;=\; \left(\frac{(X-\tau)(X-\bar\tau)}{\tau-\bar\tau}\right)^{k-1}.
\]

Finally, let $M^{A}_{2k}(\Gamma)$ be the space of (single-valued) holomorphic functions $f$ with isolated singularities which transform like modular forms of weight $2k$ with respect to $\Gamma$ and the $V^{\C}_{2k-2}$-valued differential form $\omega_f=f(\tau)(X-\tau)^{2k-2} d\tau$ satisfies the following two conditions: 
\begin{enumerate}
 \item If we integrate $\omega_f$ around a singularity we obtain an element of the abelian group $2\pi\I (2k-2)! V_{2k-2}^A$.
 \item The cohomology class $[\omega_f]\in H^1(\Gamma, V_{2k-2}^{\C/2\pi\I(2k-2)! A})$ of $\omega_f$ is trivial.
\end{enumerate}

Then there is a commutative diagram:
\[
\xymatrix{
M^{A}_{2-2k}(\Gamma) \ar[r]^{\delta^{k-1}} \ar[rd]_{\left(\frac{d}{d\tau}\right)^{2k-1}}& 
M^{A}_{0, k}(\Gamma) \ar[d]^{\delta^k} \\
& M^{A}_{2k}(\Gamma)
}
\]
Here the horizontal arrow is the operator $\delta_{-2}\cdots\delta_{2-2k}$ and the vertical one is $\delta_{2k-2} \cdots \delta_{0}$ ($\delta_w = \frac{d}{d\tau} + \frac{w}{\tau-\bar\tau}$). Also the horizontal arrow is an isomorphism and the vertical one is surjective with kernel the finite group $(V_{2k-2}/V_{2k-2}^A)^\Gamma$. In particular if $A$ is $\Q$ or $\R$, then the horizontal arrow is also an isomorphism. When $A=\R$ we also have a fourth group in our system. Denote by $C^{\omega}_{0, k}(\Gamma)$ the space of real-analytic (single-valued) functions with isolated singularities satisfying the same differential equation as for the space $M^{A}_{0, k}(\Gamma)$ and invariant (in weight $0$) with respect to the action of $\Gamma$. Then the following diagram is commutative and all its arrows are isomorphisms:
\[
\xymatrix{
M^{\I^{k-1}\R}_{2-2k}(\Gamma) \ar[r]^{\delta^{k-1}} \ar[d]_{\left(\frac{d}{d\tau}\right)^{2k-1}}& 
M^{\I^{k-1}\R}_{0, k}(\Gamma) \ar[d]^{2\Re(\cdot)}\\
M^{\I^{k-1}\R}_{2k}(\Gamma) &
C^{\omega}_{0,k}(\Gamma) \ar[l]^{\delta^k}
}
\]

Now we give two ways to construct elements of the spaces above. The first way takes the Green function as input. Fix a point $\tau_0\in\HH$. The function $G_k^{\HH/\Gamma}(\tau, \tau_0)$ belongs to the space $C^{\omega}_{0,k}(\Gamma)$. Therefore our construction produces elements in the spaces $M^{\I^{k-1}\R}_{2-2k}(\Gamma)$, $M^{\I^{k-1}\R}_{0, k}(\Gamma)$, $M^{\I^{k-1}\R}_{2k}(\Gamma)$. Denote by $g_{\tau_0}=g_{\tau_0, k}^{\HH/\Gamma}$ the corresponding element in $M^{\I^{k-1}\R}_{2k}(\Gamma)$. We prove that $g_{\tau_0}$ is meromorphic, zero at the cusps, and its principal part at $\tau_0$ is $m (-1)^{k-1} (k-1)! Q_{\tau_0}(\tau)^{-k}$ (we denote $Q_{\tau_0}(X) = \frac{(X-\tau_0)(X-\bar\tau_0)}{\tau_0-\bar\tau_0}$). The integral around $\tau_0$ of $\omega_{g_{\tau_0}} = g_{\tau_0}(\tau)(X-\tau)^{2k-2} d\tau$ is $2\pi\I \frac{(2k-2)!}{(k-1)!} Q_{\tau_0}(X)^{k-1}$. Therefore the product $(k-1)!D_0^{\frac{k-1}2} g_{\tau_0}$ satisfies all the requirements of $M^{\Z}_{2k}(\Gamma)$ except, possibly, the last one. The last requirement is automatically satisfied in the case $S_{2k}(\Gamma)=\{0\}$, at least up to torsion, i.e. there exists an integer $N_1$ (depending only on $k$ and $\Gamma$) such that $N_1 D_0^{\frac{k-1}2} g_{\tau_0}\in M^{\Z}_{2k}(\Gamma)$. The lift of $N_1 D_0^{\frac{k-1}2} g_{\tau_0}$ to $M^{\Z}_{2-2k}(\Gamma)$ and $M^{\Z}_{0, k}(\Gamma)$ is defined up to a finite group, therefore its product with certain integer number $N_2$ (which depends only on $\Gamma$ and $k$) is perfectly defined. Let $N=N_1 N_2$. We see that the construction canonically gives a function $\widehat{G}_{k,\tau_0}^{\HH/\Gamma} \in \frac{1}{N} D_0^{\frac{1-k}2} M^{\Z}_{0, k}(\Gamma)$ such that $\delta^k \widehat{G}_{k,\tau_0}^{\HH/\Gamma} = g_{\tau_0}$. If, moreover, $\tau$ is a CM point of discriminant $D$, then $\widehat{G}_{k,\tau_0}^{\HH/\Gamma} (\tau)$ is defined up to $2\pi\I \frac{1}{N} (D_0 D)^{\frac{1-k}2} \Z$. This is called ``the lifted value of the Green function'' and we conjecture that it equals $\frac{1}{N} (D_0 D)^{\frac{1-k}2} \log \alpha$ for some $\alpha\in\overline{Q}^\times$. We emphasize that the Green function produces elements of spaces $M^{\I^{k-1}\R}_{2-2k}(\Gamma)$, $M^{\I^{k-1}\R}_{0, k}(\Gamma)$, $M^{\I^{k-1}\R}_{2k}(\Gamma)$, and only if $\tau_0$ is a CM point we can ``lift'' these elements to spaces $M^{\Z}_{2-2k}(\Gamma)$, $M^{\Z}_{0, k}(\Gamma)$, $M^{\Z}_{2k}(\Gamma)$. 

There is another way to obtain functions as above, which always produces elements of spaces $M^{\Z}_{2-2k}(\Gamma)$, $M^{\Z}_{0, k}(\Gamma)$, $M^{\Z}_{2k}(\Gamma)$. Suppose we have a family of elliptic curves $\{E_s\}_{s\in S}$ over an algebraic variety $S$ defined over $\C$ and an algebraic family of higher cycles $x = \{x_s\in Z^k(E_s^{2k-2},1)\}_{s\in S}$. Let $x_s = \sum_i(W_{i,s}, f_{i,s})$ (the dimension of $W_{i,s}$ is $k-1$ and $\sum_i \Div f_{i,s} = 0$). Suppose also a map $\varphi: S\rightarrow \HH/\Gamma$ is given which is dominant, i.e. its image is $\HH/\Gamma$ without a finite number of points, and suppose that $\varphi$ makes the following diagram commutative ($j$ is the $j$-invariant):
\[
\xymatrix{
S \ar[r]_{\varphi} \ar[rd]_{j} & \HH/\Gamma \ar[d] \\
& \HH/PSL_2(\Z).
}
\]

\begin{defn*}
A triple $\XX =(\{E_s\}_{s\in S}, \{x_s\}_{s\in S}, \varphi)$ as above is called \emph{modular} if for any two points $s_1, s_2$ such that $\varphi(s_1) = \varphi(s_2)$ there exists an isomorphism $\rho: E_{s_1} \rightarrow E_{s_2}$ with the following properties:
\begin{enumerate}
 \item $\rho^{2k-2}(W_{i,s_1}) = W_{i,s_2}$,
 \item $\rho^{2k-2}(f_{i,s_1}) = \beta_i f_{i,s_2}$ for $\beta_i\in \C^\times$.
\end{enumerate} 
\end{defn*}

If we have a modular triple then we construct an element of $M^{\Z}_{2-2k}(\Gamma)$ as follows. Let $\tau\in\HH$ be such that its projection to $\HH/\Gamma$ belongs to the image of $\varphi$, say $\tau = \varphi(s)$. Choose a differential form $\omega$ on $E_s$ and suppose its period lattice is generated by $\Omega_1$ and $\Omega_2$ with $\frac{\Omega_2}{\Omega_1} = \tau$. Then put 
\[
A_\XX(\tau) = \frac{1}{\Omega_1^{2k-2}} \langle AJ^{k,1}[x_s], [\omega^{\otimes 2k-2}]\rangle.
\]
In this case the Abel-Jacobi map reduces to the following integral:
\[
A_\XX(\tau) = 2\pi\I \frac{1}{\Omega_1^{2k-2}} \int_\xi \omega^{\otimes 2k-2},
\]
where $\xi$ is a smooth $2k-2$-chain whose boundary is $\sum_i f_{i,s}^*[0,\infty]$ and $[0,\infty]$ is a path from $0$ to $\infty$ on $\C P^1$. It is clear that $A_\XX(\tau)$ thus defined does not depend on the choice of $s$ with $\varphi(s)=\tau$, and $A_\XX\in M^{\Z}_{2-2k}(\Gamma)$. We will show that 
\[
\delta^{k-1} A_\XX (\tau) = \frac{(2k-2)!}{(k-1)! \left(\int_{E_s} \omega \wedge \ol\omega\right)^{k-1}} \;\langle\, AJ^{k,1}[x],\, [\omega^{\otimes k-1} \overset{sym}\otimes \ol{\omega}^{\otimes k-1}] \,\rangle.
\]
If $\tau$ is a CM point with minimal equation $a \tau^2 + b \tau + c = 0$, $D=b^2-4ac$ and $\varphi(s)=\tau$, then there are endomorphisms of $E_s$ which act on the tangent space as the multiplications by $a\tau$ and $a\bar\tau$. Denote their graphs by $Y_{a\tau}$ and $Y_{a\bar\tau}$ correspondingly. Then one can check that the Poincar\'e dual cohomology class of the difference $Y_{a\tau} - Y_{a\bar\tau}$ is represented by the form $-2 \sqrt{D} \frac{\omega \overset{sym}\otimes \ol{\omega}}{\int_{E_s} \omega \wedge \ol\omega}$. Therefore we can construct a variety whose Poincar\'e dual class is represented by
\[
(-1)^{k-1} D^{\frac{k-1}2} \frac{(2k-2)!}{(k-1)! \left(\int_{E_s} \omega \wedge \ol\omega\right)^{k-1}} \;\omega^{\otimes k-1}\overset{sym}\otimes \ol{\omega}^{\otimes k-1},
\]
namely, we take the product of $k-1$ copies of $Y_{a\tau} - Y_{a\bar\tau}$ for each possible splitting of the product $E_s^{2k-2}$ into pairs and add them up. Note that there are precisely $\frac{(2k-2)!}{(k-1)! 2^{k-1}}$ such splittings. Denote the variety obtained in this way by $Z_s$. Then we obtain
\[
\delta^{k-1} A_\XX (\tau) \;=\; (-1)^{k-1} D^{\frac{1-k}2} \log \left( x_s \cdot Z_s\right),
\]
where the intersection number $x_s \cdot Z_s$ is defined as
\[
x_s\cdot Z_s \;=\; \prod_i \prod_{p\in W_{i,s}\cap Z_s} f_{i,s}(p)^{\ord_p W_{i,s} \cdot Z_s}
\]
if $Z_s$ does not intersect the divisors of $f_{i,s}$, singularities of $W_{i,s}$ and intersects $W_{i,s}$ properly. If this is not the case one can still define the intersection number, for example, by ``shifting'' $Z_s$ using the the addition law.

The second construction gives examples of functions in $M_{0,k}^{\Z}$ whose values at CM points are ``algebraic''. We prove that for such $A_\XX$, obtained via the second construction, the derivative $\left(\frac{d}{d\tau}\right)^{2k-1} A_\XX (\tau)$ is meromorphic. In fact we give a method to compute this derivative purely algebraically.

\begin{defn*}
A meromorphic modular form $f$ which belongs to $M_{0,k}^{\Z}$ is called \emph{geometrically representable} if $f=\left(\frac{d}{d\tau}\right)^{2k-1} A_\XX (\tau)$ for some modular triple $(\{E_s\}_{s\in S}, \{x_s\}_{s\in S}, \varphi)$ with $\{E_s\}$, $\varphi$ and $\{x_s\}$ are defined over $\ol{Q}$.
\end{defn*}

Then we have
\begin{lem*}
Suppose $\Gamma$ is a congruence subgroup, $k>1$ and $f\in M_{0,k}^{\Z}(\Gamma)$ is such that $(-1)^{k-1} N_1 \delta^k f$, for some integer $N_1$, is geometrically representable by a modular triple $(\{E_s\}_{s\in S}, \{x_s\}_{s\in S}, \varphi)$. Suppose $\tau\in\HH$ is a CM point which belongs to the image of $\varphi$. Let $D$ be the discriminant of $\tau$. Then $D^{\frac{k-1}2} f(\tau)\equiv \frac{1}{N_1} \log \alpha \mod \frac{2\pi\I}{N}$, where $\alpha$ is an algebraic number, which can be computed as the intersection $x_s \cdot Z_s$ as above with $\varphi(s) = \tau$. Here $N_2$ is the exponent of the group $H^0(\Gamma, V_{2k-2}^{\Q/\Z})$, $N=N_1 N_2$.
\end{lem*}

In particular, if a multiple of $D_0^{\frac{k-1}2} g_{k, \tau_0}$ is geometrically representable, then the conjecture is true for $\tau_1=\tau_0$ and $\tau_2$~--- any other CM point except, possibly, a finite number of points (those which do not belong to the image of $\varphi$). We give an example when this occurs. Let $\Gamma = PSL_2(\Z)$, $k=2$. Let $S$ be the open subset of $\C\times\C$ of pairs $(a,b)$ satisfying $4 a^3 + 27 b^2 \neq 0$ and $b\neq 0$. Let $\{E_s\}_{s\in S}$ be the Weierstrass family, which is defined by the projective version of the equation $y^2 = x^3 + a x + b$. Let $W_s$ be the subvariety of $E_s\times E_s$ which consists of $(x_1, y_1, x_2, y_2)$ with $x_1 + x_2 = 0$. Let $f_s$ be the function $y_1-\I y_2$. It is easy to check that $(W_s, f_s)\in Z^2(E_s\times E_s, 1)$. Let $\varphi$ be the $j$-invariant. One can verify that this gives a modular triple. Denote it by $\XX_{-4}$ (since the discriminant of the point $\tau=\I$, which is the only point which does not belong to the image of $\varphi$, is $-4$). We prove
\begin{thm*}
Consider the following modular form of weight $4$:
\[
-\frac12 \sqrt{-4}\; \delta^2 G_2^{\HH/PSL_2(\Z)}(\tau, \I) \;=\; (2\pi\I)^2 \sqrt{-4}\; \frac{432 E_4(\tau)}{j(\tau) - 1728}.
\]
This form is geometrically representable by the modular triple $\XX_{-4}$.
\end{thm*}

This implies the conjecture for $G_2^{\HH/PSL_2(\Z)}(\tau, \I)$ (keep in mind that $G = 2\Re \widehat{G}$), namely
\begin{cor*}
For any CM point $\tau$ which is not equivalent to $\I$ one has
\[
\sqrt{-4 D} \;\widehat{G}_2^{\HH/PSL_2(\Z)}(\tau, \I) \equiv 2 \log\left((W_s, f_s)\cdot (Y_{A \tau} - Y_{A\bar\tau})\right)\mod{\pi\I\Z},
\]
where $s=(a,b)$, the curve $y^2 = x^3 +ax+b$ corresponds to $\tau$, $Y_{A \tau}$ and $Y_{A \bar\tau}$ are the graphs of the endomorphisms of this curve which act on the tangent space as $A\tau$ and $A\bar\tau$, $A\tau^2 + B\tau + C$ is the minimal equation of $\tau$ with $A>0$, and $D=B^2 - 4 A C$.
\end{cor*}

As an example we verify
\begin{cor*}
\[
G_2^{\HH/PSL_2(\Z)}\left(\frac{-1+\sqrt{-7}}{2}, \I\right) \;=\;  \frac{8}{\sqrt{7}} \log(8-3\sqrt{7}).
\]
\end{cor*}

\vspace{5mm}
The text is organized in five chapters. Each chapter is provided with a more detailed introduction and reading the introduction is highly recommended for understanding the chapter, especially in the case of Chapter 3. We briefly discuss contents of each chapter here. The first chapter studies various functions on the upper half plane and differential operators. The main result of this chapter is the lifting of the values of the Green function at CM points from real numbers to the elements of 
\[
\C/2\pi\I (D_1 D_2)^{\frac{1-k}2} \Q,
\]
and, related to this, the refined version of the algebraicity conjecture (see Corollary \ref{thm:global_st:6} and Conjecture (2)). Also in this chapter we prove that the $k$-th non-holomorphic derivative of $G_k^{\HH/\Gamma}$ is a meromorphic modular form characterized by certain properties (Theorem \ref{thm:global_study3}) and discuss ways to compute the Green function from this modular form (see the discussion in the end of Section \ref{global_green}). The results of this chapter where known to experts but do not exist in print.

The second chapter contains a definition of the higher Chow groups and a construction of the Abel-Jacobi map as in \cite{GL1} and \cite{GL2}, as well as a proof that the definition of the Abel-Jacobi map is correct. We prove a formula (Theorem \ref{thm:spec_val1}) which in certain cases  relates the value of the Abel-Jacobi map with certain intersection number which takes values in the multiplicative group.

The third chapter provides a way to compute the derivative of the Abel-Jacobi map (by this we mean $\left(\frac{d}{d\tau}\right)^{2k-1} A_\XX (\tau)$) for a family of higher cycles. The answer is expressed as a certain extension of $\D$-modules. The result is more general than we need and works for any families of products of curves. In the case of the power of a family of elliptic curves we obtain an invariant of this extension which for modular triples coincides with $\left(\frac{d}{d\tau}\right)^{2k-1} A_\XX (\tau)$ and is a meromorphic modular form.

In the fourth chapter we study various objects on the Weierstrass family of elliptic curves. We provide representatives for cohomology classes which are necessary for later computations.

The fifth chapter contains the main result, namely, the construction of a family of higher cycles which computes the values of the function $G_2^{\HH/\Gamma}(\tau, \I)$ and the proof of the algebraicity conjecture in this case (Theorem \ref{thm:grfunc:1}). There we show that the meromorphic modular form $(2\pi\I)^2 \sqrt{-4}\,  \frac{432 E_4(\tau)}{j(\tau) - 1728}$, which equals $-\frac{1}{2} \sqrt{-4} \, \delta^2 G_2^{\HH/\Gamma}(\tau, \I)$, is geometrically representable.

The construction of the family of higher cycles used in the proof was inspired by other constructions of higher cycles on products of elliptic curves in \cite{GL1}, \cite{GL2}. It seems, however, that our family (see Section \ref{second-cycle} and the construction of $\XX_{-4}$ above) was not known before, though its construction is surprisingly simple. 

This text is submitted as a PhD thesis to the University of Bonn. The work was done at the Max-Planck-Institute for Mathematics. The author is grateful to the institute for hospitality and good working atmosphere. Also I wish to thank S. Bloch, J. Bruinier, N. Durov, A. Goncharov, G. Harder, D. Huybrechts, C. Kaiser, Yu. Manin, R. Sreekantan with whom I had interesting discussions on the subject of this thesis. Special thanks to M. Vlasenko who read drafts of this text and made useful remarks, and to D. Zagier who introduced me to number theory, proposed the problem, provided me with inspiration and support, and was a very good supervisor.

\tableofcontents

\chapter{Modular forms}
In this chapter we first fix some notation for the representations of $SL_2(\R)$ on the space of functions on the upper half plane. We obtain a representation for each integer $w$ which is called ``weight''. We introduce three differential operators on these representations which are intertwining: the non-holomorphic derivative $\delta$ or the ``raising operator'', the ``lowering operator'' $\delta^-$ and the Laplacian $\Delta$. We also recall the construction of the standard finite-dimensional representations $V_m$ of $SL_2(\R)$ and vectors in these representations which correspond to complex multiplication (CM) points. This is done in Section \ref{notations}.

Next, in Section \ref{eigenvalues}, we study in more detail functions which are eigenfunctions for the Laplacian with eigenvalue of the form $k(1-k)$, $k\in \Z$. The situation with these particular eigenvalues is special. In particular, to each eigenfunction $f$ with such eigenvalue corresponds a certain ``extended function'' $\wt f$ which is harmonic, but takes values in the space $V_{2k-2}$. The derivatives of these extended functions are proportional to holomorphic functions of weight $2k$ which we also call ``derivatives''.

In Section \ref{integrating} we study the inverse problem of recovering $\wt f$ from its derivatives in the case when they are invariant under the action of a congruence subgroup of $PSL_2(\Z)$. We obtain that the obstructions to solving this problem lie in certain cohomology groups.

Then, in Section \ref{local_green}, we define the Green functions $G_k^\HH$ for the upper half plane (without a congruence subgroup), the so-called ``local Green functions''. These are defined as functions of two variables $z_1$, $z_2$ in the upper half plane which are $SL_2(\R)$-invariant (diagonal action) eigenfunctions for the Laplacian with eigenvalue $k(1-k)$ and have only a logarithmic singularity on the diagonal. We show how these functions can be explicitly written using Legendre's functions. We also evaluate the actions of powers of the operator $\delta$ on them. We obtain particularly nice expressions for the $k$-th power which also gives formulae for the derivatives of the corresponding extended Green functions.

Section \ref{global_green} studies the ``global Green functions'' $G_k^{\HH/\Gamma}$ for a quotient of the upper half plane. Since these Green functions can be obtained by averaging the local ones, results of Section \ref{local_green} give us information about the singularities of the global Green functions and their derivatives. We obtain a characterisation of the derivative $g_{k,z_0}^{\HH/\Gamma}$ of the extended global Green function as a meromorphic modular form with described type of poles and having trivial class in a certain Eichler-Shimura cohomology group. Then the original Green function can be obtained from the modular form using the integration procedure described in Section \ref{integrating}.

We suppose that there are no cusp forms of weight $2k$ for the group $\Gamma$. Then being applied to CM points the integration procedure allows us to lift the value of the Green function, which is real, to a complex number defined up to an algebraic multiple of $2\pi\I$ (see Corollary \ref{thm:global_st:6}). This supports and refines the {\em algebraicity conjecture} formulated in papers \cite{GZ} and \cite{GKZ}, which says that in this case the value must be an algebraic multiple of the logarithm of an algebraic number (see Conjectures (1) and (2)). In the end of the section we discuss ways to compute the lifted Green function, and we show that it is a period (Theorem \ref{thm:global_st:7}).

\section{Notations}\label{notations}

We will consider the group $SL_2(\R)$. Elements of this group will be usually denoted by $\gamma$, and matrix elements by $a$, $b$, $c$, $d$:
\[
\gamma \eq \begin{pmatrix} a & b \\ c & d \end{pmatrix} \in SL_2(\R).
\]
This groups acts on the upper half plane $\HH$.

\subsection{Representations of $SL_2(\R)$}
The group $SL_2(\R)$ naturally acts on the following linear spaces:
\begin{itemize}
\item $V=\C^2$~--- the space of column vectors of length 2,
\item $V_m$~--- the symmetric $m$-th power of $V$.
\end{itemize}

We explain the way we view elements of $V_m$. Let $e_1$, $e_2$ be the natural basis on $V$. Then $V_m$ is the space of homogeneous polynomials in $e_1$, $e_2$ of degree $m$. If we substitute $e^1$ by a new variable $X$ and $e^2$ by $1$ we obtain a non-homogeneous polynomial in one variable of degree less or equal $m$. We represent elements of $V_m$ as polynomials in the variable $X$ of degree less or equal $m$, i.e.
\[
p\in V_m, \; p(X) \eq p_0 + p_1 X + \dots + p_m X^m.
\]
The group acts on $V_m$ on the right by
\[
(p \gamma)(X) \eq (p |_{-m} \gamma)(X) \eq p(\gamma X) (cX + d)^m, \qquad \gamma\eq\begin{pmatrix} a & b \\ c & d \end{pmatrix} \in SL_2(\R).
\]
There is the corresponding action on the left
\[
(\gamma p)(X) \eq (p \gamma^{-1})(X) \eq p(\gamma^{-1} X) (-cX + a)^m.
\]

Now consider the dual space to $V_m$, $V_m^*$. Any $v\in V_m^*$ is a functional on $V_m$. Suppose its value on $p$ is
\[
v_0 p_0 + v_1 p_1 + \dots v_m p_m,
\]
then we represent $v$ as a raw vector $(v_0, v_1, v_2, \dots, v_m)$. For any two numbers $x, y \in \C$ we form a vector
\[
v_{x, y} \eq (y^m, y^{m-1} x, y^{m-2} x^2, \dots, x^m) \in V_m.
\]
Then for any $p\in V_m$
\[
(v_{x, y}, p) \eq p\left(\frac{x}{y}\right) y^m =: p(x, y),
\]
hence
\[
\gamma v_{x, y} \eq v_{\gamma (x, y)} \eq v_{a x + b y, c x + d y}.
\]
In this way the action of $SL_2(\R)$ is given on a subset of $V_m^*$ which spans $V_m^*$.

The isomorphism between $V_m$ and $V_m^*$ can be given as follows. Define an invariant pairing between elements of the form $v_{x,y} \in V_m$:
\[
(v_{x,y}, v_{x', y'}) \eq (x y' - x' y)^m \eq \sum_{i=0}^m (-1)^i \binom{m}{i} x^{m-i} y^i x'^i y'^{m-i}.
\]
Since it is a homogeneous polynomial of degree $m$ in each pair of variables, this induces an equivariant linear map from $V_m^*$ to $V_m$
\[
(v_0, v_1, \dots, v_m) \mapsto \sum_{i=0}^m (-1)^i \binom{m}{i} v_{m-i} X^i.
\]
This is an isomorphism of representations. It induces invariant pairings on $V_m$
\[
\left(\sum_{i=0}^m p_i X^i, \sum_{i=0}^m p_i' X^i\right) \eq \sum_{i=0}^m \frac {(-1)^i p_i p_{m-i}'}{\binom{m}{i}}
\]
and on $V_m^*$
\[
((v_i), (v_i')) \eq \sum_{i=0}^m (-1)^i \binom{m}{i} v_i v_{m-i}'.
\]

The pairings above are symmetric for $m$ even and antisymmetric for $m$ odd.
The following identities hold:
\[
((z-X)^m, p) \eq p(z) \eq (p, (X-z)^m)   \qquad            (p\in V_m, z\in\C).
\]

\subsection{Differential operators}\label{sec:diff-op}
Let $S$ be a discrete subset of the upper half plane $\HH$ and $f(\z)$ be a function on $\HH-S$ with values in $\C$ and $w\in \Z$. Define differential operators
\begin{align*}
\delta_w f &\eq \frac{\partial f}{\partial z} \;+\; \frac{w}{\z - \zc} f,\\
\delta_w^- f &\eq (\z - \zc)^2 \frac{\partial f}{\partial \zc},\\
\Delta_w f &\eq (\z - \zc)^2 \frac {\partial}{\partial \z} \frac {\partial}{\partial \zc} f \;+\; w (\z - \zc) \frac {\partial}{\partial \zc}.
\end{align*}
We think about $w$ as the weight, attached to the function $f$. The weight will be always clear from context, so we will omit the subscript $w$. 
We will follow the following agreement: the operator $\delta$ increases weight by $2$, the operator $\delta^-$ decreases weight by $2$ and the operator $\Delta$ leaves weight untouched. Taking into account this agreement the following identities can be proved:
\begin{align*}
\delta^- \delta \;-\; \delta\, \delta^- &\eq w,\\
\delta\, \delta^- &\eq \Delta,\\
\delta^- \delta &\eq \Delta \;+\; w.
\end{align*}

Let the group $SL_2(\R)$ act on functions of weight $w$ by the usual formula:
\[
(f|_w \gamma)(\z) \eq f(\gamma \z) (c\z + d)^{-w}.
\]
This is a right action. We also define the corresponding left action
\[
(\gamma f)(\z) \eq (f|_w \gamma^{-1})(\z) \eq f(\gamma^{-1} \z) (-c\z + a)^{-w}.
\]
Note that this action commutes with the operators $\delta$, $\delta^-$, $\Delta$, it maps functions defined on $\HH-S$ to functions defined on $\HH-\gamma S$.

It is also convenient to modify the complex conjugation for functions with weight to make it commuting with the action of the group. For $f$ of weight $w$ we put
\[
f^*(\z) \eq (\z - \zc)^w \overline{f(\z)}.
\]
Assign to $f^*$ weight $-w$. We can check that
\begin{align*}
f^{**} &\eq (-1)^w f, \\
\delta (f^*) &\eq (\delta^- f)^*, \\
\delta^- (f^*) &\eq (\delta f)^*, \\
\Delta (f^*) &\eq ((\Delta + w) f)^*, \\
\gamma (f^*) &\eq (\gamma f)^*.
\end{align*}
We remark that for the weight $0$ the operator $*$ is the usual complex conjugation, the operator $\delta$ is the usual $\frac{\partial}{\partial z}$, and the operator $\Delta$ is the usual Laplace operator for the hyperbolic metric $-y^2(\frac{\partial^2}{\partial x^2} + \frac{\partial^2}{\partial y^2})$.

We list several formulae, which are convenient to use in computations. We assume that the constant function $1$ has weight $0$. Consider the functions
\[
X-z,\qquad \frac{X-\zc}{z-\zc},
\]
which are thought as functions in $z$ with values in $V$ of weights $-1$ and $1$ respectively. Then
\begin{gather*}
\delta 1 \eq \delta^- 1 \eq 0, \\
(X-z)^* \eq \frac{X-\zc}{z-\zc},\qquad \left(\frac{X-\zc}{z-\zc}\right)^* \eq -(X-z), \\
\delta (X-z) \eq -\frac{X-\zc}{z-\zc},\qquad \delta \frac{X-\zc}{z-\zc} \eq 0, \\
\delta^- (X-z) \eq 0,\qquad \delta^- \frac{X-\zc}{z-\zc} \eq X-z, \\
\delta \left((X-z)^a \left(\frac{X-\zc}{z-\zc}\right)^b\right) \eq -a (X-z)^{a-1} \left(\frac{X-\zc}{z-\zc}\right)^{b+1}, \\
\delta^- \left( (X-z)^a \left(\frac{X-\zc}{z-\zc}\right)^b\right) \eq b (X-z)^{a+1} \left(\frac{X-\zc}{z-\zc}\right)^{b-1}, \\
\Delta \left((X-z)^a \left(\frac{X-\zc}{z-\zc}\right)^b\right) \eq -b (a+1) (X-z)^a \left(\frac{X-\zc}{z-\zc}\right)^b.
\end{gather*}

\subsection{CM points and quadratic forms}
We will frequently use the following notation:
\[
Q_z \eq \frac{(X-z)(X-\zc)}{z-\zc}.
\]
As a function of $z$ $Q_z$ is a function with values in $V_2$ of weight $0$.

A CM point is a point $z\in \HH$ which satisfies a quadratic equation of degree $2$ with integer coefficients. Let $z$ be a CM point. Write the minimal equation for $z$:
\[
a z^2 + b z + c \eq 0,\qquad a,b,c\in \Z,\;\; a>0.
\]
The discriminant of the minimal equation $D_z=b^2-4 a c$ is called the discriminant of $z$. An elementary computation gives
\[
Q_z \eq \frac {a X^2 + b X + c}{\sqrt{D}}.
\]

It is clear that a point $z\in \HH$ is a CM point if and only if $Q_z$ is proportional to a polynomial with integer coefficients.

\section{Eigenfunctions of the Laplacian}\label{eigenvalues}
We consider functions on $\HH-S$ for discrete subsets $S\subset\HH$.
For integers $k, w$ we denote by $F_{k, w}$ the space of functions of weight $w$ satisfying
\[
\Delta f \eq (k (1-k) + \tfrac {w (w-2) }{4}) f.
\]
It is easy to check the following properties of the spaces $F_{k, w}$:
\begin{prop}
\begin{enumerate}
\item The space $F_{k, w}$ is invariant under the action of the group $SL_2(\R)$ (meaning, of course, that $\gamma\in SL_2(\R)$ changes $S$ to $\gamma S$).
\item The operator $*$ maps $F_{k, w}$ to $F_{k, -w}$.
\item The operator $\delta$ maps $F_{k, w}$ to $F_{k, w+2}$. It is invertible for all values of $w$, except, possibly, $2k - 2$ and $-2k$ with the inverse given by
\[
\delta_w^{-1} \eq \frac4{(w+2k)(w-2k+2)} \delta_{w+2}^-.
\]
\item The operator $\delta^-$ maps $F_{k, w}$ to $F_{k, w-2}$. It is invertible for all values of $w$, except, possibly, $2k$ and $2-2k$ with the inverse given by
\[
(\delta_w^-)^{-1} \eq \frac4{(w-2k)(w+2k-2)} \delta_{w-2}.
\]
\end{enumerate}
\end{prop}

We will occasionally use negative powers of $\delta$ and $\delta^-$ when they can be defined using this proposition.

Next we state some basic facts
\begin{prop}\label{prop2_2}
For integer numbers $k$, $l$ we have
\[
(X-z)^{k-l-1} \left(\frac{X-\zc}{z-\zc}\right)^{k+l-1} \in F_{k, 2l}\otimes V_{2k-2}.
\]
\end{prop}
\begin{proof}
Use the formulae listed in the end of Section \ref{sec:diff-op}.
\end{proof}

\begin{prop}
For $f$, $g$ in $F_{k,0}$ and $1-k\leq l \leq k-1$ we have
\[
\delta^{-l} f \, \delta^l g \eq (\delta^-)^l f \, (\delta^-)^{-l}g.
\]
\end{prop}
\begin{proof}
Note that
\[
(\delta^-)^l f \eq (-1)^l \frac{(k+l-1)!}{(k-l-1)!}\; \delta^{-l} f,
\]
and analogously
\[
\delta^l g \eq (-1)^l \frac{(k+l-1)!}{(k-l-1)!}\; (\delta^-)^{-l} g,
\]
which implies the required identity.
\end{proof}

Let us introduce the following operation. For two functions $f$, $g$ from $F_{k,0}$ we put
\[
f * g \eq \sum_{l=1-k}^{k-1} \delta^{-l}\! f \, \delta^l \! g,
\]
which has weight $0$.
Note that the previous proposition implies
\[
f*g \eq \sum_{l=1-k}^{k-1} (-1)^l \; (\delta^-)^{-l} \! f\; (\delta^-)^l \! g,
\]
so
\[
\overline{f*g} \eq \overline{f}*\overline{g}.
\]
It is also easy to see that 
\[
\frac{\partial}{\partial z} (f*g) \eq (-1)^{k-1}(\delta^{1-k}\! f\, \delta^k\! g \;+\; \delta^k\! f \,\delta^{1-k}\! g).
\]
Consider the function $Q_z^{k-1}$. By Proposition \ref{prop2_2} 
\[
    Q_z^{k-1} \in F_{k, 0}\otimes V_{2k-2} \qquad \text{as a function of $z$.}
\]
One can compute that for $1-k \leq l \leq k-1$
\begin{align*}
\delta^l Q_z^{k-1} &\eq (-1)^l\, (k-1)! \frac{(X-z)^{k-1-l} }{(k-l-1)!} \left(\frac{X-\zc}{z-\zc}\right)^{k-1+l} \; (1-k \leq l \leq k-1), \\
\delta^k Q_z^{k-1} &\eq 0.
\end{align*}

For $f\in F_{k, 0}$ we denote
\begin{align*}
\wt f &\eq (-1)^{k-1} \binom{2k-2}{k-1} f * Q_z^{k-1}  
\\ &\eq (-1)^{k-1} \frac{(2k-2)!}{(k-1)!} \sum_{l=1-k}^{k-1} \frac{(X-z)^{k+l-1}}{(k+l-1)!} \left(\frac{X-\zc}{z-\zc}\right)^{k-l-1} \delta^l f.
\end{align*}
It is easy to check that
\[
\wt{\gamma f} \eq \gamma \wt f\qquad \text{for $\gamma\in SL_2(\R)$,}
\]
where $\gamma$ acts on $\wt f$ by the simultaneous action on $z$ in weight $0$ and $X$ in weight $2-2k$. Also
\[
\overline{\wt f} \eq (-1)^{k-1} \wt{\overline f}.
\]
One can compute the scalar product of $\delta^i Q_z^{k-1}$ and $\delta^j Q_z^{k-1}$ as follows:
\begin{lem}
\[
(\delta^i Q_z^{k-1}, \delta^j Q_z^{k-1}) \eq \begin{cases}
0, & \text{if $i\neq -j$} \\
(-1)^{k-1-i} \binom{2k-2}{k-1}^{-1} & \text{if $i=-j$.}
\end{cases}
\]
\end{lem}
\begin{proof}
We prove by induction on $i$ starting from $1-k$.
Since 
\[
\delta^{1-k} Q_z^{k-1} \eq (-1)^{k-1} \frac{(k-1)!}{(2k-2)!}\, (X-\z)^{2k-2},
\]
for any polynomial $p\in V_{2k-2}$ we have
\[
(\delta^{1-k} Q_z^{k-1}, p) \eq (-1)^{k-1} \frac{(k-1)!}{(2k-2)!}\, p(\z).
\]
Hence $(\delta^{1-k} Q_z^{k-1}, \delta^j Q_z^{k-1})$ is not zero only for $j=k-1$ and in this case
\[
(\delta^{1-k} Q_z^{k-1}, \delta^{k-1} Q_z^{k-1}) \eq \binom{2k-2}{k-1}^{-1}.
\]
If the statement is true for $i$ then for any $j$, taking into account that the weight of $(\delta^i Q_z^{k-1}, \delta^j Q_z^{k-1})$ is $2i+2j$,
\[
0 \eq \delta(\delta^i Q_z^{k-1}, \delta^j Q_z^{k-1}) 
\eq (\delta^{i+1} Q_z^{k-1}, \delta^j Q_z^{k-1}) + (\delta^i Q_z^{k-1}, \delta^{j+1} Q_z^{k-1}).
\]
Hence
\[
(\delta^{i+1} Q_z^{k-1}, \delta^j Q_z^{k-1}) \eq -(\delta^i Q_z^{k-1}, \delta^{j+1} Q_z^{k-1}).
\]
We see that if $i+j\neq -1$, this is zero. If $i+j=-1$, this equals exactly
\[
-(-1)^{k-1-i} \binom{2k-2}{k-1}^{-1}.
\]
\end{proof}

Therefore the original function $f$ can be recovered from $\wt f$ as
\[
f \eq (\wt f, Q_z^{k-1}).
\]

Note also that $\Delta \wt f=0$. This is true because
\[
\delta^- \delta^k f \eq (\Delta+2 k - 2)\, \delta^{k-1} f \eq 0.
\]

Let us summarize.
\begin{thm} \label{thm:eigenfunc:5}
Let $f\in F_{k, 0}$ for $k\geq 1$. Then the function $\wt f$ satisfies the following properties (note that $F_{0,0}$ is the space of harmonic functions):
\begin{align*}
\wt f &\;\in\; F_{0,0}\otimes V_{2k-2},\\
f &\eq (\wt f, Q_z^{k-1}),\\
\delta^l f &\eq (\wt f, \delta^l Q_z^{k-1}), \\
\frac{\partial \wt f}{\partial z} &\eq (X-z)^{2k-2}\frac{(-1)^{k-1}\delta^k f}{(k-1)!}\, ,\\
\frac{\partial \wt f}{\partial \zc} &\eq (X-\zc)^{2k-2}\frac{\bar{\delta}^k f}{(k-1)!}
\qquad (\bar\delta g := \ol{\delta \bar g}).
\end{align*}
In the opposite way, if $g\in F_{0,0}\otimes V_{2k-2}$ is such that $\frac{\partial g}{\partial z}$ is divisible by $(X-z)^{2k-2}$ and $\frac{\partial g}{\partial \zc}$ is divisible by $(X-\zc)^{2k-2}$, then the function $f(z) = (g(z), Q_z^{k-1})$ satisfies $f\in F_{k,0}$ and $g=\wt f$.
\end{thm}
\begin{proof}
It is only not clear how to prove the ``opposite way'' part of the statement. Suppose $g$ satisfies the conditions above. Consider the function
\[
g_0 \eq (g, \delta^{k-1} Q_z^{k-1}).
\]
Then 
\[
\delta g_0 \eq (\delta g,\, \delta^{k-1} Q_z^{k-1}) \;+\; (g,\, \delta^{k} Q_z^{k-1})  \eq \left( \frac{\partial}{\partial z} g, \delta^{k-1} Q_z^{k-1} \right).
\]
This shows that $\delta g_0$ is a product of a holomorphic function and the function 
\[
\left( (X-z)^{2k-2},\, \delta^{k-1} Q_z^{k-1} \right),
\]
which is constant. Therefore $\delta g_0$ is itself holomorphic. Hence 
\[
\Delta g_0 \eq (\delta^- \delta - (2k-2)) g_0 = (2-2k) g_0,
\]
which means $g_0\in F_{k,2k-2}$. Next we use that $\delta^- g$ is divisible by $(X-\zc)^{2k-2}$ to prove that
\[
(g, \, \delta^l Q_z^{k-1}) \eq \delta^{l+1-k} g_0 \in F_{k, 2l}.
\]
This shows that $f=\delta^{1-k} g_0 \in F_{k,0}$ and $g=\wt f$.
\end{proof}

We mention one more formula for $\wt f$, this time expressing it as a polynomial in $(X-z)$:
\begin{prop} \label{prop:eigenfunc:6}
Let $f\in F_{k, 0}$ for $k\geq 1$. Then 
\[
\wt f = (-1)^{k-1}\, \frac{(2k-2)!}{(k-1)!} \sum_{i=0}^{2k-2} \frac{(X-z)^i}{i!} \left(\frac{\partial}{\partial z}\right)^i \delta^{1-k} f,
\]
and
\[
\left(\frac{\partial}{\partial z}\right)^{2k-1} \delta^{1-k} f = \delta^k f.
\]
\end{prop}
\begin{proof}
The first formula can be verified in two steps. First, check that both sides coincide when $X=z$. Second, take the derivative of both sides with respect to $z$ and notice that derivative of the left hand side must be divisible by $(X-z)^{2k-2}$.

The second step of the proof also gives the second formula.
\end{proof}

Suppose we know $\delta^k f$, $\bar{\delta}^k f$ and we want to recover $f$. Using the theorem above one can transform the latter problem into the problem of recovering $\wt f$ from $\frac{\partial \wt f}{\partial z}$ and $\frac{\partial \wt f}{\partial \zc}$. The next section explains how this can be done.

\section{Integrating modular forms}\label{integrating}
Let $\Gamma$ be a congruence subgroup of $SL_2(\R)$, $S$~--- a finite union of orbits of $\Gamma$ in $\HH$, $U=\HH-S$. Consider a smooth closed differential $1$-form $\omega$ on $U$ with coefficients in $V_{2k-2}$, where $V_{2k-2}$ is the space of polynomials in one variable $X$ of degree not greater than $2k-2$. $V_{2k-2}$ is equipped with the natural action of $SL_2(\R)$ (even $SL_2(\C)$). Suppose moreover that $\omega$ is equivariant for the action of $\Gamma$ and let $A$ be a $\Gamma$-submodule of $V_{2k-2}$. We are looking for a function 
\[
I_\omega^{A,\Gamma}: U\longrightarrow V_{2k-2}/A,
\]
which satisfies the following properties:
\begin{enumerate}
\item $I_\omega^{A,\Gamma}$ is smooth, which means that for any point $z\in U$ there is a neighborhood $W$ of $z$ and a smooth function $g:W\longrightarrow V_{2k-2}$, such that $g$ coincides with $I_\omega^{A,\Gamma}$ modulo $A$.
\item $d I_\omega^{A,\Gamma} = \omega$,
\item $I_\omega^{A,\Gamma}$ is equivariant for the action of $\Gamma$.
\end{enumerate}
Our basic example of $\omega$ is 
\[
\omega \eq f(z) (X-z)^{2k-2} dz,
\]
for some meromorphic modular form $f$ of weight $2k$. The module $A$ may be the submodule of $V_{2k-2}$ consisting either of polynomials with real coefficients, polynomials with imaginary coefficients or of polynomials with coefficients in some discrete subgroup of $\C$.

It is clear that for existence of $I_\omega^{A,\Gamma}$ the following condition is necessary:
\begin{condition}[Residue condition]
For every point $s\in S$ the integral of $\omega$ along a small loop around $s$ belongs to $A$.
\end{condition}

Choose a basepoint $a \in U$. The differential $\omega$ defines a cocycle with values in $V^{2k-2}/A$ in the following way:
\[
\sigma_a(\gamma) \eq \int_{a}^{\gamma a} \omega.
\]
It satisfies 
\[
\sigma_a(\gamma \gamma') \eq \sigma_a(\gamma) \;+\; \gamma \sigma_a(\gamma').
\]
The integral of $\omega$ is a function with values in $V_{2k-2}/A$:
\[
I_a(x) \eq \int_{a}^x \omega.
\]
If we change the basepoint from $a$ to $a'$ the integral and the cocycle change in the following way:
\[
I_{a'}(x) \eq I_a(x) \;+\; \int_{a'}^a \omega,
\]
\[
\sigma_{a'}(\gamma) \eq \sigma_a(\gamma) \;+\; \int_{a'}^a \omega \;-\; \int_{\gamma a'}^{\gamma a} \omega,
\]
so if we introduce an element $v_{a a'}$ of $V_{2k-2}/A$ by 
\[
v_{a a'} \eq \int_{a'}^a \omega,
\]
then
\[
I_{a'}(x) \eq I_a(x) \;+\; v_{a a'},
\]
\[
\sigma_{a'} \eq \sigma_a \;-\; \delta v_{a a'},
\]
where $\delta$ denotes the differential 
\[
(\delta v)(\gamma) \eq \gamma v \,-\, v.
\]
Therefore $\sigma_a$ defines a class in $\sigma_\omega^{A,\Gamma}\in H^1(\Gamma, V_{2k-2}/A)$, which does not depend on the base point.

If such an $I_\omega^{A,\Gamma}$ as explained in the beginning of the section exists, it must differ from $I_a$ by a constant, say $v_a\in V_{2k-2}/A$. Let us describe all $v_a\in V_{2k-2}$ such that $I_a+v_a$ is equivariant for the action of $\Gamma$. This means that for any $\gamma\in\Gamma$
\[
I_a(\gamma x) \,+\, v_a \eq \gamma\, (I_a(x) \,+\, v_a).
\]
Splitting the path of integration and using the equivariance of $\omega$ we conclude:
\[
I_a(\gamma x) \eq I_{\gamma a}(\gamma x) \,+\, \sigma_a(\gamma) \eq \gamma I_a(x) \,+\, \sigma_a(\gamma),
\]
so the last equation is equivalent to the following:
\[
\sigma_a(\gamma) \eq \gamma v_a - v_a \eq \delta v_a.
\]
Therefore
\begin{prop}
Suppose $\omega$ satisfies the residue condition. Then the function $I_\omega^{A,\Gamma}$ satisfying the properties (i)-(iii) exists if and only if the class $\sigma_\omega^{A,\Gamma}\in H^1(\Gamma, V_{2k-2}/A)$ is trivial. If such a function exists, then it is unique up to addition of an element of $H^0(\Gamma,V_{2k-2}/A)$ and satisfies
\[
\gamma I_\omega^{A,\Gamma}(z) \,-\, I_\omega^{A,\Gamma}(z) \eq \int_z^{\gamma z} \omega \mod A \qquad \text{for all $\gamma\in\Gamma$.}
\]
\end{prop}

For any $\Gamma$-module $M$ we denote by $C^1(\Gamma, M)$ the abelian group of cocycles, i.e. maps $\sigma:\Gamma \longrightarrow M$ such that
\[
\sigma(\gamma_1\gamma_2) \eq \sigma(\gamma_1) \,+\, \gamma_1\sigma(\gamma_2), \qquad \text{for all $\gamma_1, \gamma_2 \in \Gamma$.}
\]
Then we have the following exact sequence:
\[
0\longrightarrow H^0(\Gamma, M) \To M \xrightarrow{\;\,d\;\,} C^1(\Gamma, M) \To H^1(\Gamma, M)\to 0,
\]
where $d$ is given by the usual formula
\[
(d m)(\gamma) \eq \gamma m \,-\, m, \qquad m\in M,\gamma\in\Gamma
\]
Denote by $C_1(\Gamma, M)$ the abelian group of cycles, by which we mean the quotient group of $\Z\Gamma\otimes M$ by the subgroup generated by elements of the form
\[
\gamma_1\gamma_2\otimes m \,-\, \gamma_2\otimes m \,-\, \gamma_1\otimes\gamma_2 m, \qquad \text{for $m\in M$, $\gamma_1, \gamma_2 \in \Gamma$.}
\]
We have the corresponding exact sequence for homology:
\[
0\longrightarrow H_1(\Gamma, M) \To C_1(\Gamma, M) \xrightarrow{\;\,d\;\,} M  \To H_0(\Gamma, M)\to 0,
\]
where $d$ is given by 
\[
d (\gamma\otimes m) \eq \gamma m - m, \qquad m\in M,\gamma\in\Gamma.
\]
Suppose we have an invariant biadditive pairing
\[
(\cdot,\cdot): M\otimes M' \To N,
\]
where $M$ and $M'$ are $\Gamma$-modules and $N$ is an abelian group.
Denote by $C^1(\Gamma, M)_0$ the group of cocycles which map to zero in $H^1(\Gamma, M)$, and by $M'_0$ the group of elements in $M'$ which map to zero in $H_0(\Gamma, M')$.

Note that there is a canonical pairing 
\[
C^1(\Gamma, M)\otimes C_1(\Gamma, M') \To N,
\]
given by
\[
(\sigma, \gamma\otimes m) \eq (\sigma(\gamma^{-1}), m),
\]
since
\begin{multline*}
(\sigma,\, \gamma_1\gamma_2\otimes m \,-\, \gamma_2\otimes m \,-\, \gamma_1\otimes \gamma_2 m) \\ \eq
(\sigma(\gamma_2^{-1}\gamma_1^{-1}), m) \,-\, (\sigma(\gamma_2^{-1}), m) \,-\, (\sigma(\gamma_1^{-1}), \gamma_2 m)\eq 0.
\end{multline*}

This pairing has the following property which can be easily checked:
\begin{prop}
For any $m\in M$, $c\in C_1(\Gamma, M')$
\[
(dm, c) \eq (m, dc).
\]
\end{prop}

This implies that $(\cdot, \cdot)$ induces pairings
\[
H^0(\Gamma, M)\otimes H_0(\Gamma, M') \To N, \qquad 
H^1(\Gamma, M)\otimes H_1(\Gamma, M') \To N.
\]

If we have $\sigma\in C^1(\Gamma, M)_0$ and $m'\in M'_0$, we can either represent $\sigma$ as a coboundary, i.e. $\sigma = d m_0$, and then consider $(m_0, m')$, or represent $m'$ as a boundary, $m' = d c$, and then consider $(\sigma, c)$. 
\begin{prop}
Either of these two approaches defines a pairing 
\[
C^1(\Gamma, M)_0 \otimes M'_0 \To N,
\]
moreover the resulting two pairings coincide.
\end{prop}

We put $M=V_{2k-2}/A$, $M' = B$, where $B$ is a $\Gamma$- invariant subgroup of $V_{2k-2}$ and $N=\C / (A, B)$. 

\begin{thm}\label{int_pairing}
Let $S\subset \HH$ be a finite union of orbits of $\Gamma$, $A$ and $B$ be $\Gamma$-invariant subgroups of $V_{2k-2}$, $\omega$ be a smooth closed invariant differential $1$-form on $U=\HH-S$ with coefficients in $V_{2k-2}$ whose integrals along small loops around points of $S$ belong to $A$, $z$ be a point in $U$, $v$ be an element in $B$ such that:
\begin{enumerate}
\item The class of $\omega$ in $H^1(\Gamma, V_{2k-2}/A)$ is $0$,
\item the class of $v$ in $H_0(\Gamma, B)$ is $0$.
\end{enumerate}
Then the following two approaches lead to the same element of $\C/(A, B)$, which we denote by $I^{A, B, \Gamma}(\omega, z, v)$:
\begin{enumerate}
\item First represent $\omega$ as a differential of an invariant $V_{2k-2}/A$-valued function $I^{A, \Gamma}_\omega$, and then put
\[
I^{A, B, \Gamma}(\omega, z, v) \eq (I_\omega^{A, \Gamma}(z), v).
\]
\item First represent $v$ as
\[
v \eq \sum_{i=1}^n (\gamma_i u_i - u_i)\qquad \text{for $\gamma_i\in\Gamma$, $u_i\in B$,}
\]
and then put
\[
I^{A, B, \Gamma}(\omega, z, v) \eq \sum_{i=1}^n \left(\int_z^{\gamma_i^{-1} z} \omega,\; u_i \right).
\]
\end{enumerate}
\end{thm}


\section{Local study of Green's functions}\label{local_green}

Let $k$ be an integer, $k>1$.
We consider the Green function for the upper half plane of weight $2 k$, which we denote by $G_{k}^{\HH}$ (note that this is not yet the Green function which is our main object of study, the property (ii) will be different). It is the unique function which satisfies the following properties:

\begin{enumerate}
\item $G_k^\HH$ is a smooth function on $\HH\times\HH-\{\z_1=\z_2\}$ with values in $\R$.
\item $G_k^\HH(\gamma \z_1, \gamma \z_2) = G_k^\HH(\z_1, \z_2)$ for all $\gamma\in SL_2(\R)$.
\item $\Delta_i G_k^\HH = k(1-k) G_k^\HH$, where $\Delta_i$ denotes the Laplace operator with respect to $\z_i$.
\item $G_k^\HH = \log|\z_1-\z_2|^2 + O(1)$ when $\z_1$ tends to $\z_2$.
\item $G_k^\HH$ tends to $0$ when $\z_1$ tends to infinity.
\end{enumerate}

In this section we obtain two formulae. The first formula is for $\delta_1^n\delta_2^m G_k^\HH$, and it involves hypergeometric series. The second formula is a particular case of the first for $n=k$, and the resulting expression for this case is a rational function of $\z_1, \z_2, \zc_2$. Note that because of the symmetry between $\z_1$ and $\z_2$ this case is similar to the case $m=k$. 

Because of the second property the function $G_k^\HH$ is a function of the hyperbolic distance. Denote by $t(\z_1, \z_2)$ the hyperbolic cosine of the hyperbolic distance, i.e.
\[
\begin{split}
t(\z_1, \z_2) \eq 1 + 2\frac{(\z_1-\z_2)(\zc_2-\zc_1)}{(\z_1-\zc_1)(\z_2-\zc_2)} \eq -1 + 2\frac{(\z_1-\zc_2)(\z_2-\zc_1)}{(\z_1-\zc_1)(\z_2-\zc_2)}\\
\eq\frac{-2\z_1\zc_1 + (\z_1+\zc_1)(\z_2+\zc_2) - 2\z_2\zc_2}{(\z_1-\zc_1)(\z_2-\zc_2)}.
\end{split}
\]
Then 
\[
G_k^\HH(\z_1, \z_2) \eq -2 \calQ_{k-1}(t),
\]
where $\calQ_{k-1}$ is the Legendre's function of the second kind. The function $\calQ_{k-1}$ has the following integral definition (see \cite{WW}):
\[
\calQ_{k-1}(t) = 2^{-k} \int_{-1}^1 (1-x^2)^{k-1} (t-x)^{-k} dx.
\]

The function $\calQ_{k-1}$ has the following two expansions at infinity, which can be obtained from the integral representation above expanding the integrand into power series (the first expansion can also be found in \cite{WW}):
\[
\begin{split}
\calQ_{k-1}(t) \eq \frac{2^{k-1} (k-1)!^2}{(2k-1)!}\, t^{-k} F\left(\frac{k}2, \frac{k+1}2; k+\frac12; t^{-2}\right) \\
\eq \frac{2^{k-1} (k-1)!^2}{(2k-1)!} (t+1)^{-k}\, F\left(k, k; 2k; \frac2{1+t}\right),
\end{split}
\]
here $F$ denotes the hypergeometric series.
We are going to compute various derivatives of $G_k^\HH$ using the second expansion. For this purpose we first compute:
\[
\delta_1 t \eq \frac{\partial t}{\partial\z_1} \eq 2\frac{(\zc_1-\z_2)(\zc_1-\zc_2)}{(\z_1-\zc_1)^2(\z_2-\zc_2)},
\]
\[
\delta_2 t \eq \frac{\partial t}{\partial\z_2} \eq 2\frac{(\z_1-\zc_2)(\zc_1-\zc_2)}{(\z_1-\zc_1)(\z_2-\zc_2)^2},
\]
noting that $\delta_1 t$ has weight $2$ in $\z_1$ and weight $0$ in $\z_2$ we compute:
\[
\delta_1^2 t \eq \delta_2^2 t \eq 0,
\]
\[
\delta_1 \delta_2 t \eq -2\left(\frac{\zc_1-\zc_2}{(\z_1-\zc_1)(\z_2-\zc_2)}\right)^2 \eq \frac{\delta_1 t\, \delta_2 t}{t+1}.
\]
We will use the following formula for the derivative of the hypergeometric series:
\[
\frac{\partial F(a,b;c;x)}{\partial x} \eq a \frac{F(a+1,b;c;x)-F(a,b;c;x)}{x}.
\]
We find that
\begin{equation*}
\frac{\partial ((t+1)^{-m}\, F(m, n; c; \frac2{t+1}))}{\partial t} \eq 
-m (t+1)^{-m-1}\, F(m + 1, n; c; \tfrac 2 {t+1}), 
\end{equation*}
so
\begin{multline*}
\delta_1^n G_k^\HH(\z_1, \z_2) \eq (-1)^{n+1}\, 2^k \frac{(k-1)! (k+n-1)!}{(2k-1)!} \times \\ 
(t+1)^{-k-n}\, F(k+n,k;2k;\tfrac 2{t+1})\, (\delta_1 t)^n.
\end{multline*}
To apply $\delta_2^m$ we rewrite the last expression as
\begin{multline*}
\delta_1^n G_k^\HH(\z_1, \z_2) \eq (-1)^{n+1}\, 2^k \frac{(k-1)! (k+n-1)!}{(2k-1)!} \times \\
(t+1)^{-k}\, F(k,k+n;2k;\tfrac 2 {t+1})\, (\delta_1\delta_2 t)^n (\delta_2 t)^{-n},
\end{multline*}
so we can again apply the same formula, since $\delta_2$ of $\delta_1\delta_2 t$ and $\delta_2 t$ is zero:
\begin{multline*}
\delta_2^m\delta_1^n G_k^\HH(\z_1, \z_2) \eq (-1)^{m+n+1}\, 2^k \frac{(k+m-1)!(k+n-1)!}{(2k-1)!} \times \\
(t+1)^{-k-m}\, F(k+m,k+n;2k;\tfrac 2 {t+1})\, (\delta_1\delta_2 t)^n (\delta_2 t)^{m-n}.
\end{multline*}
To make the formula symmetric in $m$ and $n$ we rewrite it as
\begin{multline*}
\delta_2^m\delta_1^n G_k^\HH(\z_1, \z_2) \eq (-1)^{m+n+1}\, 2^k \frac{(k+m-1)!(k+n-1)!}{(2k-1)!} \times \\
(t+1)^{-k-m-n}\, F(k+m,k+n;2k;\tfrac 2 {t+1})\, (\delta_1 t)^n (\delta_2 t)^m.
\end{multline*}
Let us introduce the following function of weight $-2$ in $\z_1$ and $0$ in $\z_2$:
\[
Q_{\z_2}(\z_1) \eq \frac{(\z_1 - \z_2)(\z_1 - \zc_2)}{\z_2-\zc_2},
\]
there is a corresponding function $Q_{\z_1}(\z_2)$. One can check:
\[
\delta_1 t \eq \frac{t^2 - 1}2\, Q_{\z_2}(\z_1)^{-1},
\]
\[
\delta_2 t \eq \frac{t^2 - 1}2\, Q_{\z_1}(\z_2)^{-1},
\]
so our first formula is
\begin{multline*}
\delta_2^m\delta_1^n G_k^\HH(\z_1, \z_2) \eq (-1)^{m+n+1}\, \frac{(k+m-1)!(k+n-1)!}{(2k-1)!} \left(\frac{t+1}2\right)^{-k}\! \times \\
 \left(\frac{t-1}2\right)^{m+n} \! F(k+m,k+n;2k;\tfrac 2 {t+1})\; Q_{\z_2}(\z_1)^{-n}\, Q_{\z_1}(\z_2)^{-m} \qquad (m, n \geq 1-k).
\end{multline*}

In particular, when $n=k$ we obtain
\begin{multline*}
\delta_2^m\delta_1^k G_k^\HH(\z_1, \z_2) \eq (-1)^{m+k+1}\, 2^{-m} (k+m-1)! (t+1)^{-k} (t-1)^{k+m} \times\\
F(k+m, 2k; 2k; \tfrac 2 {t+1})\; Q_{\z_2}(\z_1)^{-k}\, Q_{\z_1}(\z_2)^{-m},
\end{multline*}
and using the identity
\[
F(a, b; b; x) \eq (1-x)^{-a}
\]
we get the second formula
\begin{align*}
\delta_1^k \delta_2^m G_k^\HH(\z_1, \z_2) &\eq (-1)^{m+k+1}\, (k+m-1)! \left(\frac{t+1}2\right)^m\!  Q_{\z_1}(\z_2)^{-m} Q_{\z_2}(\z_1)^{-k}\\
&\eq (-1)^{k-1}\, (k+m-1)!\, \frac{(\z_2-\zc_2)^{k-m}}{(\z_1-\z_2)^{k+m}(\z_1-\zc_2)^{k-m}}\\
&\eq \frac{(-1)^{k-1}\, (k+m-1)!}{(z_1-z_2)^{2m}}\; Q_{z_2}(z_1)^{m-k} \qquad (m\geq 1-k).
\end{align*}

\section{Global study of Green's functions}\label{global_green}
Let $\Gamma$ be a congruence subgroup of $PSL_2(\Z)$ and $k>1$. The Green's function on $\HH/\Gamma$ of weight $2k$ is the unique function $G_k^{\HH/\Gamma}$ with the following properties:

\begin{enumerate}
\item $G_k^{\HH/\Gamma}$ is a smooth function on $\HH\times\HH-\{\z_1=\gamma\z_2\ \,|\, \gamma\in\Gamma\}$ with values in $\R$.
\item $G_k^{\HH/\Gamma}(\gamma_1 \z_1, \gamma_2 \z_2) = G_k^{\HH/\Gamma}(\z_1, \z_2)$ for all $\gamma_1, \gamma_2\in \Gamma$.
\item $\Delta_i G_k^{\HH/\Gamma} = k(1-k) G_k^{\HH/\Gamma}$.
\item $G_k^{\HH/\Gamma} = m\log|\z_1-\z_2|^2 + O(1)$ when $\z_1$ tends to $\z_2$ ($m$ is the order of the stabilizer of $z_2$, which is almost always $1$).
\item $G_k^{\HH/\Gamma}$ tends to $0$ when $\z_1$ tends to a cusp.
\end{enumerate}

The series
\[
\sum_{\gamma\in\Gamma} G_k^\HH(z_1, \gamma z_2)
\]
is convergent and satisfies the properties above, so
\[
G_k^{\HH/\Gamma}(z_1, z_2) \eq \sum_{\gamma\in\Gamma} G_k^\HH(z_1, \gamma z_2).
\]

Consider the function $G_k^{\HH/\Gamma}(z, z_0)$ for a fixed $z_0\in\HH$. We put 
\begin{align*}
\G_{k, z_0}^{\HH/\Gamma}(z) &\eq \wt{G_k^{\HH/\Gamma}(z, z_0)} \\
&\eq (-1)^{k-1} \binom{2k-2}{k-1}\sum_{l=1-k}^{k-1} (-1)^l\, \delta^{-l}(Q_z^{k-1}) \, \delta^l G_k^{\HH/\Gamma}(z, z_0) \\
&\eq (-1)^{k-1} \frac{(2k-2)!}{(k-1)!} \sum_{l=1-k}^{k-1} \frac{(X-z)^{k+l-1}}{(k+l-1)!} \left(\frac{X-\zc}{z-\zc}\right)^{k-l-1} \delta^l G_k^{\HH/\Gamma}(z, z_0).
\end{align*}
recall that negative powers of $\delta$ can be defined for eigenfunctions of the Laplacian. 

Note that since $G_k^{\HH/\Gamma}(z, z_0)$ has real values and $Q_z$ has imaginary values (we let complex conjugation act on $X$ identically), we have the following
\begin{prop}
The values of the function $i^{k-1} \G_{k, z_0}^{\HH/\Gamma}(z)$ are real polynomials.
\end{prop}
Since 
\[
\delta^k(Q_z^{k-1}) \eq 0,
\]
it is easy to compute that
\begin{align*}
\frac{\partial \G_{k, z_0}^{\HH/\Gamma}(z)}{\partial z} \eq \delta \G_{k, z_0}^{\HH/\Gamma}(z) 
&\eq \binom{2k-2}{k-1}\, \delta^{1-k}(Q_z^{k-1})\, \delta^k G_k^{\HH/\Gamma}(z, z_0) \\
&\eq (X-z)^{2k-2}\, \frac{(-1)^{k-1}\delta^k\, G_k^{\HH/\Gamma}(z, z_0)}{(k-1)!}.
\end{align*}
On the other hand, because of the proposition above
\[
\frac{\partial \G_{k, z_0}^{\HH/\Gamma}(z)}{\partial \zc} \eq (X-\zc)^{2k-2} \frac{\overline{\delta^k G_k^{\HH/\Gamma}(z, z_0)}}{(k-1)!}.
\]

Consider the function 
\[
g_{k, z_0}^{\HH/\Gamma}(z) \eq \delta^k G_k^{\HH/\Gamma}(z, z_0).
\]
This is a meromorphic modular form in $z$. There is a corresponding differential $1$-form with coefficients in $V_{2k-2}$
\[
(X-z)^{2k-2} g_{k, z_0}^{\HH/\Gamma}(z) dz.
\]
Let $V_{2k-2}^\R$ denote the space of polynomials in $V_{2k-2}$ which have real coefficients.
\begin{prop}\label{cohProp}
The class of the differential form
\[
(X-z)^{2k-2} g_{k, z_0}^{\HH/\Gamma}(z) dz
\]
in the cohomology group
\[
\I^{k-1} H^1(\Gamma, V_{2k-2}^\R) \eq H^1(\Gamma, V_{2k-2}/\I^k V_{2k-2}^\R)
\]
is trivial and the function
\[
(-1)^{k-1} (k-1)!\, \frac{\G_{k, z_0}^{\HH/\Gamma}(z)}2
\]
is an integral of $\omega$.
\end{prop}
\begin{proof}
Let us denote
\[
\omega \eq (X-z)^{2k-2} g_{k, z_0}^{\HH/\Gamma}(z) dz.
\]
We have proved before that
\[
(-1)^{k-1} (k-1)!\; d \G_{k, z_0}^{\HH/\Gamma}(z, X) \eq \omega + (-1)^{k-1} \bar\omega,
\]
so that
\[
(-1)^{k-1} (k-1)!\; \I^{k-1} d \G_{k, z_0}^{\HH/\Gamma}(z, X) \eq \I^{k-1}\omega + \overline{\I^{k-1}\omega}.
\]
It implies that the integral of $i^{k-1}\omega$ around a pole of $\omega$ is in $\I V_{2k-2}^\R$, so the integral of $\omega$ around a pole is in $\I^k V_{2k-2}^\R$, this is why $\omega$ satisfies the residue condition of Section \ref{integrating} for $A=\I^k V_{2k-2}^\R$. Hence the class of $\omega$ in $H^1(\Gamma, V_{2k-2}/\I^k V_{2k-2}^\R)$ is correctly defined. Moreover, 
\[
(-1)^{k-1} (k-1)!\;\frac12 d \G_{k, z_0}^{\HH/\Gamma}(z, X) \eq \omega - \frac{(-1)^{k} \bar\omega + \omega}2,
\]
so
\[
(-1)^{k-1} (k-1)!\; \frac{\G_{k, z_0}^{\HH/\Gamma}(a, X) - \G_{k, z_0}^{\HH/\Gamma}(b, X)}2 \equiv
\int_{b}^a \omega \mod{i^k V_{2k-2}^\R},
\]
which implies that the class of $\omega$ is trivial:
\[
\sigma_\omega^{i^k V_{2k-2}^\R, \Gamma} \equiv 0 \mod{i^k V_{2k-2}^\R},
\]
and $(-1)^{k-1} (k-1)! \,\frac12 \G_{k, z_0}^{\HH/\Gamma}(z)$ is an integral of $\omega$.
\end{proof}

\begin{thm}\label{thm:global_study3}
For any $z_0\in\HH$ the function $g_{k, z_0}^{\HH/\Gamma}(z)$ is the unique function which satisfies the following properties:
\begin{enumerate}
\item It is a meromorphic modular form of weight $2k$ in $z$ whose set of poles is $\Gamma z_0$ which is zero at the cusps.
\item In a neighborhood of $z_0$
\[
g_{k, z_0}^{\HH/\Gamma}(z) \eq m (-1)^{k-1} (k-1)!\; Q_{z_0}(z)^{-k} + O(1)\qquad (m=|\Stab_\Gamma z_0|)
\]
\item The class of the corresponding differential form in the cohomology group
\[
H^1(\Gamma, V_{2k-2}/\I^{k} V_{2k-2}^\R) \eq \I^{k-1} H^1(\Gamma, V_{2k-2}^\R)
\]
is trivial as in Proposition \ref{cohProp}.
\end{enumerate}
\end{thm}
\begin{proof}
First we prove that the function $g_{k, z_0}^{\HH/\Gamma}(z)$ actually satisfies these conditions. The first two conditions follow from the local study. In fact, we have
\[
\delta^k G_k^{\HH/\Gamma}(z, z_0) \eq \sum_{\gamma\in\Gamma} \delta^k G_k^{\HH}(z, \gamma z_0),
\]
because for all $0\le l \le k$
\[
\delta^l G_k^{\HH/\Gamma}(z, z_0) \eq O(t(z, z_0)^{l-k}\, |Q_{z_0}(z)|^{-l}),
\]
and using the inequality
\[
|Q_{z_0}(z)| \;\geq\; \frac{t(z, z_0)-1}2\, |z-\zc|
\]
we obtain that the series
\[
\sum_{\gamma\in\Gamma} \delta^l G_k^{\HH}(z, \gamma z_0)
\]
is locally uniformly majorated by the series
\[
\sum_{\gamma\in\Gamma} t(z, \gamma z_0)^{-k},
\]
which converges locally uniformly in $z$.

This already implies that the function $g_{k, z_0}^{\HH/\Gamma}(z)$ is meromorphic and has poles of the specified type. The transformation property follows from the invariance of $G_k^{\HH/\Gamma}(z, z_0)$. Since we can move any cusp to $\infty$ by an element of $SL_2(\Z)$ it is enough to check the cuspidality at $\infty$. This is clear from the following (recall that $k>1$):
\begin{multline*}
\delta^k G_k^{\HH/\Gamma}(z, z_0) \eq\\
(-1)^{k-1} (k-1)! (z_0-\zc_0)^k \sum_{\gamma\in\Gamma} \frac1{(\gamma z- z_0)^k (\gamma z - \zc_0)^k (cz+d)^{2k}}.
\end{multline*}

The last condition is precisely Proposition \ref{cohProp}. It only remains to prove the uniqueness.

Now suppose there are two different functions which satisfy the conditions. Then their difference is a cusp form of weight $2k$ which has either purely real or purely imaginary cohomology depending on whether $k$ is odd or even. This contradicts the Eichler-Shimura theorem, which says that
\[
H_{parabolic}^1(\Gamma, V_{2k-2}) \cong S_{2k} \oplus \overline{S_{2k}},
\]
so any nontrivial parabolic cohomology class which is either purely real or purely imaginary cannot be represented by a cusp form.
\end{proof}

The value of the Green function can be recovered in the following way:
\begin{thm}
Let $k>1$, $z_0\in\HH$ and $g_{k, z_0}^{\HH/\Gamma}(z)$ be the function that satisfies the conditions of the theorem above. Put
\[
w \eq (X-z)^{2k-2} g_{k, z_0}^{\HH/\Gamma}(z) dz
\]
and apply one of the two approaches formulated in Theorem \ref{int_pairing} for $A = \I^k V_{2k-2}^\R$, $B = \I^{k-1} V_{2k-2}^\R$, $v=Q_z^{k-1}$ to get an element 
\[
I^{\I^k V_{2k-2}^\R, \I^{k-1} V_{2k-2}^\R, \Gamma}(\omega, z, Q_z^{k-1}) \in \C/\I\R\eq\R,
\]
for some $z\in \HH$.
Then
\[
\frac12 G_k^{\HH/\Gamma}(z, z_0) \;=\; \frac{(-1)^{k-1}}{(k-1)!}\; \Re I^{\I^k V_{2k-2}^\R, \I^{k-1} V_{2k-2}^\R, \Gamma}(\omega, z, Q_z^{k-1}).
\]
\end{thm}
\begin{proof}
We note that Theorem \ref{int_pairing} can be applied since
\begin{enumerate}
\item The class of $\omega$ in $H^1(\Gamma, V_{2k-2}/A)$ is trivial by the third property of the function $g_{k, z_0}^{\HH/\Gamma}(z)$.
\item The whole homology group $H_0(\Gamma, B)$ is trivial because $k>1$.
\end{enumerate}

We can put
\[
I^{A, \Gamma}(z) \eq (-1)^{k-1} (k-1)! \; \frac12 \G_{k,z_0}^{\HH/\Gamma}(z)
\]
because of Proposition \ref{cohProp}. This implies that
\begin{multline*}
I^{i^k V_{2k-2}^\R, i^{k-1} V_{2k-2}^\R, \Gamma}(\omega, z, Q_z^{k-1}) \\ \equiv (-1)^{k-1}(k-1)!\,\left(\frac12 \G_{k,z_0}^{\HH/\Gamma}(z), Q_z^{k-1}\right) \mod \C/\I\R,
\end{multline*}
so the statement follows from the identity
\[
(\G_{k,z_0}^{\HH/\Gamma}(z), Q_z^{k-1}) \eq G_k^{\HH/\Gamma}(z, z_0),
\]
which was proved in Section \ref{eigenvalues}.
\end{proof}

Let us use Theorem \ref{thm:global_study3} to compute the residues of the form
\[
\omega\eq(X-z)^{2k-2} g_{k, z_0}^{\HH/\Gamma}(z) dz.
\]
The residue at $z=z_0$ can be computed as follows:
\begin{multline*}
\res (X-z)^{2k-2} Q_{z_0}(z)^{-k} dz \eq \res \frac{(X-z)^{2k-2} (z_0-\zc_0)^k}{(z-z_0)^k (z-\zc_0)^k} dz\\
\eq \res \frac{(X-z_0 - (z-z_0))^{2k-2}}{(z-z_0)^k(1 + \frac{z-z_0}{z_0-\zc_0})^k} dz,
\end{multline*}
so we have to compute the residue of the following Laurent series
\begin{multline*}
\sum_{i,j} (-1)^{i+j} \binom{2k-2}{i} \binom{k+j-1}{j} (X-z_0)^i (z-z_0)^{k-2-i+j} (z_0-\zc_0)^{-j} dz.
\end{multline*}
This gives (recall that $m$ is the order of the stabilizer of $z_0$.
\begin{multline*}
(-1)^{k-1} \sum_j \binom{2k-2}{k+j-1} \binom{k+j-1}{j} (X-z_0)^{k+j-1} (z_0-\zc_0)^{-j} \\
\eq(-1)^{k-1} \binom{2k-2}{k-1} \sum_j \binom{k-1}{j} (X-z_0)^{k+j-1} (z_0-\zc_0)^{-j} \\
\eq (-1)^{k-1} \binom{2k-2}{k-1} \frac{(X-z_0)^{k-1} (X-\zc_0)^{k-1}}{(z_0-\zc_0)^{k-1}} \eq (-1)^{k-1} \binom{2k-2}{k-1}  Q_{z_0}^{k-1}.
\end{multline*}

Therefore we obtain
\[
\res \omega \eq m \,\frac{(2k-2)!}{(k-1)!}\; Q_{z_0}^{k-1}.
\]
This suggests that if $z_0$ is a CM point we should put
\[
A\eq 2\pi\I \,\frac{(2k-2)!}{(k-1)!}\; D_0^{\frac{1-k}2} V_{2k-2}^\Z,
\]
where $D_0$ is the discriminant of $z_0$. Note that $A\subset \I^k V_{2k-2}^\R$ and the integrals of $\omega$ along small loops around $z_0$ are still in $A$.

Suppose $z$ is also a CM point. Then we may refine $B$ in the following way:
\[
B\eq D^{\frac{1-k}2} V_{2k-2}^\Z.
\]
We still have $Q_{z}^{k-1}\in B$.

Consider the groups $H_0(\Gamma, B)$ and $H^1(\Gamma, V_{2k-2}/A)$. Since $H_0(\Gamma, V_{2k-2})=0$ and $B$ is finitely generated over $\Z$ and spans $V_{2k-2}$ over $\C$, $H_0(\Gamma, B)$ is a finite group.

The form $\omega$ defines certain class in $H^1(\Gamma, V_{2k-2}/A)$, which is cuspidal. We have a direct sum decomposition $V_{2k-2} = \I^k V_{2k-2}^\R \oplus \I^{k-1} V_{2k-2}^\R$ and $A\subset  \I^k V_{2k-2}^\R$. Therefore $V_{2k-2}/A = \I^k V_{2k-2}^\R/A \oplus \I^{k-1}V_{2k-2}^\R$. By Proposition \ref{cohProp} the second component of the class of $\omega$ is trivial, but there is no reason, in general, for the class in $\I^k V_{2k-2}^\R/A$ to be trivial. However the following proposition says

\begin{prop}
Suppose there are no cusp forms for the group $\Gamma$ of weight $2k$. Then the group of classes in $\I^k V_{2k-2}^\R/A$ which are cuspidal is a finite group.
\end{prop}
\begin{proof}
The $\Gamma$-module $\I^k V_{2k-2}^\R/A$ can be obtained from the $\Gamma$-module $V_{2k-2}^\Z$ by extension of scalars. Since the corresponding cuspidal cohomology group of $V_{2k-2}^\Z$ is finite, the statement follows.
\end{proof}

Suppose there are no cusp forms for the group $\Gamma$ of weight $2k$. Let $N_A$ be the exponent of the finite group in the proposition above. Let $N_B$ be the exponent of $H_0(\Gamma, B)$. Then applying Theorem \ref{int_pairing} to $N_A \omega$ and $N_B Q_z^{k-1}$ we obtain an element
\[
I^{A, B, \Gamma}(N_A \omega, z, N_B Q_z^{k-1}) \in \C/(A,B).
\]

Note that
\[
(A,B)\subset 2\pi\I (D_0 D)^{\frac{1-k}2} m \,\frac{(2k-2)!}{(k-1)!}\; \langle V_{2k-2}^\Z, V_{2k-2}^\Z\rangle.
\]
It is easy to see that
\[
\langle V_{2k-2}^\Z, V_{2k-2}^\Z\rangle \subset \frac{(k-1)!}{(2k-2)!} \Z.
\]
Therefore
\[
(A,B)\subset 2\pi\I (D_0 D)^{\frac{1-k}2} \Z.
\]

Let us denote $N=N_A N_B (k-1)!$. Then we may define
\[
\widehat{G}_k^{\HH/\Gamma}(z, z_0)  \eq \frac{(-1)^{k-1}}{N}\, I^{A, B, \Gamma}(N_A \omega, z, N_B Q_z^{k-1}) \in \C/2\pi\I (D_0 D)^{\frac{1-k}2} \frac{1}{N} \Z.
\]
This construction proves the following

\begin{cor}\label{thm:global_st:6}
Let $\Gamma$ be a congruence subgroup of $PSL_2(\Z)$, $k>1$ an integer such that there are no cusp forms of weight $2k$ for $\Gamma$. Let $N_A$ be the exponent of the group $H^1_{parabolic}(\Gamma, V^{\Q/\Z}_{2k-2})$, and $N_B$ be the exponent of the group $H_0(\Gamma, V^\Z_{2k-2})$, $N=N_A N_B (k-1)!$. Then for any two non-equivalent CM points $z_1$, $z_2$ of discriminants $D_1$, $D_2$ the pair $N_A (X-z)^{2k-2} g_{k, z_1}^{\HH/\Gamma} dz$, $N_B Q_{z_2}^{k-1}$ satisfies the conditions of Theorem \ref{int_pairing} for $A=2\pi\I\frac{(2k-2)!}{(k-1)!} D_1^{\frac{1-k}2} V_{2k-2}^\Z$, $B=D_2^{\frac{1-k}2}V_{2k-2}^\Z$. The corresponding value
\[
\widehat{G}_k^{\HH/\Gamma}(z_1, z_2) \,=\, \frac{(-1)^{k-1}}{N}\, I^{A,B,\Gamma}(\,N_A (X-z)^{2k-2} g_{k, z_1}^{\HH/\Gamma}(z) dz,\; z_2,\; N_B Q_{z_2}^{k-1})
\]
is defined modulo $2\pi\I \frac{1}{N} \Z$ and its real part equals $\frac12 G_k^{\HH/\Gamma}(z_1, z_2)$.
\end{cor}

The algebraicity conjecture was formulated in \cite{GZ} (p. 317) and extended in \cite{GKZ} (p. 556). In the case when there are no cusp forms it says the following:
\begin{conjecture}[1]
Suppose there are no cusp forms of weight $2k$ for $\Gamma$. Then for any two CM points $z_1$, $z_2$ of discriminants $D_1$, $D_2$ there is an algebraic number $\alpha$ such that
\[
G_k^{\HH/\Gamma}(z_1, z_2) \eq (D_1 D_2)^{\frac{1-k}2} \log \alpha.
\]
\end{conjecture}

We can refine the conjecture now:

\begin{conjecture}[2]
Suppose there are no cusp forms of weight $2k$ for $\Gamma$. Then for any two CM points $z_1$, $z_2$ of discriminants $D_1$, $D_2$ there is an algebraic number $\alpha$ such that
\[
(D_1 D_2)^{\frac{k-1}2} \widehat{G}_k^{\HH/\Gamma}(z_1, z_2) \eq  \log \alpha \mod \frac{2\pi\I}{N} \Z.
\]
\end{conjecture}

The second version of the conjecture is stronger since it gives not only the absolute value of $\alpha$, but also its argument.

\begin{rem}
The conjecture as it is formulated in \cite{GZ} and \cite{GKZ} is more general since it also deals with the case when there are cusp forms of weight $2k$ for $\Gamma$. It is not difficult to extend Corollary \ref{thm:global_st:6} and Conjecture (2) to the more general case.
\end{rem}

Both ways to define $I^{A,B,\Gamma}$ (see Theorem \ref{int_pairing}) are useful. The first way, which was also explained in the introduction, shows that for a given CM point $z_0$ for any invariant $V_{2k-2}/\frac{1}{N_A} A$~-valued function $g$ such that
\[
d g = (X-z)^{2k-2} g_{k, z_0}^{\HH/\Gamma}(z) dz,
\]
one has
\[
\widehat{G}_k^{\HH/\Gamma}(z, z_0) \;\equiv\; \frac{(-1)^{k-1}}{(k-1)!} \;(g, Q_z^{k-1})\qquad \mod 2\pi\I(D_0 D)^{\frac{1-k}2} \frac{1}{N} \Z.
\]
This gives us a way to prove the conjecture. Namely, the idea is to construct some function $g$ such that $g$ satisfies the condition above and $(g, Q_z^{k-1})$ can be computed algebraically. Later we will indeed construct such $g$ for the point $z_0=\I$ using Abel-Jacobi maps for higher Chow groups.

The second way is to represent $Q_z^{k-1}$ as a sum
\[
Q_z^{k-1} \eq \sum_{i=1}^n (\gamma_i u_i - u_i) \qquad \text{for $\gamma_i\in\Gamma$, $u_i\in \frac{1}{N_B} B$,}
\]
and then put
\begin{multline*}
\widehat{G}_k^{\HH/\Gamma}(z, z_0)  \;\equiv\; \frac{(-1)^{k-1}}{(k-1)!} \; \sum_{i=1}^n \left( \int_z^{\gamma_i^{-1} z} (X-z)^{2k-2} g_{k, z_0}^{\HH/\Gamma}(z) dz, u_i \right) \\
 \mod 2\pi\I(D_0 D)^{\frac{1-k}2} \frac{1}{N} \Z.
\end{multline*}
This approach provides a way to compute $\widehat{G}_k^{\HH/\Gamma}(z, z_0)$ approximately by numerical integration and proves the following important fact:

\begin{thm}\label{thm:global_st:7}
Let $\Gamma$ be a congruence subgroup of $PSL_2(\Z)$, $k>1$ an integer such that there are no cusp forms of weight $2k$ for $\Gamma$. Then for any two non-equivalent CM points $z_1$, $z_2$ the value $\widehat{G}_k^{\HH/\Gamma}(z_1, z_2)$ is a {\em period} (see \cite{KZ}).
\end{thm}


\chapter{Higher Chow groups and Abel-Jacobi maps}
In this chapter we assume $X$ is a smooth projective variety over $\C$ and $k$ is an integer. We give definition of the first higher Chow group $CH^k(X,1)$ as the group of cycles modulo the group of boundaries. The cycles are formal sums 
\[
\sum_i (W_i, f_i)
\]
where each $W_i$ is a subvariety of codimension $k-1$ in $X$ and each $f_i$ is a non-zero rational function on $W_i$ such that the sum of divisors of $f_i$ is $0$ as a cycle of codimension $k$. The boundaries are cycles which can be constructed from pairs of functions $(f, g)$ as $(\Div f, g) - (\Div g, f)$ (see Section \ref{chow_groups} for details).

In Section \ref{abel-jacobi} we give a construction of the Abel-Jacobi map from the first higher Chow group to the quotient of the cohomology group $H^{2k-2}(X, \C)$ by a certain subspace and a certain lattice. The subspace is the $k$-th step in the Hodge filtration and the lattice is the integral cohomology. We prove that the definition does not depend on various choices in the remaining part of Section \ref{abel-jacobi}.

It is clear that it makes sense to consider the pairing of the value of the Abel-Jacobi map with the cohomology class of an ordinary cycle $Z$ (linear combination of irreducible subvarieties) of the ``complementary dimension'' (dimension is $k-1$). This gives a complex number defined up to an integer. In Section \ref{spec_val_abel-jacobi}, Theorem \ref{thm:spec_val1}, we prove that the resulting number times $2\pi\I$ is the logarithm of the intersection number, which can be computed explicitly as the product over the intersection points of $Z$ and $W_i$ of values of $f_i$ over all $i$ (taking multiplicities into account). Although this statement follows from compatibility of Abel-Jacobi maps with intersection products, we give two independent proofs. One uses currents, and another uses relative cohomology and a construction of logarithmic representatives of fundamental classes in relative cohomology.

\section{Notation}
The symbol ``$\I$'' denotes $\sqrt{-1}$ to distinguish it from ``$i$'', which will usually be used as an index.

\section{The Hodge theory}
Let $X$ be a smooth projective variety over $\C$. For each $k$ we denote the sheaf of smooth $k$-forms by $\A^k_X$. We have the usual decomposition
\[
\A^k_X\eq\bigoplus_{p+q=k}\A^{p, q}_X,
\]
where $\A^{p, q}_X$ is the sheaf of smooth $(p,q)$-forms on $X$. We have the Hodge filtration on $\A^k_X$ defined as
\[
F^j\A^k_X \;:=\; \bigoplus_{p+q=k,\, p\ge j}\A^{p, q}_X.
\]
One can compute the cohomology groups of $X$ by taking the cohomology of the complex of global forms:
\[
H^k(X, \C) \eq H^k(\A^\bullet(X)),
\]
where for a sheaf $\F$ and an open set $U$ $\F(U)$ denotes the sections of the sheaf over $U$ and we omit the subscript $X$ when we write $\A^\bullet(X)$. One obtains the Hodge filtration on $H^k(X, \C)$ as the one induced by the filtration on $\A^k_X$, i.e.
\[
F^jH^k(X,\C) \eq \frac{\kernel(d:\A^k(X)\To\A^{k+1}(X))\cap F^j\A^k(X)}{\image(d:\A^{k-1}(X)\To\A^k(X))\cap F^j\A^k(X)}.
\]
By the Hodge theory there is a canonical decomposition
\[
H^k(X,\C)\eq\bigoplus_{p+q=k} H^{p, q}(X,\C)
\]
with
\[
F^jH^k(X,\C)\eq\bigoplus_{p+q=k,\,p\ge j} H^{p, q}(X,\C).
\]
We will frequently use the following consequence of the Hodge theory:
\begin{prop}
If $\omega\in F^j\A^k(X)$ is exact, i.e. there exists $\eta\in\A^{k-1}(X)$ with $d\eta=\omega$, then $\omega=d\eta'$ for some $\eta'\in F^j\A^{k-1}(X)$.
\end{prop}
\begin{proof}
Let 
\[
\omega\eq\sum_{p+q=k,\,p\ge j} \omega_{p, q},\qquad \eta\eq\sum_{p+q=k-1}\eta_{p, q},\qquad\omega_{p, q},\eta_{p, q}\in\A^{p, q}(X).
\]
Put
\[
\eta^\#\eq\sum_{p+q=k-1,\,p\ge j} \eta_{p, q}.
\]
Then the following form is exact:
\[
\omega^\#\eq\omega-d\eta^\#\eq\omega_{j, k-j}-\bar\partial\eta_{j, k-j-1}.
\]
Since $d\omega=0$ and $\omega\in F^j\A^k(X)$ we have $\bar\partial\omega_{j, k-j}=0$. Therefore $\bar\partial\omega^\#=0$, so $\omega^\#$ gives some class in the Dalbeaut cohomology group $H^{j, k-j}_{\bar\partial}(X) \cong H^{j,k-j}(X, \C)$. This class is trivial because of the Hodge decomposition and the fact that the class of $\omega^\#$ is trivial in $H^k(X,\C)$. Hence $\omega^\#=\bar\partial\eta'$ for some $\eta'\in\A^{p,q-1}(X)$. This implies that $\omega-d(\eta^\#+\eta')\in F^{j+1}\A^k(X)$, so we can use induction on $j$ to complete the proof.
\end{proof}
Note that there is a canonical homomorphism
\[
\iota:H^k(X,\Z)\To H^k(X,\C).
\]
To simplify the notation we will sometimes write
\[
H^k(X,\Z)\cap F^j H^k(X,\C),\qquad H^k(X,\Z)+F^j H^k(X,\C),
\]
respectively, instead of
\[
H^k(X,\Z)\cap \iota^{-1}(F^j H^k(X,\C)),\qquad H^k(X,\Z)+\iota(F^j H^k(X,\C)).
\]

\section{Higher Chow groups}\label{chow_groups}
Let $X$ be a smooth complex projective variety of dimension $n$. Recall that the ordinary Chow group $CH^k(X)$ of codimension $k$ cycles is the quotient group
\[
CH^k(X) \;:=\; Z^k(X)/B^k(X),
\]
where $Z^k(X)$ is the free abelian group generated by irreducible algebraic subvarieties of $X$ of codimension $k$ and $B^k(X)$ is the subgroup generated by principal divisors on subvarieties of $X$ of codimension $k-1$.

To define the first higher Chow group $CH^k(X,1)$ (see \cite{GL2}) consider the group $C^k(X,1)$ which is the free abelian group generated by pairs $(W,f)$, where $W$ is an irreducible algebraic subvariety of $X$ of codimension $k-1$ and $f$ is a non-zero rational function on $W$, modulo the relations
\[
(W, f_1 f_2) \eq (W,f_1)\,+\,(W, f_2),
\]
where $f_1$ and $f_2$ are two rational functions on $W$. The group $Z^k(X,1)$ is defined to be the kernel of the map 
\[
C^k(X,1)\To B^k(X)
\]
sending $(W,f)$ to $\Div f$, the divisor of $f$. Define the group $B^k(X,1)$ as the subgroup of $Z^k(X,1)$ generated by elements of the form
\[
(\Div g, h|_{\Div g}) \;-\; (\Div h, g|_{\Div h}),
\]
where $g$, $h$ are non-zero rational functions on some $V\subset X$ of codimension $k-2$ whose divisors have no component in common. Here we extend the notation $(W,f)$ to linear combinations of subvarieties, i.e. if $W=\sum_j n_j W_j$ is a linear combination of irreducible subvarieties with integer coefficients and $f$ is a non-zero rational function on some bigger subvariety which restricts to a non-zero rational function on each $W_j$, then 
\[
(W,f) \;:=\; \sum_j (W_j,f^{n_j}|_{W_j}).
\]
We put
\[
CH^k(X,1) \;:=\; Z^k(X,1)/B^k(X,1)
\]
so that any element of $CH^k(X,1)$ has the form
\[
\sum_i (W_i, f_i),
\]
where 
\[
\sum_i \Div f_i \eq 0.
\]
\begin{example}
If $k=1$ then the corresponding higher Chow group is simply the multiplicative group of complex numbers, $CH^1(X,1)=\C^\times$.
\end{example}
\begin{example}
The group $CH^{n+1}(X,1)$ is generated by pairs $(W,f)$ where $W$ is a point and $f$ is a non-zero complex number. In fact one can see that there is a surjective homomorphism of abelian groups
\[
CH^n(X)\otimes \C^\times \To CH^{n+1}(X,1).
\]
Using the Weil reciprocity law one can also check that there exists a homomorphism from $CH^{n+1}(X,1)$ to $\C^\times$ which sends $(W,f)$ to $f$.
\end{example}

\section{The Abel-Jacobi map}\label{abel-jacobi}
\subsection{Abel-Jacobi for the ordinary Chow group}
Recall that for the ordinary Chow group we have the cycle class map
\[
\cl^k:CH^k(X) \To H^{2k}(X, \Z)\cap F^k H^{2k}(X,\C),
\]
which sends $V$, a subvariety of codimension $k$, to its class 
\[
[V] \;\in\; H_{2n-2k}(X,\Z) \;\cong\; H^{2k}(X,\Z).
\]
Denote by $CH^k(X)_0$ the kernel of $\cl^k$. We have the Abel-Jacobi map
\[
AJ^k:CH^k(X)_0 \To \frac{H^{2k-1}(X, \C)}{F^k H^{2k-1}(X,\C)+H^{2k-1}(X,\Z)}.
\]
This is defined as follows. Let $\gamma$ be an algebraic cycle of codimension $k$ whose homology class is $0$. It follows that $\gamma=\partial \xi$ for some $2n-2k+1$-chain $\xi$. Choosing such $\xi$ we obtain a linear functional on the space of $2n-2k+1$-forms given by integrating a form against $\xi$. 

Let us show that this defines a linear functional on $F^{n-k+1}H^{2n-2k+1}$. Indeed, if $\omega$ is an exact form from $F^{n-k+1}\A^{2n-2k+1}$ then $\omega=d\eta$ for $\eta\in F^{n-k+1}\A^{2n-2k}$ by Hodge theory. Hence
\[
\int_\xi\omega \eq \int_\gamma\eta
\]
is zero. 

Choosing another $\xi'$ such that $\partial \xi'=\gamma$ we have $\xi-\xi'$ closed, so the corresponding functionals for $\xi$ and $\xi'$ differ by the functional induced by the corresponding element of $H_{2n-2k+1}(X, \Z)$. So we obtain a map
\begin{multline*}
AJ^k:CH^k(X)_0 \To \frac{(F^{n-k+1}H^{2n-2k+1}(X,\C))^*}{H_{2n-2k+1}(X, \Z)} \\
\cong \frac{H^{2k-1}(X, \C)}{F^k H^{2k-1}(X,\C)+H^{2k-1}(X,\Z)}.
\end{multline*}
\subsection{Abel-Jacobi for the first higher Chow group}
Let $x$ represent an element of $CH^k(X,1)$, $k\in Z^i(X,1)$ i.e.
\[
x\eq\sum_i (W_i, f_i)
\]
with
\[
\sum_i \Div f_i \eq 0.
\]
We denote the corresponding element in $CH^k(X,1)$ by $[x]$.
We choose a path $[0,\infty]\subset \C \PP^1$. Let us denote 
\[
\gamma_i \eq f_i^* [0,\infty],
\]
which is a $2n-2k+1$-chain on $X$ whose boundary is $-\Div f_i$. This implies that the chain
\[
\gamma \eq \sum_i \gamma_i
\]
is a cycle. By the Poincar\'e duality $\gamma$ has a class $[\gamma]\in H^{2k-1}(X,\Z)$. 
\begin{prop}
The map $x\rightarrow [\gamma]$ defines a cycle class map
\[
\cl^{k,1}:CH^k(X,1)\To H^{2k-1}(X, \Z) \cap F^k H^{2k-1}(X, \C).
\]
\end{prop}
\begin{proof}
We have to check two things:
\begin{enumerate}
\item For any $x\in Z^k(X,1)$ the image of $[\gamma]$ belongs to $F^k H^{2k-1}(X, \C)$.
\item If $x=(\Div g, h|_{\Div g})-(\Div h, g|_{\Div h})$
for $V$, $g$ and $h$ as in the definition of $B^k(X,1)$ in Section \ref{chow_groups}, then $\gamma$ is homologically trivial.
\end{enumerate}

For (i) it is enough to prove that the pairing of $\gamma$ with any element of $F^{n-k+1} H^{2n-2k+1}(X,\C)$ is zero. Take a closed form $\omega\in F^{n-k+1} \A^{2n-2k+1}(X)$. Recall that
\[
\int_\gamma \omega \eq \sum_i \int_{\gamma_i} \omega.
\]
Let $n(\gamma_i)$ be a small tubular neighborhood of $\gamma_i$ in $W_i$. Then, up to terms which tend to $0$ as the radius of the neighborhood tends to $0$, we have
\begin{align*}
\int_{\gamma_i} \omega 
\eq \frac{1}{2\pi \I}\int_{\partial n(\gamma_i)} \omega \log f_i 
&\eq -\frac{1}{2\pi \I}\int_{W_i-n(\gamma_i)} d(\omega \log f_i) \\
&\eq -\frac{1}{2\pi \I} \int_{W_i-n(\gamma_i)} \frac{d f_i}{f_i} \wedge \omega.
\end{align*}
The form in the last integral belongs to $F^{n-k+2}\A^{2n-2k+2}$ and $W_i$ has complex dimension $n-k+1$, so the integral is zero.

For proving (ii), if $x=(\Div g, h|_{\Div g})-(\Div h, g|_{\Div h})$, then the chain $-(g\times h)^* ([0,\infty]\times[0,\infty])$ has boundary $\gamma$.
\end{proof}

On the other hand since complex conjugation acts trivially on the group $H^{2k-1}(X, \Z)$ and
\[
F^k H^{2k-1}(X,\C) \cap \overline{F^k H^{2k-1}(X,\C)} \eq \{0\},
\]
we have the following
\begin{prop}
For any $[x]\in CH^k(X,1)$ the class $\cl^{k,1} [x]$ is torsion.
\end{prop}
Thus the cycle class map is a map
\[
\cl^{k,1}:CH^k(X,1)\To H^{2k-1}(X, \Z)_{tors}.
\]
Let us denote the kernel of this map by $CH^k(X,1)_0$.
To construct the Abel-Jacobi map 
\[
AJ^{k,1}:CH^k(X,1)_0\To \frac{H^{2k-2}(X,\C)}{F^k H^{2k-2}(X,\C)+2\pi\I\, H^{2k-2}(X,\Z)}
\]
we first identify
\[
\frac{H^{2k-2}(X,\C)}{F^k H^{2k-2}(X,\C)+2\pi\I H^{2k-2}(X,\Z)}\cong
\frac{(F^{n-k+1}H^{2n-2k+2}(X,\C))^*}{2\pi\I\, H^{2k-2}(X,\Z)}.
\]
If $[x]\in CH^k(X,1)_0$ with $x\in Z^k(X,1)$ then $\gamma=\partial \xi$ for some $2n-2k+2$-chain $\xi$. Then for any $\omega\in F^{n-k+1}\A^{2n-2k+2}(X)$ with $d\omega=0$ we take the following number:
\[\tag{*}
\langle AJ^{k,1}[x], [\omega] \rangle\eq \sum_i\int_{W_i\setminus\gamma_i} \omega\log f_i \,+\, 2\pi\I\,\int_\xi \omega,
\]
where the logarithm on $\C \PP^1$ is defined using the cut along the chosen path $[0,\infty]$.

To prove that this correctly defines a map 
\[
AJ^{k,1}:CH^k(X,1)_0\To \frac{F^{n-k+1}H^{2n-2k+2}(X,\C)^*}{2\pi\I\, H^{2k-2}(X,\Z)}
\]
we need to show that the construction does not depend on the following choices:
\begin{itemize}
\item the choice of the path $[0,\infty]$;
\item the choice of the branch of the logarithm on $\C \PP^1-[0,\infty]$;
\item the choice of the representative of $x$, which is defined up to an element of $B^k(X,1)$;
\item the choice of $\xi$, which is defined up to a $2n-2k+2$-cycle;
\item the choice of $\omega$, which is defined up to a coboundary.
\end{itemize}
We prove this in a series of propositions.
\begin{prop}
The value of (*) does not depend on the choice of the path $[0,\infty]$.
\end{prop}
\begin{proof}
Let $p$ and $p'$ be two different paths on $\C \PP^1$ from $0$ to $\infty$. Let 
\begin{align*}
\gamma_i &\eq f_i^* p, & \gamma_i' &\eq f_i^*p', \\
\gamma &\eq \sum_i \gamma_i, & \gamma' &\eq\sum_i \gamma_i'.
\end{align*}
Choose a 2-chain $q$ on $\C \PP^1$ whose boundary is $p'-p$, let
\[
\eta_i \eq f_i^* q, \qquad \eta\eq\sum_i \eta_i.
\]
We choose $\xi$ such that $\partial \xi = \gamma$ and put $\xi'=\xi+\eta$ so that $\partial \xi' = \gamma'$. Let $l$ be a branch of the logarithm on $\C \PP^1-p$. Then the function
\[
l'(t)\eq l(t)-2\pi\I{\mathbf 1}_q(t)
\]
is a branch of the logarithm on $\C \PP^1-p'$, where ${\mathbf 1}_q$ is the characteristic function of $q$. Then we compare
\begin{align*}
\sum_i\int_{W_i-\gamma_i'} \omega \,l'(f_i) -
\sum_i\int_{W_i-\gamma_i} \omega \,l(f_i) &\eq -2 \pi\I \int_\eta \omega, \\
\int_{\xi'}\omega - \int_{\xi}\omega &\eq \int_\eta\omega,
\end{align*}
so the value of (*) does not change.
\end{proof}
\begin{rem}
In fact this proof also shows that we can even vary each $\gamma_k$ as long as its class in the homology $H_{2n-2k+1}(W_i,|\Div f_i|)$ remains constant.
\end{rem}
\begin{prop}
Changing the branch of the logarithm changes the value of (*) by an element from $2\pi\I\,H^{2k-2}(X,\Z)$.
\end{prop}
\begin{proof}
Changing the branch of the logarithm amounts to adding $2\pi\I m$ for $m\in \Z$, which changes the value of (*) by
\[
2\pi\I m \sum_i \int_{W_i}\omega,
\]
which is a functional induced by the image of $m W_i$ in $H^{2k-2}(X,\Z)$.
\end{proof}
\begin{prop}
If
\[
x\eq\Div g \otimes h|_{\Div g} \,-\, \Div h \otimes g|_{\Div h},
\]
for $V$, $g$ and $h$ as in the definition of $B^k(X,1)$ in Section \ref{chow_groups}, then the value of (*) is zero.
\end{prop}
\begin{proof}
We may take
\[
\xi \eq -(g\times h)^* ([0,\infty]\times[0,\infty])
\]
since
\[
\partial ([0,\infty]\times[0,\infty]) \eq [0,\infty]\times([0]-[\infty]) \;-\; ([0]-[\infty])\times [0,\infty].
\]
Inside $V$ we consider $\Div g$ which has real codimension $2$ and $\gamma_h=h^*[0,\infty]$ which has real codimension $1$.
Consider a small neighborhood of the divisor of $g$ which we denote by $n(\Div g)$ and a small neighborhood of $\gamma_g$ which we denote by $n(\gamma_g)$. We will use corresponding notation for $h$, i.e. $n(\Div h)$ and $n(\gamma_h)$. Then we may rewrite
\[
\int_{\Div g-\gamma_h} \omega\log h
\eq\frac{1}{2\pi\I}\int_{\partial n(\Div g)\cap n(\gamma_h)^c} \frac{dg}{g}\wedge\omega \log h.
\]
Take a chain $S=V\cap n(\Div g)^c \cap n(\gamma_h)^c$. Then 
\[
\partial S \eq -\partial n(\Div g)\cap n(\gamma_h)^c \;-\; \partial n(\gamma_h)\cap n(\Div g)^c,
\]
so using the Stokes formula we obtain
\[
-\int_{\partial n(\Div g)\cap n(\gamma_h)^c} \frac{dg}{g}\wedge\omega \log h - \int_{\partial n(\gamma_h)\cap n(\Div g)^c} \frac{dg}{g}\wedge\omega \log h \eq 
\int_{S} \frac{dh}{h}\wedge\frac{dg}{g}\wedge\omega. 
\]
Since $\omega\in F^{n-k+1}\A^{2n-2k+2}(X)$, $\frac{dh}{h}\wedge\frac{dg}{g}\wedge\omega\in F^{n-k+3}\A^{2n-2k+4}(S)$ and its integral is zero because $S\subset V$, which is a $n-k+2$-dimensional complex variety.
Hence
\begin{equation*}
\int_{\Div g-\gamma_h} \omega\log h \eq -\frac{1}{2\pi\I} \int_{\partial n(\gamma_h)\cap n(\Div g)^c} \frac{dg}{g}\wedge\omega \log h
\eq-\int_{\gamma_h\cap n(\Div g)^c} \frac{dg}{g}\wedge\omega.
\end{equation*}
We apply the Stokes formula again for $T=\gamma_h\cap n(\gamma_g)^c$ and $\omega\log g$. We have
\[
\partial T \eq -\Div h\cap n(\gamma_g)^c - (\partial n(\gamma_g)\cap\gamma_h),
\]
so
\[
-\int_{\Div h\cap n(\gamma_g)^c} \omega\log g-\int_{\partial n(\gamma_g)\cap\gamma_h} \omega\log g \eq \int_{\gamma_h-n(\Div g)} \frac{dg}{g}\wedge\omega.
\]
This implies that
\[
\int_{\Div g-\gamma_h} \omega\log h-\int_{\Div h-\gamma_g}\omega\log g \eq \int_{\partial n(\gamma_g)\cap\gamma_h} \omega\log g \eq 2\pi\I\int_{\gamma_g\cap\gamma_h} \omega.
\]
Since $\xi=-\gamma_g\cap\gamma_h$ the statement follows.
\end{proof}
\begin{prop}
Changing $\xi$ by a $2n-2k+2$-cycle changes the value of the functional (*) by an element of $2\pi\I\,H^{2k-2}(X,\Z)$.
\end{prop}
\begin{proof}
This is clear.
\end{proof}
\begin{prop}
Changing $\omega$ by a coboundary does not change the value of the functional (*).
\end{prop}
\begin{proof}
Indeed, if $\omega$ is a coboundary, then $\omega=d \eta$ with $\eta\in F^{n-k+1}\A(X)$ by the Hodge theory. This implies
\[
(*)\eq\sum_i\int_{W_i-\gamma_i} d(\eta \log f_i) \,-\, \sum_i\int_{W_i-\gamma_i}\frac{df_i}{f_i}\wedge\eta \,+\, 2\pi\I\int_\gamma\eta.
\]
Since $\frac{df_i}{f_i}\wedge\eta\in F^{n-i+2}\A(X)$ and $W_i$ is a $n-k+1$-dimensional complex variety, the second summand is zero. Applying the Stokes formula to the first summand we obtain
\[
\sum_i\int_{W_i-\gamma_i} d(\eta \log f_i)\eq-\sum_i\int_{\partial n(\gamma_i)} \eta \log f_i\eq-2\pi\I\sum_i\int_{\gamma_i}\eta,
\]
where $n(\gamma_i)$ denotes a small neighborhood of $\gamma_i$ inside $W_i$. Hence the statement.
\end{proof}

\section{Special values of the Abel-Jacobi map}\label{spec_val_abel-jacobi}
Let $x=\sum_i (W_i, f_i)$ represent an element $[x]\in CH^k(X,1)_0$.
Then 
\[
AJ^{k,1} [x]\in \frac{F^{n-k+1}H^{2n-2k+2}(X,\C)^*}{2\pi\I\,H^{2k-2}(X,\Z)}.
\]
Given a subvariety $Z\subset X$ of codimension $n-k+1$ we may consider $\cl^{n-k+1} Z\in H^{2n-2k+2}(X,\Z)\cap F^{n-k+1} H^{2n-2k+2}(X,\C)$. Then
\[
(AJ^{k,1}[x], \cl^{n-k+1} Z) \;\in\; \C/2\pi\I\,\Z,
\]
so this number can be written as $\log\alpha$ for a unique $\alpha\in\C$. We now show how to construct this number in a different way.
Consider the cycle $S\subset X$ which is the union of all $|\Div f_i|$ and singular parts of $W_i$.
We say that $Z$ intersects $x$ properly if $Z$ properly intersects $S$ and all $W_i$. This means that $Z$ does not intersect $S$ and intersects each $W_i$ in several points. 
Note that by the Moving Lemma for any given $Z$ there exists $Z'$ which is rationally (hence homologically) equivalent to $Z$ and intersects $x$ properly. One defines the intersection number
\[
x\cdot Z = \prod_i \prod_{p\in W_i\cap Z} f_i(p)^{\ord_p(W_i\cdot Z)}.
\]
\begin{thm}\label{thm:spec_val1}
Let $x=\sum_i (W_i,f_i)$ be a representative of $[x]\in CH^k(X,1)_0$.
Let $Z\subset X$ be a smooth subvariety of dimension $k-1$ intersecting $x$ properly.
Then
\[
(AJ^{k,1}[x], \cl^{n-k+1} Z) \;\equiv\; \log(x\cdot Z) \mod 2\pi\I.
\]
\end{thm}
\begin{proof}
By the definition (*) of $AJ^{k,1}$ we have
\[
(AJ^{k,1}[x], \cl^{n-k+1} Z) \eq \sum_i\int_{W_i\setminus\gamma_i}\omega\log{f_i} \,+\, 2\pi\I\int_\xi\omega,
\]
where $\omega\in F^{n-k+1}\A^{2n-2k+2}(X)$ is a form whose class $[\omega]$ equals to the Poincar\'e dual of the class $[Z]\in H_{2k-2}(X,\Z)$.

The current
\[
\omega \;-\; \delta_Z.
\]
is homologically trivial since both $\omega$ and $\delta_Z$ represent the same class in the cohomology. Hence there exists a current $\eta\in F^{n-k+1}\D^{2n-2k+1}(X)$, smooth outside $|Z|$, such that
\[
d\eta \eq \omega \,-\, \delta_Z.
\]
If we denote by $\eta_0$ the corresponding form on $X \setminus |Z|$, 
\[
\eta_0 \;\in\;  F^{n-k+1}\A^{2n-2k+1}(X\setminus |Z|),
\]
we obtain an identity
\[
d\eta_0 \eq \omega,
\]
which is true outside $|Z|$. 

In the definition of the Abel-Jacobi map we choose $\gamma_k$ and $\xi$ to be transversal to $|Z|$. This means that $\gamma_k$ does not intersect $|Z|$ for each $k$ and $\xi$ intersects $|Z|$ only in several points.

Choosing small neighborhoods of $|\Div f_i|$, $|Z|\cap W_i$, $\gamma_i$ inside $W_i$ and $|Z|\cap \xi$ inside $\xi$ and denoting them by $n(\Div f_i)$, $n(|Z|\cap W_i)$, $n(\gamma_i)$, $n(|Z|\cap \xi)$ respectively we may write
\begin{multline*}
\sum_i\int_{W_i\setminus(n(\gamma_i)\cup n(|Z|\cap W_i))}\omega\log{f_i}+2\pi\I\int_{\xi\setminus n(|Z|\cap\xi)}\omega \eq\\ 
\sum_i\int_{W_i\setminus(n(\gamma_i)\cup n(|Z|\cap W_i))}d \eta_0\log{f_i}+2\pi\I\int_{\xi\setminus n(|Z|\cap\xi)}d \eta_0.
\end{multline*}
We transform the second term into
\[
2\pi\I\int_{\gamma} \eta_0 - 2\pi\I\int_{\partial n(|Z|\cap\xi)} \eta_0
\]
and the $i$-th summand in the first term into
\[
-\int_{\partial n(\gamma_i)}\eta_0\log{f_i}-\int_{\partial n(|Z|\cap W_i)}\eta_0\log{f_i}-\int_{W_i\setminus(n(\gamma_i)\cap n(|Z|\cap W_i))}\frac{d f_i}{f_i}\wedge\eta_0,
\]
where the last term equals to $0$ 
because $\frac{d f_i}{f_i}\wedge\eta_0\in F^{n-k+2}\A^{2n-2k+2}$.

Consider the integral
\[
\int_{\partial n(|Z|\cap W_i)}\eta_0\log{f_i}.
\]
Let $p$ be an intersection point of $|Z|$ and $W_i$. Then there is a neighborhood $U$ of $p$ which is analytically isomorphic to the product of open balls $B_W$ and $B_Z$, and $W_i\cap U$ maps to $B_W\times \{0\}$,  $|Z|\cap U$ maps to $\{0\}\times B_Z$. Let $n(|Z|\cap W_i)$ have only one connected component in $U$ and this is $D_W\times \{0\}$ where $D_W\subset B_W$ is a closed ball. We extend $f_i$ to $U$ by means of the projection $U\To B_W$. Let $\chi_\varepsilon(t)$ be a family of smooth functions on $\R$ which approximate $\delta_0$ as in \cite{GH}. For any current $T$ on $U$ we put
\[
T_\varepsilon\eq T*(\chi_\varepsilon^{2n})
\]
where $*$ denotes the convolution and 
\begin{multline*}
\chi_\varepsilon^{2n}(w_1, \dots, w_{2n-2k+2},b_1,\dots,b_{2k-2}) \\ \eq \chi_\varepsilon(w_1)\dots\chi_\varepsilon(w_{2n-2k+2})\chi_\varepsilon(b_1)\dots\chi_\varepsilon(b_{2k-2}).
\end{multline*}
Then
\[
\int_{\partial D_W \times \{0\}} \eta_0\log{f_i} \eq \lim_{\varepsilon\rightarrow 0} \int_U (\delta_{\partial D_W \times \{0\}})_\varepsilon\wedge\eta_0\log{f_i}
\]
since $\eta_0$ is smooth outside $\{0\}\times B_Z$. Applying the identity
\[
(\delta_{\partial D_W \times \{0\}})_\varepsilon \eq -d (\delta_{D_W \times \{0\}})_\varepsilon
\]
the last expression equals to
\[
-\lim_{\varepsilon\rightarrow 0} \int_U d (\delta_{D_W \times \{0\}})_\varepsilon \wedge \eta_0 \log{f_i} \eq \lim_{\varepsilon\rightarrow 0}\, (\eta_0 \log{f_i}, d (\delta_{D_W \times \{0\}})_\varepsilon),
\]
where we treat $\eta_0 \log{f_i}$ as a current and $d (\delta_{D_W \times \{0\}})_\varepsilon$ as a form. This is, by the definition of the differential for currents,
\[
\lim_{\varepsilon\rightarrow 0}\, (d(\eta_0 \log{f_i}), (\delta_{D_W \times \{0\}})_\varepsilon).
\]
Now we expand
\[
d(\eta_0 \log{f_i}) \eq (\omega-\delta_Z)\log{f_i} \,+\, \frac{d f_i}{f_i}\wedge\eta_0,
\]
thus obtaining
\begin{multline*}
\lim_{\varepsilon\rightarrow 0}\, ((\omega-\delta_Z)\log{f_i} + \tfrac{d f_i}{f_i}\wedge\eta_0, (\delta_{D_W \times \{0\}})_\varepsilon) \\ \eq\int_{D_W\times\{0\}} \omega \log{f_i} - \log{f_i(p)}\cdot\ord_p(Z\cdot W_i) + \lim_{\varepsilon\rightarrow 0} \, (\tfrac{d f_i}{f_i}\wedge\eta_0, (\delta_{D_W \times \{0\}})_\varepsilon),
\end{multline*}
the last summand being zero because 
\[
\eta_0 \;\in\; F^{n-k+1}\D^{2n-2k+1}, \qquad
(\delta_{D_W \times \{0\}})_\varepsilon \;\in\; F^{k-1}\A^{2k-2}.
\]
Note that when the radius of the ball $D_W$ tends to zero the first summand tends to zero, so can be neglected. Therefore the limit value of 
\[
\int_{\partial n(|Z|\cap W_i)}\eta_0\log{f_i}
\]
is
\[
-\sum_{p\in W_i\cap |Z|} \log{f_i(p)}\cdot\ord_pZ.
\]

The sum
\[
\sum_i\int_{\partial n(\gamma_i)}\eta_0\log{f_i}
\]
annihilates (in the limit) the integral
\[
2\pi\I\int_{\partial\xi} \eta_0.
\]

The remaining summand
\[
-2\pi\I\int_{\partial n(|Z|\cap\xi)} \eta_0
\]
tends to $2\pi\I$ times the intersection number of $Z$ and $\xi$ according to a reasoning similar to the one used above.
\end{proof}

\subsection{Construction of the fundamental class}
We would like to produce another proof of this theorem without the use of currents. For this we recall some cohomology constructions. What follows can be done for any smooth variety $X$ over $\C$.

Let $\Omega^\bullet_X$ be the holomorphic de Rham complex of $X$. For any integer $j$ denote by $F^j\Omega^\bullet_X$ the subcomplex
\[
F^j\Omega^i_X\eq\begin{cases}0& \text{if $i<j$,}\\ \Omega^i_X& \text{if $i\ge j$.}\end{cases}
\]
There are natural maps $F^j\Omega^\bullet_X\rightarrow\Omega^\bullet_X$ and $F^j\Omega^\bullet_X\rightarrow\Omega^j_X$.

Let $Y\subset X$ be a smooth subvariety of codimension $j$. Recall the construction of the fundamental class of $Y$ in $\HC^{2j}_Y(X,F^j\Omega^\bullet_X)$ (see \cite{groth:fga}, Expos\'e 149 and \cite{bloch:semireg}). 

We first construct the Hodge class $c^H(Y)\in H^j_Y(X,\Omega^j_X)$. There is a spectral sequence
\[
E^{p,q}_2\eq H^p(X,\SC_Y^q(\Omega^j_X)) \;\Rightarrow\; H_Y^{p+q}(X,\Omega^j_X).
\]
Since $\Omega^j_X$ is locally free $\SC_Y^q(\Omega^j_X)=0$ for $q<j$. This implies that 
\[
H_Y^j(X,\Omega^j_X) \eq \Gamma(X,\SC_Y^j(\Omega^j_X)).
\]
Let $V$ be an open subset of $X$ on which $Y$ is a complete intersection, so there exist regular functions $f_1$, $f_2$,\dots, $f_j$ on $V$ which generate the ideal of $Y\cap V$. Put $V_i=V\setminus \{f_i=0\}$. Then $V_i$ form a covering of $V\setminus(V\cap Y)$. So we can consider the \v Cech cohomology and the section
\[
(2\pi\I)^{-j}\frac{d f_1}{f_1}\wedge\dots\frac{d f_j}{f_j} \;\in\; \Gamma(\bigcap V_i, \Omega_X^j)
\]
produces an element of $H^{j-1}(V\setminus(V\cap Y), \Omega_X^j)$. We obtain an element of $H^j_{V\cap Y}(V,\Omega_X^j)$ by applying the boundary map of the long exact sequence
\[
\dots\To H^{j-1}(V\setminus(V\cap Y), \Omega_X^j)\To H^j_{V\cap Y}(V,\Omega_X^j)\To H^j(V,\Omega_X^j) \To\dots.
\]
This element is a section of the sheaf $\SC_Y^j(\Omega^j_X)$ over $V$. These local sections glue together to produce a global section
\[
c^H(Y) \;\in\; H^j_Y(X,\Omega^j_X).
\]

The differential $d:\Omega^j_X\rightarrow\Omega^{j+1}_X$ induces the differential on cohomology $d:H^j_Y(X,\Omega^j_X)\rightarrow H^j_Y(X,\Omega^{j+1}_X)$. We have $d c^H(Y)=0$ since for $V$, $f_i$ as above
\[
d\left( (2\pi\I)^{-j}\frac{d f_1}{f_1}\wedge\dots\frac{d f_j}{f_j}\right) \eq 0.
\]
There is a spectral sequence for the hypercohomology
\[
E^{p,q}_2\eq H^p(H^q_Y(X,F^j\Omega_X^\bullet)) \;\Rightarrow\; \HC^{p+q}_Y(X,F^j\Omega_X^\bullet),
\]
which shows that 
\begin{align*}
\HC^{2j}_Y(X,F^j\Omega_X^\bullet)
&\eq H^j(H^j_Y(X,F^j\Omega_X^\bullet))\\
&\eq\kernel(d: H^j_Y(X,\Omega_X^j)\rightarrow H^j_Y(X,\Omega_X^{j+1})).
\end{align*}
Therefore the natural map $\HC^{2j}_Y(X,F^j\Omega_X^\bullet)\rightarrow H^j_Y(X, \Omega_X^j)$ is an injection and $c^H(Y)$ lifts to a unique $c^F(Y)\in\HC^{2j}_Y(X,F^j\Omega_X^\bullet)$. The natural map 
\[
\HC^{2j}_Y(X,F^j\Omega_X^\bullet)\To \HC^{2j}_Y(X,\Omega_X^\bullet) 
\;\cong\; H^{2j}_Y(X,\C)
\]
sends $c^F(Y)$ to an element $c^{DR}(Y)\in H^{2j}_Y(X,\C)$. 

We will prove now that $c^{DR}(Y)$ is the Thom class of $Y$, i.e. that its value on a class in $H_{2j}(X,X\setminus Y)$ equals to the intersection number of this class with $Y$. Since the real codimension of $Y$ is $2j$, $\SC^p_Y(\C)=0$ for $p<2j$. Therefore $H^{2j}_Y(X,\C)=\Gamma(X,\SC^{2j}_Y(\C))$ and $c^{DR}(Y)$ can be described locally. Let $V$ be a neighborhood of $X$ which is isomorphic to a product of unit disks, $V\cong \DD^n$, and such that $V\cap Y$ is given by equations $z_i=0$ for $i=1,\dots j$, where $z_i$ is the coordinate on $i$-th disk. Then 
\[
H^{2j}_{Y\cap V}(V,\C)
\eq H^{2j}(\DD^n,\, (\DD^j\setminus\{0\}^j)\times \DD^{n-j};\, \C)
\eq\C.
\]
So to evaluate the restriction $c^{DR}(Y)|_V$ it is enough to evaluate $c^{DR}(Y)|_V$ on the generator of the homology $H_{2j}(\DD^n,(\DD^j\setminus\{0\}^j)\times \DD^{n-j};\C)$, the transverse class $\DD^j\times \{0\}^{n-j}$. Since the class $c^{H}(Y)|_V$ is coming from the boundary map $H^{j-1}((\DD^j\setminus\{0\}^j)\times \DD^{n-j},\Omega^j)\rightarrow H^j(\DD^n,(\DD^j\setminus\{0\}^j)\times \DD^{n-j};\Omega^j)$
\[
\langle c^{DR}(Y)|_V, \DD^j\times \{0\}^{n-j}\rangle \eq 
\langle c^{H}(Y)|_V, \DD^j\times \{0\}^{n-j}\rangle \eq 
\langle c_0, \partial(\DD^j\times \{0\}^{n-j})\rangle,
\]
where $c_0$ is the corresponding class in $H^{j-1}((\DD^j\setminus\{0\}^j)\times \DD^{n-j},\Omega^j)$. The class $c_0$, in its turn, comes from the map
\[
H^0(\DD_0^j\times\DD^{n-j}, \Omega^j) \To H^{j-1}((\DD^j\setminus\{0\}^j)\times \DD^{n-j},\Omega^j).
\]
This map comes from \v Cech cohomology and can also be constructed using successive application of Mayer-Vietoris exact sequences. Hence we have the dual map in homology
\[
H_{2j-1}((\DD^j\setminus\{0\}^j)\times \DD^{n-j},\C)\To
H_j(\DD_0^j\times\DD^{n-j},\C)
\]
and one can see that the image of $\partial(\DD^j\times \{0\}^{n-j})$ under this map is $\SB^j\times\{0\}^{n-j}$, where $\SB$ is the unit circle. Indeed, let 
\[
U_k\eq\DD_0^k\times(\DD^{j-k}\setminus\{0\}^{j-k})\times\{0\}^{n-j}.
\]
Then $U_0=(\DD^j\setminus\{0\}^j$ and $U_{j-1}=\DD_0^j\times\DD^{n-j}$. Put 
\[
X_k^1\eq\DD_0^k\times\DD^{j-k}\times\{0\}^{n-j},\qquad
X_k^2\eq\DD_0^{k-1}\times\DD\times(\DD^{j-k}\setminus\{0\}^{j-k})\times\{0\}^{n-j}.
\]
Then
\[
X_k^1\cap X_k^2 \eq U_k,\qquad X_k^1\cup X_k^2 \eq U_{k-1},
\]
so there is a boundary map in the Mayer-Vietoris sequence associated to $X_k^1$, $X_k^2$ which goes from $H_{2j-k}(U_{k-1})$ to $H_{2j-k-1}(U_{k})$. We are going to show, by induction, that the $k$-th iterated image of $\partial(\DD^j\times \{0\}^{n-j})$ under these maps is 
\[
\SB^k\times\partial(\DD^{j-k})\times\{0\}^{n-j} \;\in\; H_{2j-k-1}(U_{k},\C).
\]
The boundary map in the Mayer-Vietoris sequence can be decomposed (see \cite{dold}, p. 49) as
\[
H_{2j-k}(U_{k-1})\To H_{2j-k}(U_{k-1},X_k^2) 
\;\cong\; H_{2j-k}(X_k^1, U_k) \To H_{2j-k-1}(U_k).
\]
Writing
\[
\SB^{k-1}\times\partial(\DD^{j-k+1})\times\{0\}^{n-j} \eq \SB^k\times\DD^{j-k}\times\{0\} \,+\, \SB^{k-1}\times\DD\times\partial(\DD^{j-k})\times\{0\}^{n-j}
\]
we see that the second summand is contained in $X_k^2$, so is trivial in homology group $H_{2j-k}(U_{k-1},X_k^2)$. The first summand belongs to $H_{2j-k}(X_k^1, U_k)$, so it remains to take its boundary, which is exactly $\SB^k\times\partial(\DD^{j-k})\times\{0\}^{n-j}$.

Therefore
\[
\langle c_0, \partial(\DD^j\times \{0\}^{n-j})\rangle \eq \int_{\SB^j} (2\pi\I)^{-j} \frac{d z_1}{z_1}\wedge\dots\frac{d z_j}{z_j} \eq 1,
\]
which means that the constructed class $c^{DR}(Y)$ is indeed the Thom class of $Y$.

It is clear now that the following theorem is true:
\begin{thm}\label{thm52}
Let $X$ be a smooth variety over $\C$ of dimension $n$. Let $Y\subset X$ be a smooth closed subvariety of $X$ of codimension $j$. Then there is a class $c^F(Y)\in \HC_Y^{2j}(X, F^j\Omega_X^\bullet)$ which satisfies the following conditions:
\begin{enumerate}
\item The image $c^{DR}(Y)$ of $c^F(Y)$ in $H_Y^{2j}(X, \C)$ is the Thom class of $Y$.
\item The class $c^H(Y)$, which is the image of $c^F(Y)$ in $H_Y^j(X,\Omega_X^j)$, is logarithmic, i.e. for any open set $V\subset X$ and a holomorphic function $f$ on $V$ which is zero on $Y\cap V$ the product $f\cdot c^H(Y)|_V$ is $0$.
\end{enumerate}
\end{thm}

We extend the definition of $c^F$, $c^{DR}$, $c^H$ to formal linear combinations of subvarieties in the obvious way.

\subsection{Dolbeault local cohomology}
We show how to interpret the results of the previous section using smooth forms. For this we show how local cohomology can be computed using Dolbeault resolutions. Let $X$ be a smooth variety over $\C$ and $Y$ be a subvariety of codimension $j$, $U=X\setminus Y$. Let $\jmath$ be the inclusion $U\rightarrow X$.

\begin{prop}
Let $S$ be a soft sheaf on $X$ which locally has no nonzero sections supported on $Y$. Then $H^p_Y(X,S)$ and $\SC^p_Y(S)$ are zero unless $p=1$,
\begin{align*}
H^1_Y(X,S) &\eq \cokernel(\Gamma(X,S)\To\Gamma(U,S)), \\
\SC^1_Y(S) &\eq \cokernel(S\To \jmath(S|_U)).
\end{align*}
\end{prop}
\begin{proof}
For any open set $V\subset X$ we have the long exact sequence for local cohomology:
\[
\dots\To H^p_{V\cap Y}(V,S) \To H^p(V,S)\To H^p(V\cap U,S)\To\dots.
\]
The groups $H^p(V,S)$ and $H^p(V\cap U,S)$ vanish for $p>0$. Hence $H^p_{V\cap Y}(V,S)$ vanish for $p>1$. This group also vanishes for $p=0$ by the condition on $S$. For $V=X$ this implies the statement about the groups $H^p_Y(X,S)$. Since $\SC^p_Y(S)$ is the sheaf associated to the presheaf $(V\rightarrow H^p_{V\cap Y}(V,S))$, this implies the statement for the groups $\SC^p_Y(S)$.
\end{proof}

Therefore one can compute local cohomology in the following way:
\begin{prop}
Let $F^\bullet$ be a bounded complex of sheaves on $X$ and $F^\bullet\rightarrow S^\bullet$ a bounded soft resolution, such that each $S^i$ locally has no nonzero sections supported on $Y$. Then
\[
\HC^i_Y(X, F^\bullet) \;\cong\; H^i(S^\bullet(X,U)),\qquad 
\SC^i_Y(F^\bullet) \;\cong\; H^i(S^\bullet_{X,U}),
\]
where $S^\bullet_{X,U}$ is the complex of sheaves defined as follows:
\[
S^i_{X,U}\eq S^i \,\oplus\, \jmath_*(S^{i-1}|_U), \qquad 
d(a,b) \eq (da, -db-a|_U),
\]
where $a$, $b$ are sections of $S^i$, $\jmath_*(S^{i-1}|_U)$, and $S^\bullet(X,U)=\Gamma(X, S^\bullet_{X,U})$.
\end{prop}
\begin{proof}
Consider the following spectral sequence:
\[
E_2^{pq}\eq H^p(H^q_Y(X,S^\bullet)) \;\Rightarrow\; \HC^{p+q}_Y(X,S^\bullet).
\]
By the proposition above the spectral sequence degenerates and 
\begin{align*}
\HC^i(X,S^\bullet) 
& \;\cong\; H^{p-1}(\cokernel(\Gamma(X,S^\bullet)\To \Gamma(U,S^\bullet))) 
\\ &\eq H^p(\cone(\Gamma(X,S^\bullet)\To \Gamma(U,S^\bullet))[-1]) 
\eq H^p(S^\bullet(X,U)).
\end{align*}
The statement about $\SC^i_Y(F^\bullet)$ can be proved similarly.
\end{proof}

In particular the proposition works for the following resolutions:
\[
\Omega_X^j\To (\A_X^{j\bullet}, \bar\partial),\qquad F^j\Omega^\bullet\To (F^j\A_X^\bullet, d).
\]

\subsection{A proof using relative cohomology}
We now give the second proof of Theorem \ref{thm:spec_val1}.
\begin{proof}
Let $j=n-k+1$, the codimension of $Z$ and the dimension of $W_i$.
By the definition (*) of $AJ^{k,1}$ we have
\[
(AJ^{k,1}[x], \cl^{n-k+1} Z) \eq \sum_i\int_{W_i\setminus\gamma_i}\omega\log{f_i} \,+\, 2\pi\I\int_\xi\omega,
\]
where $\omega\in F^{n-k+1}\A^{2n-2k+2}(X)$ is a form whose class $[\omega]$ equals to the Poincar\'e dual of the class $[Z]\in H_{2k-2}(X,\Z)$.

In the definition of the Abel-Jacobi map we choose $\gamma_k$ and $\xi$ to be transversal to $|Z|$. This means that $\gamma_k$ does not intersect $|Z|$ for each $k$ and $\xi$ intersects $|Z|$ only in several points. Let $U=X\setminus |Z|$.

Let $(\omega,\eta)\in F^j\A^{2j}(X, U)$ be a representative of $c^F(Z)$. Then $(\omega,\eta)$ is also a representative of $c^{DR}(Z)$, the Thom class of $Z$. This means that $\omega\in F^{j}\A^{2j}(X)$, $\eta\in F^{j}\A^{2j-1}(U)$, $\omega=-d\eta$ on $U$, and for any $2j$-chain $c$ on $X$ with boundary in $U$
\[
\int_c \omega + \int_{\partial c} \eta = [c]\cdot Z.
\]
It is clear that we may choose $\omega$ as a representative of $[Z]$.

Since $(\omega,\eta)$ represent the Thom class of $Z$, we have the following identity:
\[
\int_\xi\omega \,+\, \int_{\partial\xi}\eta \eq [\xi] \cdot Z,
\]
where $[\xi]\in H_{2j}(X, U)$ is the class of $\xi$.
Therefore
\[
2\pi\I\int_\xi\omega \;\equiv\; -2\pi\I\int_\gamma\eta \mod{2\pi\I\Z}
\]
and we have the following decomposition:
\[
(AJ^{k,1}[x], \cl^{n-k+1} Z) \eq \sum_i(\int_{W_i\setminus\gamma_i}\omega\log{f_i} - 2\pi\I\int_{\gamma_i}\eta).
\]
Consider each summand separately. Let $n(\gamma_i)$ be a small neighborhood of $\gamma_i$ in $W_i$. Then, up to terms which tend to zero when the radius of the neighborhood tends to zero, we may write
\[
2\pi\I\int_{\gamma_i}\eta \eq \int_{\partial n(\gamma_i)} \eta \log{f_i}
\]
and
\[
\int_{W_i\setminus\gamma_i}\omega\log{f_i} - 2\pi\I\int_{\gamma_i}\eta
\eq\int_{W_i\setminus n(\gamma_i)} \omega\log{f_i} + \int_{\partial(W_i\setminus n(\gamma_i))} \eta \log{f_i}.
\]
We see that this is nothing else than the pairing of classes
\begin{align*}
[W_i\setminus n(\gamma_i)] & \;\in\; H_{2j}(W_i\setminus \gamma_i, W_i\setminus (\gamma_i\cup Z_i)), \\
[(\omega|_{W_i}\log{f_i},\, \eta|_{W_i\setminus Z_i}\log{f_i})] &\;\in\; H^{2j}_{Z_i}(W_i\setminus \gamma_i),
\end{align*}
where $Z_i=|Z|\cap W_i$. Indeed, 
\[
d(\eta|_{W_i\setminus Z_i}\log{f_i}) \eq -\omega|_{W_i\setminus Z_i} \log{f_i} 
\]
since $\eta\in F^j(U)$, so the pair $(\omega|_{W_i}\log{f_i}, \eta|_{W_i\setminus Z_i}\log{f_i})$ defines a class in the local cohomology.

For each point $p$ in the finite set $W_i\cap |Z|$ we choose a small neighborhood $n(p)$ of $p$ in $W_i$ in such a way that the $n(p)$ do not intersect each other and do not intersect $\gamma_i$. Then
\[
\int_{W_i\setminus n(\gamma_i)} \omega\log{f_i} + \int_{\partial(W_i\setminus n(\gamma_i))} \eta \log{f_i} \eq \sum_{p\in |Z|\cap W_i} (\int_{n(p)} \omega\log{f_i} + \int_{\partial n(p)} \eta\log{f_i}).
\]

Since the dimension of $n(p)$ is $j$, $F^j\Omega^\bullet_{n(p)}$ is just $\Omega^j_{n(p)}$ in degree $j$. Consider the multiplication by the function $\log{f_i}-\log{f_i(p)}$, which acts on $\Omega^j_{n(p)}$. Since the class $[(\omega, \eta)|_{n(p)}]\in H^j_p(n(p), \Omega^j)$ is the restriction of the class $c^H(Z)$, which is logarithmic, it is logarithmic itself, so the multiplication by $\log{f_i}-\log{f_i(p)}$ kills it. Hence
\begin{align*}
\int_{n(p)} \omega\log{f_i} + \int_{\partial n(p)} \eta\log{f_i} 
&\eq \left(\int_{n(p)} \omega + \int_{\partial n(p)} \eta \right) \log{f_i(p)} \\
&\eq \ord_p(Z\cdot W_i)\log{f_i(p)}
\end{align*}
and the assertion follows.
\end{proof}

\begin{rem}
One could also consider the product 
\[
CH^k(X,1)\times CH^{n-k+1}(X) \rightarrow CH^{n+1}(X,1)
\]
which sends $(\sum_i(W_i,f_i) , Z)$ to 
\[
\sum_i\sum_{p\in W_i\cap |Z|} \ord_p(Z\cdot W_i) \, (p, f_i(p)).
\]
The Abel-Jacobi map acts
\[
CH^{n+1}(X,1)\To \frac{H^0(X,\C)^*}{2\pi\I\, H^0(X,\Z)} \;\cong\; \C/2\pi\I\,\Z
\]
sending $(p, a)$, $p\in X$, $a\in \C^\times$ to $\log a$, so the theorem proved above simply means that the following diagram commutes:
\[
\begin{small}
\begin{CD}
CH^k(X,1)\times CH^{n-k+1}(X) @>>> CH^{n+1}(X,1)\\
@V{AJ^{k,1}\times cl^{n-k+1}(X)}VV @V{AJ^{n+1,1}}VV\\
\begin{matrix} 
{\displaystyle \frac{(F^{n-k+1}H^{2n-2k+2}(X,\C))^*}{2\pi\I\,H^{2k-2}(X,\Z)}} \times \\
\times (F^{n-k+1}H^{2n-2k+2}(X,\C)\cap H^{2n-2k+2}(X,\Z)) 
\end{matrix}
@>>> \C/\Z
\end{CD}
\end{small}
\]
There is a more general statement that the regulator map from higher Chow groups into the Deligne cohomology is compatible with products. This is mentioned in \cite{bloch:acbc}. The construction of higher Chow groups there is different from the one above, but it can be proved that they are canonically isomorphic and the regulator map corresponds to  the Abel-Jacobi map (see \cite{KLM}).
\end{rem}

\chapter{Derivatives of the Abel-Jacobi map}
This chapter consists of two parts. In the first part we explain how we represent cohomology classes using algebraic forms. 

For any open cover of a topological space $X$ and any sheaf $F$ of abelian groups on $X$ one has the \v{C}ech complex, which is a complex formed by the groups of sections of $F$ over all possible intersections of the sets of the cover. We use a slightly more general approach which involves hypercovers. Roughly speaking, hypercover is a system of open subsets of $X$ organised in a certain way so that for any sheaf of abelian groups $F$ we obtain a complex of abelian groups formed by the groups of sections of $F$ over the sets of the hypercover. Note that the open subsets in the definition of the \v{C}ech complex are naturally parametrised by the faces of a simplex. Similarly, the open subsets of a hypercover are parametrised by the cells of a certain cell complex $\sigma$. In our examples $\sigma$ will be a cube.

If $X$ is an algebraic variety with a hypercover $\U$ for the Zariski topology, one obtains a complex $\Omega^i(\U)^\bullet$ of groups corresponding to the sheaf $\Omega^i_X$ of algebraic de Rham forms of degree $i$ for each $i$. Thus we obtain a bicomplex. The total complex of this bicomplex is denoted $\Omega^\bullet(\U)$ and the elements of $\Omega^i(\U)$ are called ``hyperforms'' of degree $i$. Thus a hyperform of degree $i$ is given by several forms of degree $i$ on certain open subsets of $X$, several forms of degree $i-1$, etc. There is a canonical homomorphism from $H^i(\Omega^\bullet(\U))$ to the algebraic de Rham cohomology $\HC^i(X, \Omega_X)$ of $X$, which is an isomorphism if $\U$ is composed of affine open sets. Dually one obtains a notion of ``hyperchains'' so that one can integrate a hyperform along a hyperchain. We obtain the Stokes theorem for hyperforms, and certain structures on hyperforms, namely the Hodge filtration, products, and the Gauss-Manin connection if $X$ varies in a family. We emphasize that the hyperchains are composed of \emph{topological} chains on the manifold $X(\C)$ with its \emph{analytic} topology, while the hyperforms are composed of \emph{algebraic} forms. This way of representing cohomology classes is extremely convenient when $X$ is a product of curves because one can easily write down representatives for classes obtained by external multiplication of classes on curves, while the classes on curves can be represented by differentials of \emph{second kind}.

The definition and properties of hyperforms are given in Sections \ref{hypercovers}~--\ref{products} and \ref{gauss-manin}.

In Sections \ref{residues}~--\ref{comp_res} we study certain residues and trace maps associated with hypercovers, which are first defined as certain integrals, and then, in the case when our variety is embedded in a product of curves and the hypercover is the product hypercover, we give a way to compute these residues and trace maps in an algebraic way using iterated residues. 

Our trace maps are defined in the case when $X$ is a product of curves and on each curve a Zariski cover is given. Suppose a finite family $M$ of irreducible subvarieties of $X$ is given. The elements of $M$ will be called \emph{special}. All further constructions do not change if one adds varieties to $M$, so one possible way to think is to consider limits over all possible $M$, but we will not pursue this point of view in order to simplify the exposition. Let $Z$ be a special variety and let $Z^0$ denote the complement in $Z$ of all proper special subvarieties of $Z$. Then we construct a family of classes $h_a(Z)$ of dimension $\dim Z$ in the homology of $Z^0$ ($a$ runs over the set of $\dim Z$-dimensional cells of $\sigma$). Then the residues $\res_{Z,a}^{\Int}$ are defined for forms from $\Omega^{\dim Z} (Z^0)$ as the integrals along $h_a(Z)$ (``$\Int$'' stands for ``integral''). The trace map $\Tr_Z^{\Int}$ is a linear functional on the space of hyperforms $\Omega^{2\dim Z}(\U\cap Z^0)$, obtained by summing up the residues of the corresponding components of a hyperform (see Section \ref{residues}).

Then we compute these residues algebraically in terms of iterated residues. We introduce $\Tr_Z$ (without ``$\Int$'') as a certain combination of iterated residues so that $\Tr_Z^\Int$ is a simple multiple of $\Tr_Z$ (see Section \ref{comp_res}). In the end of Section \ref{comp_res} we summarize by formulating a recipe for computation of the trace map.

The second part of this chapter is motivated by the following idea. Suppose we have a family of varieties $\{X_s\}_{s\in S}$ depending on certain parameters. Suppose we have a family of higher cycles $\{x_s\in Z^k(X_s,1)\}_{s\in S}$. Let us denote the value of the Abel-Jacobi map on $x_s$ by $AJ_s$. The element $AJ_s$ belongs to the quotient of a cohomology group by a subspace and a lattice. One notices that one can ``kill'' the lattice by differentiating $AJ_s$ with respect to parameters.

The objects $AJ_s$ are defined in terms of their pairing with cohomology classes from certain step in the Hodge filtration. In fact for any hyperform $\omega$ of right degree (without the condition that $\omega$ belongs to a certain step in the Hodge filtration) we may define $(AJ_s, \omega)$, but this value depends on some choices. We obtain two natural formulae (Propositions \ref{prop:diff_aj:1}, \ref{prop:diff_aj:2}). The first one gives $d (AJ_s, \omega) - (AJ_s, \nabla\omega)$ ($\nabla$ is the Gauss-Manin derivative), and the second one gives $(AJ_s, d \eta)$, where $\eta$ is a hyperform of degree one less. Note that our Gauss-Manin derivative is actually a lift of the ordinary Gauss-Manin derivative, which is defined on cohomology classes, to hyperforms, and it depends on some choices.

In Section \ref{subs:ext_dmod} we define two maps $\Psi_0$ and $\Psi_1$, which give values of $(AJ_s, d \eta)$ and $d (AJ_s, \omega) - (AJ_s, \nabla\omega)$ correspondingly. Although the definition of $AJ_s$ is purely transcendental, the maps $\Psi_0$ and $\Psi_1$ can be computed algebraically using the trace map. The map $\Psi_1$ is defined on hyperforms of degree $2n-2k+2$ ($n$ is the dimension of $X$) and takes values in the sheaf of differential $1$-forms on the base $S$. The map $\Psi_0$ is defined on hyperforms of degree $2n-2k+1$ and takes values in the sheaf of functions on the base $S$. The maps $\Psi_0$ and $\Psi_1$ satisfy a list of axioms (see Propositions \ref{prop:ext-dmod:1}, \ref{prop:ext-dmod:2}) which allows us to construct a sheaf of $\D$-modules, which is an extension of the cohomology sheaf $H^{2n-2k+2}(X_s, \C)$ by the structure sheaf of $S$. This extension contains all the information about the derivative of the Abel-Jacobi map.

The results mentioned above are stated only for the case when $X_s$ is a product of curves because they depend on our definition of trace maps.

Next, in Section \ref{prod-ellc} we specialise to the case when $X_s$ is the product of an elliptic curve by itself $2k-2$ times. In this case to the extension of $\D$-modules mentioned above corresponds a certain invariant $\Psi_{\analyt}'(B_{2k-2})$, which becomes a meromorphic modular form of weight $2k$ if the family is an open subset of a modular family. 

Finally, in Section \ref{subs:anal_comp} we show that if the modular form obtained from the extension of $\D$-modules obtained from a family of higher cycles is proportional to the modular form $g_{z_0}^{k, \HH/\Gamma}$, obtained by taking derivatives of the Green function, then (in the case when there are no cusp forms) we obtain a formula for the values of the Green function, which proves the algebraicity conjecture (Section \ref{global_green}) for the given $k$, $\Gamma$, $z_0$, and arbitrary second CM point $z$ (see Theorem \ref{thm:anal-comp:9}).

\section{Hypercovers}\label{hypercovers}
\subsection{Abstract cell complex}
\begin{defn}
An \emph{abstract cell complex} is a graded set
\[
\sigma \eq \bigcup_{i\geq 0} \sigma_i
\]
together with homomorphisms $d_i:\Z[\sigma_i]\To\Z[\sigma_{i-1}]$ such that $d_i\circ d_{i+1}=0$ and $\veps\circ d_0=0$, where $\veps$ is the \emph{augmentation map} $\Z[\sigma_0]\To\Z$. The elements of the set $\sigma_i$ are called \emph{cells} of dimension $i$ and the homomorphisms $d_i$ are called \emph{boundary maps}.
\end{defn}
\[
\begin{CD}
\cdots @>d>>\Z[\sigma_2] @>d>> \Z[\sigma_1] @>d>> \Z[\sigma_0] \\
& & & & & & @V{\veps}VV\\
& & & & & & \Z
\end{CD}
\]

If an abstract cell complex is given we denote it usually by $\sigma$ and encode the boundary maps by the coefficients $D_{ab}$:
\[
d a \eq \sum_{b\in\sigma_{i-1}} D_{ab} b, \qquad (a\in\sigma_i).
\]

The basic example is the standard simplex of dimension $n$, $\Delta_n$. In this case $\sigma_i$ is the set of all subsets of $\{0,1,\ldots,n\}$ of size $i+1$ and, writing subsets as increasing sequences, 
\[
d_i(j_0,\ldots,j_i) \eq \sum_{k=0}^i (-1)^k (j_0,\ldots,\hat{j}_k,\ldots,j_i).
\]

Given two abstract cell complexes $\sigma$, $\sigma'$ the \emph{product} 
\[
\sigma'' \eq \sigma\times \sigma'
\]
is again an abstract cell complex in the following way:
\[
\sigma''_i \eq \bigcup_{j=0}^i \sigma_j\times\sigma'_{i-j},
\]
\[
d_i(a\times b) \eq (d_j a) \times b \;+\; (-1)^j\, a\times d_{i-j} b \qquad \text{for $a\in\sigma_j$, $b\in\sigma'_{i-j}$.}
\]
For example, the cube is defined as
\[
\square_n \eq \Delta_1\times\cdots\times\Delta_1,\qquad\text{the product has $n$ terms.}
\]

\subsection{Hypercover}
Let $X$ be a topological space and $\sigma$ an abstract cell complex. 

\begin{defn}
A \emph{hypercover} of $X$ (indexed by $\sigma$) is a system of open subsets $\U=(U_a)$ indexed by elements $a\in\sigma$ such that:

\begin{enumerate}
\item Whenever a cell $b$ belongs to the boundary of a cell $a$, $U_a\subset U_b$.
\item For any point $x\in X$ consider all cells $a$ such that $x\in U_a$. These cells form a subcomplex of $\sigma$ by the first property, call it $\sigma_x$. We require the complex of abelian groups coming from $\sigma_x$ to be a resolution of $\Z$ with the augmentation map induced by $\varepsilon$.
\end{enumerate}
\end{defn}

\begin{example}
Put $\sigma=\Delta_n$ and let $X$ be a space covered by $n+1$ open subsets $U_0,\dots U_n$. For any sequence $(j_0,\dots,j_i)\in\sigma_i$ we put 
\[
U_{(j_0,\dots,j_i)}\eq U_{j_0}\cap\dots\cap U_{j_i}.
\]
This clearly gives a hypercover which is called the \v{C}ech hypercover.
\end{example}

If $\sigma$ and $\sigma'$ are two abstract chain complexes, $(U_a)$ and $(U'_{a'})$ are hypercovers of spaces $X$, $X'$ indexed by $\sigma$ and $\sigma'$ correspondingly, by making all products $(U_a\times U'_{a'})$ one gets a hypercover of $X\times X'$ indexed by $\sigma\times\sigma'$, the \emph{product} hypercover.

\subsection{Hypersection}
Let $X$ be a topological space, $\sigma$ an abstract chain complex and $\U=(U_a)$ a hypercover.

Let $\Z_a$ denote the constant sheaf with fiber $\Z$ on $U_a$ extended by zero to $X$. Let 
\[
\Z_{\U i} \eq \bigoplus_{a\in \sigma_i} \Z_a.
\]
If $U_a\subset U_b$, there is a canonical morphism $r_{ab}:\Z_a \rightarrow \Z_b$. We define $d_i:\Z_{\U i} \rightarrow \Z_{\U i-1}$. If $a\in\sigma_i$ then the morphism from $\Z_a$ to $\Z_{\U i-1}$ is 
\[
\sum_{b\in\sigma_{i-1}}  D_{ab} r_{a b},\qquad \text{where}
\]
\[
d a \eq \sum_{b\in\sigma_{i-1}} D_{ab} b.
\]
Denote the corresponding sequence of sheaves by $\Z_\U$. The augmentation map $\veps:\Z_{\U 0} \rightarrow \Z_X$ is defined as the sum of the canonical morphisms $\veps_a:\Z_a\rightarrow\Z_X$ for $a\in\sigma_0$.
\[
\begin{CD}
\cdots @>d>> \Z_{\U 2} @>d>> \Z_{\U 1} @>d>> \Z_{\U 0} \\
& & & & & & @V{\veps}VV \\
& & & & & & \Z_X
\end{CD}
\]

\begin{prop}
The sequence $\Z_\U$ is a resolution of $\Z_X$.
\end{prop}
\begin{proof}
For any point $x\in X$ the stalk of $\Z_\U$ at $x$ is simply the complex of abelian groups corresponding to $\sigma_x$. So the statement follows from the second condition of hypercover.
\end{proof}

In other words, we have obtained a \emph{quasi-isomorphism}
\[
\veps:\Z_{\U} \rightarrow \Z_X.
\]
For any sheaf or complex of sheaves $\F$ on $X$ this gives a morphism of complexes
\[
\veps^*: \HHom(\Z_X, \F) \rightarrow \HHom(\Z_\U, \F).
\]
Note that $\HHom(\Z_X, \F) = \Gamma(X, \F)$ and for any $a\in\sigma$ $\HHom(\Z_a, \F) = \Gamma(U_a, \F)$. We denote the complex $\HHom(\Z_\U, \F)$ by $\Gamma(\U, \F)$ and call elements of this complex \emph{hypersections}. 

The complex of hypersections can be defined more precisely as follows. Let $\F$ be a complex of sheaves on $X$ written as
\[
\F^0\rightarrow\F^1\rightarrow\dots.
\]
Then we put
\[
\Gamma(\U,\F)^i \eq \prod_{j\geq0,a\in \sigma_j} \Gamma(U_a, \F^{i-j}).
\]
The coboundary of a hypersection $s=(s_a)\in \Gamma(\U,\F)^i$ is a hypersection $ds=(s'_a)\in \Gamma(\U,\F)^{i+1}$, where
\begin{equation}\label{dsa}
s'_a \eq d s_a + (-1)^{i-j+1} \sum_{b\in \sigma_{j-1}} D_{ab} s_b|_{U_a},\qquad a\in\sigma_j.
\end{equation}
One can check directly that $d^2=0$. Let $d^2s=(s''_a)$, $a\in\sigma_j$.
\begin{align*}
s''_a 
&\eq d s'_a \;+\; (-1)^{i-j} \sum_{b\in \sigma_{j-1}} D_{ab} s'_b|_{U_a} \\
&\eq
(-1)^{i-j+1} \sum_{b\in\sigma_{j-1}} D_{ab} d s_b|_{U_a} 
+ (-1)^{i-j} \sum_{b\in \sigma_{j-1}} D_{ab} d s_b|_{U_a} \\
& \phantom{XXXXXXXXXXXXX} 
+ \sum_{b\in \sigma_{j-1}} \sum_{c\in \sigma_{j-2}} D_{ab} D_{bc} s_c|_{U_a} \\
& \eq 0.
\end{align*}

The augmentation morphism $\veps^*:\Gamma(X,\F)\rightarrow \Gamma(\U,\F)$ is defined by sending a section $s\in\Gamma(X,\F^i)$ to the hyperchain $(s'_a)$ where $s'_a$ is not zero only for $a\in \sigma_0$ and equals to the restriction of $s$ to $U_a$.

\begin{rem}
In the case of the \v{C}ech hypercover (corresponding to an open cover) we obtain the \v{C}ech complex.
\end{rem}
\begin{rem}
If I took the definition of $\HHom$ from \cite{harts:rd}, p. 64 the sign in (\ref{dsa}) would be $(-1)^{i+1}$. I assume that switching from a chain complex $\Z_\U$ to the corresponding cochain complex adds the multiplier $(-1)^j$.  My choice of sign is made in order to get the sign  correct in the case of the \v{C}ech hypercover.
\end{rem}

\subsection{Hyperchain}
We assume some class of topological chains on $X$ is given, i.e. semialgebraic chains. For any open set $U\subset X$ we denote by $C_i(U)$ the group of chains with support on $U$ of dimension $i$, which is the free abelian group generated by maps of the chosen class from the standard simplex of dimension $i$ to $U$. Suppose an abstract cell complex $\sigma$ and a hypercover $(U_a)$ is given. 

The complex of \emph{hyperchains} is defined similarly to the complex of hypersections, but with a change of sign. We put
\begin{align*}
C_i(\U) &\eq \bigoplus_{j\geq 0} C_{i-j}(\U, j), \qquad\text{where} \\
C_{i}(\U,j) &\eq \bigoplus_{a\in \sigma_j} C_{i}(U_a).
\end{align*}
Given a hyperchain $\xi\in C_{i-j}(U_a)\subset C_i(\U)$ for $a\in\sigma_j$ its boundary is 
\[
\partial_h \xi \;:=\; \partial\xi \,+\, (-1)^{i-j} \sum_{b\in\sigma_{j-1}} D_{a b}\, j_b(\xi).
\]
Here $\partial\xi\in C_{i-1-j}(U_a)$ is the ordinary boundary of $\xi$ and $j_b(\xi)$ denotes the same chain as $\xi$, but considered as an element of $C_{i-j}(U_b)$ if $U_a\subset U_b$. The definition of the boundary map is then extended to $C_i(\U)$ by linearity. 

The augmentation morphism $\veps:C(\U)\rightarrow C(X)$ sends all chains in $C_{i-j}(U_a)$ with $j>0$ to $0$ and the ones with $j=0$ to themselves.

We have the following lifting property for lifting chains to hyperchains. Here the word \emph{subdivision} of a chain or of a hyperchain means some iterated barycentric subdivision of all of its simplices. It is clear that the operation of taking subdivision commutes with the boundary operation.
\begin{lem}
Let $\xi$ be a chain on $X$, $\eta=\partial \xi$, $\bar\eta$ a hyperchain such that $\veps\bar\eta=\eta$ and $\partial_h \bar\eta=0$. Then, after possibly replacing $\eta$, $\bar\eta$, $\xi$ with a subdivision, there exists a hyperchain $\bar\xi$ such that $\partial_h\bar\xi=\bar\eta$ and $\veps \bar\xi=\xi$.
\end{lem}
\[
\begin{CD}
\bar\xi\in C_i(\U) @>{\partial_h}>> \bar\eta\in C_{i-1}(\U) \\
@V{\veps}VV @V{\veps}VV \\
\xi\in C_i(X) @>{\partial}>> \eta\in C_{i-1}(X)
\end{CD}
\]
\begin{proof}
The complex $C_\bullet(\U)$ is the total complex of the bicomplex $C_\bullet(\U,\bullet)$ with the horizontal differential induced by the boundary map of chains in space $X$ and vertical one induced by the boundary map of the complex $\sigma$:
\[
\begin{CD}
& & C_{i-1}(\U, 1) @>>> C_{i-2}(\U,1) \\
& & @VVV @VVV\\
C_i(\U,0) @>>> C_{i-1}(\U,0) @>>> C_{i-2}(\U,0)\\
@V{\veps}VV @V{\veps}VV\\
C_i(X) @>>> C_{i-1}(X)
\end{CD}
\]
Clearly, it is enough to prove that the vertical complexes are exact (up to a subdivision). That is, we need to prove that if 
\[
\gamma\eq\sum_{a\in\sigma_{j}} \gamma_a \;\in\; C_i(\U,j) \quad \text{is such that}
\]
\[
\sum_{a\in\sigma_{j}} D_{ab} \gamma_a\eq0 \quad \text{for any $b\in\sigma_{j-1}$,}
\]
then there exists 
\[
\bar\gamma\eq\sum_{c\in\sigma_{j+1}} \bar\gamma_c \in C_i(\U,j+1) \quad \text{such that}
\]
\[
\sum_{c\in\sigma_{j+1}} D_{ca} \bar\gamma_c \eq \gamma'_a \quad \text{for any $a\in\sigma_{j}$,}
\]
where $\gamma'_a$ is a subdivision of $\gamma_a$.

It is enough to prove the statement for multiples of a single simplex. Let $s$ be a simplex in $X$ which enters $\gamma_a$ with coefficient $s_a\in\Z$. Then 
\[
\sum_{a\in\sigma_{j}} D_{ab} s_a\eq 0 \quad \text{for any $b\in\sigma_{j-1}$.}
\]
Therefore the cycle $\sum_a s_a a$ of $\sigma$ is closed. Since the simplex $s$ must belong to $U_a$ for all $a$ for which $s_a\neq 0$, the mentioned cycle is a closed cycle of $\sigma_x$ for all $x$ in the closure of $s$. For each $x$ one can therefore represent it as a boundary of some cycle of $\sigma_x$, say $t_x$,
\[
t_x \eq \sum_{c\in\sigma_{j+1}} t_{x c} c.
\]
Here $t_{x c}$ is zero unless $x\in U_c$. Let $V_x$ be the open set which is the intersection of all $U_c$ for which $t_{xc}$ is not zero. Then $V_x$ form a cover of the support of $s$, so there exists a finite subcover. Let it be $V_1=V_{x_1}$, $V_2=V_{x_2}$, \dots, $V_k=V_{x_k}$. One can subdivide the simplex $s$ so that each simplex of the subdivision belongs to one of the chosen open sets. Let us do this and denote the subdivision by $s'$ so that
\[
s' \eq \sum_{l=1}^k s'_l, \qquad |s'_l| \subset V_k.
\]
Now it is clear that one may put
\[
\bar\gamma_l \eq \sum_{c\in\sigma_{j+1}} t_{x_l c} j_c(s'_l)
\]
because $s'_l$ belongs to $U_c$ for every $c$ for which $t_{x_l c}$ is not zero. Then 
\[
\sum_{c\in\sigma_{j+1}} D_{ca} t_{x_l c} s'_l \eq s_a s'_l \qquad \text{for any $a\in\sigma_{j}$,}
\]
therefore $\bar\gamma:=\sum_{l=1}^k \bar\gamma_l$ satisfies the required condition.
\end{proof}

This immediately implies
\begin{cor}
If $\xi$ is a closed chain on $X$ then there exists a closed hyperchain $\bar\xi$ whose augmentation is a subdivision of $\xi$. The hyperchain $\bar\xi$ is unique up to a boundary of a hyperchain whose augmentation is $0$.
\end{cor}

\begin{rem}
One extends the construction of hyperchains to the cases of chains on an open subset of $X$ and to the relative chains. It is clear that the boundary maps commute with the natural projections from chains to relative chains.
\end{rem}

\subsection{Integration}
We again fix a space $X$, an abstract cell complex $\sigma$ and a hypercover $(U_a)$. Suppose some class of chains and some class of differential forms are given such that one can integrate a form along a chain. We have the de Rham complex of sheaves of differential forms $\Omega^\bullet$:
\[
\Omega^0\ra\Omega^1\ra\dots.
\]
Also we have the complex of chains on $X$:
\[
\dots\ra C_1(X)\ra C_0(X).
\]

The hypersections of $\Omega^\bullet$ will be called the \emph{hyperforms}. Given a hyperform $\omega$ of degree $d$ and a hyperchain $\xi$ of degree $d$ we define the \emph{integral}:
\[
\int_\xi \omega \eq \sum_{a\in\sigma} \int_{\xi_a}\omega_a.
\]
Here $\omega_a \in\Gamma(U_a, \Omega^{d-\dim a})$ and $\xi_a\in C_{d-\dim a}(U_a)$ are the components of $\omega$ and $\xi$.

The main property is the Stokes theorem:
\begin{prop}
For a hyperform $\omega$ of degree $d$ and a hyperchain $\xi$ of degree $d+1$ one has
\[
\int_{\partial_h \xi} \omega \eq \int_{\xi} d\omega.
\]
\end{prop}
\begin{proof}
By the definition
\[
\int_{\xi} d\omega \eq \sum_{a\in\sigma} \int_{\xi_a} \left(d \omega_a + (-1)^{d-\dim a+1}\sum_{b\in\sigma_{\dim a - 1}} D_{ab} \omega_b \right).
\]
Applying the Stokes formula this is further equal to 
\[
\sum_{a\in\sigma}\left(\int_{\partial\xi_a}\omega_a + (-1)^{d-\dim a+1}\sum_{b\in\sigma_{\dim a - 1}} D_{ab} \int_{\xi_a} \omega_b \right) \eq \sum_{a\in\sigma} \int_{\partial \xi_a} \omega \eq \int_{\partial_h\xi} \omega.
\]
\end{proof}

\subsection{Hodge filtration}\label{filtration}
The Hodge filtration on hyperforms is induced by the Hodge filtration on the de Rham complex:
\[
F^i \Omega^\bullet(\U) \eq \Omega^{\bullet\geq i}(\U),
\]
so a hyperform belongs to $F^i$ if all of its components which are forms of degree less than $i$ are zero.

\subsection{Products}\label{products}
For two spaces $X$, $X'$, two abstract complexes $\sigma$, $\sigma'$, two hypercovers $(U_a)$, $U_{a'}'$ we have the product hypercover. There is exterior product on chains and on sheaves. For example, consider the exterior product on chains. For two chains $c\in C_i(X)$, $c'\in C_{i'}(X)$ we obtain the product chain $c\times c'\in C_{i+i'}(X\times X')$. This operation satisfies the property
\[
\partial(c\times c')\eq(\partial c) \times c' \,+\, (-1)^i c\times \partial c'.
\]

Introduce an exterior product on hyperchains. Let $c\in C_{i-j}(U_a)\subset C_i(\U)$, $c'\in C_{i'-j'}(U'_{a'})\subset C_{i'}(\U')$. Here $\dim a = j$, $\dim a'=j'$. Then put
\[
c\times_h c' \;:=\; (-1)^{j(i'-j')} c\times c' \;\in\; C_{i+i'-j-j'}(U_a\times U'_{a'}) \;\subset\; C_{i+i'}(\U\times\U').
\]
The superscript $h$ stands for ``hyper''.

One can check that this satisfies
\begin{prop}
\[
\partial_h(c\times_h c') \eq (\partial_h c) \times_h c' \,+\, (-1)^i c\times_h \partial_h c'.
\]
\end{prop}

On hypersections the exterior product is defined in a similar way. If $s=(s_a)\in \Gamma(\U, \F^\bullet)^i$ and $s'=(s'_{a'})\in \Gamma(\U', \F'_\bullet)^{i'}$, then
\[
(s\times s')_{a\times a'} \eq (-1)^{\dim a (i'-\dim a')} s_a\times s'_{a'}.
\]
One has a similar formula
\begin{prop}
\[
d(s\times s') \eq (ds) \times s' \,+\, (-1)^i s\times ds'.
\]
\end{prop}

\subsection{Residues}\label{residues}
Let $X$ be a projective algebraic variety over $\C$, $\sigma$ an abstract cell complex and $\U=(U_a)$ a hypercover of $X$ in the Zariski topology. 

\begin{defn}
By a \emph{refinement} of $\U$ we understand any hypercover $\U'=(U_a')$ of $X$ in the analytic topology such that $U_a'\subset U_a$ for any $a\in \sigma$.
\end{defn}

Suppose a finite family $M$ of irreducible subvarieties of $X$ is given. These subvarieties will be called \emph{special}.

\begin{defn}
A refinement $\U'=(U_a')$ is called \emph{nice} if for any cell $a\in \sigma$ and any special subvariety $Z$ such that $\dim Z < \dim a$ one has $U_a'\cap Z = \varnothing$. 
\end{defn}

If we have a nice refinement $\U'=(U_a')$ we can make the following construction. Fix a special subvariety $Z$. Take the fundamental class of $Z$, represent it by a chain (of dimension $2\dim Z$) and lift it to a closed hyperchain $(c_a)$ with respect to the hypercover $(U_a' \cap Z)$. Note that for $a\in\sigma_{\dim Z+1}$ the set $U_a'\cap Z$ is empty. Therefore for every $a\in \sigma_{\dim Z}$ the chain $c_a$ is a closed chain of dimension $\dim Z$. Moreover, if we choose a different representation of the fundamental class or a different lift, the hyperchain will differ by a boundary. This implies that the chain $c_a$ for $a\in\sigma_{\dim Z}$ will change by a boundary. Thus the class of $c_a$ in the homology of $Z\cap U_a'$ does not depend on the choices made. Denote this class by $h_a(Z, \U')\in H_{\dim Z}(Z\cap U_a')$. 

\begin{defn}
For any meromorphic differential form $\omega$ on $Z\cap U_a$ of degree $\dim Z$ which is regular outside the special subvarieties of $Z$, for a cell $a$ of dimension $\dim Z$, we define its residue as
\[
\res_{a,Z,\U'}^{\Int} \omega \eq \int_{h_a(Z, \U')} \omega.
\]
\end{defn}

Correspondingly, if we have a meromorphic hyperform of degree $2 \dim Z$ on $Z$, then it is just a family of meromorphic forms of degree $\dim Z$ indexed by the cells of dimension $\dim Z$. Thus a choice of the refinement gives a trace map:
\begin{defn}
For any meromorphic hyperform $\omega$ on $Z$ of degree $2\dim Z$ which is regular outside the special subvarieties of $Z$ we define its trace as
\[
\Tr_{Z,\U'}^{\Int} \omega \eq \sum_{a\in \sigma_{\dim Z}} \res_{a,Z,\U'}^{\Int} \omega_a.
\]
\end{defn}

In the next section we construct a nice refinement for the case when $X$ is a product of curves and a finite Zariski cover is chosen for each of these curves. So we will omit the subscript $\U'$.

\subsection{Construction of nice refinements}\label{construction_refinements}
We will construct some nice refinements in certain special situation and relate the corresponding residues with the ordinary residues.

Let $X$ be a product of smooth projective curves, $X=X_1\times\cdots\times X_n$. Suppose on each curve a finite open cover (in the Zariski topology) is given. As we have seen before, this gives a hypercover of $X_k$ indexed by the standard simplex of dimension, say $m_k$, the \v{C}ech hypercover. By taking products we obtain an abstract cell complex
\[
\sigma \eq \Delta_{m_1}\times \Delta_{m_2}\times \dotsb \times\Delta_{m_n}
\]
and the product hypercover $\U$ of $X$ indexed by $\sigma$. For any finite family $M$ of irreducible subvarieties of $X$ we construct a nice refinement of $\U$. Here we describe the construction.

Without loss of generality we may assume that the set of special subvarieties $M$ satisfies the following conditions:
\begin{enumerate}
\item For any open set $U_a$ belonging to the hypercover the irreducible components of its complement are in $M$.
\item For any two sets in $M$ the irreducible components of their intersection is also in $M$.
\item If $M$ contains a set $Z$, then it also contains the singular locus of $Z$.
\item If $M$ contains a set $Z$, then for any subset $L\subset\{1,\ldots,n\}$ $M$ contains the irreducible components of the set where the projection from $Z$ to the product $\times_{k\in L} X_k$ is not \'etale (this may be the whole set $Z$).
\end{enumerate}

Let $X_k$ be one of the curves above. Choose a Riemannian metric on $X_k$. Let $n(S,\veps)$ denote the $\veps$-neighborhood of a set $S$. For any $k$ we give a refinement of the hypercover of the curve $X_k$. This depends on two real numbers $\veps>\veps'>0$. 

Suppose $X_k$ is covered by the open sets $U_{k,0},\ldots,U_{k,m_k}$. Denote the complement $X_k\setminus U_{k,0}$ by $S_k$. The set $S_k$ is a finite set of points $p_1,p_2,\ldots,p_r$. 

Consider those points among $p_1,p_2,\ldots,p_r$ which are covered by $U_{k,1}$. Suppose they are $p_1,p_2,\ldots,p_{r_1}$. Then consider those points among the remaining ones which are covered by $U_{k,2}$. Suppose they are $p_{r_1+1},\ldots,p_{r_2}$, etc. In this way we obtain a decomposition of the set $S_k$ into $m_k$ subsets, some of them empty:
\[
S_k \eq \bigcup_{i=1}^{m_k} S_{k,i}.
\]
Let $\wt{S_k} = S_k\cup \{\eta_k\}$, where $\eta_k$ is the generic point of $X_k$. Let $S_{k,0} = \{\eta_k\}$.
For $p\in\wt{S_k}$ put 
\[
U^p_k(\veps, \veps') \eq \begin{cases} n(p, \veps) & \text{if $p$ is a closed point, $p\in S_k$,}\\ X_k\setminus \overline{n(S_k, \veps')} & \text{if $p$ is the generic point.}\end{cases}
\]
For $p\in S_k$ put
\begin{align*}
R^p_k(\veps, \veps') \eq U^p_k(\veps, \veps') \cap U^{\eta_k}_k(\veps, \veps'),\qquad R_k(\veps, \veps') &\eq \bigcup_{p\in S_k} R^p_k(\veps,\veps'), 
\\ U^S_k(\veps,\veps') &\eq \bigcup_{p\in S_k} U^p_k(\veps,\veps').
\end{align*}

We define the refinement in the following way. Put 
\[
U_{k,i}'(\veps, \veps') \eq \bigcup_{p\in S_{k,i}} U^p_k(\veps, \veps'),\qquad i=1,\dotsc,m_k.
\]
This is a cover of $X_k$. We obtain $\U_k'(\veps,\veps')$ as the \v{C}ech hypercover associated to this cover. If $\veps$ is small enough, this is a refinement of the original cover, i.e. $U_{k,i}'\subset U_{k,i}$. Moreover, if $\veps$ is small enough, the open sets $U^p_k$ for $p\in S_k$ are non-intersecting disks. Let $q>0$ be a number such that as soon as $\veps<q$ these two conditions are satisfied.

For $p\in \wt{S_k}$ we define by $a_k(p)$ the $0$-cell $a$ such that $p\in S_{k,a}$. For $p\in S_k$ we define by $a_k(\eta_k,p)$ the $1$-cell $a$ which joins $a_k(\eta_k)$ and $a_k(p)$. Then for any cell $a$ of the standard simplex $\Delta_{m_k}$ the elements of the hypercover are given as follows:
\[
U_{k,a}'(\veps, \veps') \eq \begin{cases} 
  \bigcup_{p\in \wt{S_k},a_k(p)=a} U^p_k(\veps, \veps') & \text{if $\dim a=0$,}\\
  \bigcup_{p\in S_k, a_k(p, \eta_k) =a} R^p_k(\veps, \veps') & \text{if $\dim a=1$,}\\
  \varnothing & \text{otherwise.}
\end{cases}
\]

We consider vectors of real numbers $\vec\veps=(\veps_1,\ldots,\veps_n, \veps_1',\ldots,\veps_n')$ such that for each $k$, $q>\veps_k>\veps_k'>0$. Denote
\[
\U(\vec\veps) \eq \times_{k=1}^n \U_k(\veps_k,\veps_k').
\]
This is a refinement of $\U$. We denote for any $L\subset\{1,\dots,n\}$ 
\[
R_L(\vec\veps) \eq \times_{k\in L} R_k(\veps_k,\veps_k'),\qquad U^S_L(\vec\veps) \eq \times_{k\in L} U^S_k(\veps_k,\veps_k'),\qquad S_L\eq\times_{k\in L} S_k.
\]
For any $p\in S_L$ we denote $p_k=\pi_k p$ and put
\[
R_L^p(\vec\veps) \eq \times_{k\in L} R_k^{p_k}(\veps_k, \veps_k').
\]
We will simply write $R_L$, $U_L^S$, $R_L^p$ when there is no confusion.
We put 
\[
X_L\eq\times_{k\in L} X_k,\qquad \pi_L:X\ra X_L\quad \text{the projection.}
\]

\begin{defn}
We say that the choice of real numbers $\veps_m,\ldots,\veps_n, \veps_m',\ldots,\veps_n'$ is \emph{good} if $q>\veps_k>\veps_k'>0$ for $k=m,\ldots,n$ and there is an increasing sequence of positive real numbers 
\[
\delta_m, \delta_m',\delta_{m+1},\delta_{m+1}',\dotsc,\delta_n,\delta_n'
\]
such that the following conditions are satisfied for any special set $Z$, an index $k$ ($m\leq k\leq n$) and $p\in S_k$:
\begin{enumerate}
\item If $x\in Z$ is such that $\dist(\pi_k(x),p)\leq\veps_k$, then $\dist(x, Z')<\delta_k$.
\item If $x\in X$ is such that $\dist(x, Z)<\delta_{k-1}'$ and $\dist(\pi_k(x), p)<\veps_k$, then $\dist(x, Z')<\delta_k$.
\item If $x\in X$ is such that $\dist(x, Z')<\delta_{k-1}'$, then $\dist(\pi_k(x), p)\leq\veps_k'$.
\end{enumerate}
And the following condition is satisfied for any two special sets $Z_1$, $Z_2$ and an index $k$, $m\leq k\leq n$: 
\begin{enumerate}
\item[(iv)] If $x\in X$ is such that $\dist(x, Z_1)<\delta_k$ and $\dist(x, Z_2)<\delta_k$, then $\dist(x, Z_1\cap Z_2)<\delta_k'$.
\end{enumerate}
We have denoted by $\pi_k$ the projection $X\rightarrow X_k$ and by $Z'$ the intersection $Z\cap\pi_k^{-1} p$. The second and the third conditions are required only for $k>m$.
\end{defn}

Next we prove that good choices indeed exist. For this it is enough to show
\begin{lem}\label{lem:good_choices}
For every good choice of real numbers $\veps_{m+1},\ldots,\veps_n, \veps_{m+1}',\ldots,\veps_n'$ there exists a number $t>0$ such that for any $\veps_m, \veps_m'$ satisfying $t>\veps_m>\veps_m'>0$ the choice $\veps_{m},\ldots,\veps_n$, $\veps_{m}',\ldots,\veps_n'$ is good. Therefore good sequences exist.
\end{lem}
\begin{proof}
Suppose there is a sequence $\delta_{m+1},\ldots,\delta_n'$ as in the definition above. Consider $Z\in M$ and $Z'=Z\cap \pi_{m+1}^{-1} S_{m+1}$. We first show that we can choose $\delta_m'$ such that the second and the third conditions are satisfied for $k=m+1$.

Consider the set 
\[
V\eq\{x\in X: \dist(x, Z')<\delta_k\quad \text{or}\quad \dist(\pi_{m+1}(x), S_{m+1})> \veps_{m+1}\}.
\]
This is an open neighborhood of $Z$ by the first condition for $k=m+1$. Therefore if $\delta_m'>0$ is small enough, then the $\delta_m'$-neighborhood of $Z$ is also contained in $V$. Therefore the second condition is satisfied for $k=m+1$.

Since $Z'$ is a compact set and $\pi_{m+1}(Z')\subset S_{m+1}$, if $\delta_m'>0$ is small enough, the third condition is satisfied.

Therefore we can choose $\delta_m'>0$ such that $\delta_m'<\delta_{m+1}$ and both the second and the third conditions are satisfied for all $Z$ and $k=m+1$.

Next we choose $\delta_m$ such that the fourth condition is satisfied. To see that this can be done for $Z_1$, $Z_2$ we consider the compact set $X\setminus n(Z_1\cap Z_2, \delta_m')$. The two continuous functions $\dist(\bullet, Z_1)$ and $\dist(\bullet, Z_2)$ do not attain simultaneously value zero on this set. Therefore if $\delta_m$ is small enough, these functions cannot attain simultaneously value less than $\delta_m$. This is equivalent to the fourth condition.

Next consider $Z\in M$ and $Z'=Z\cap \pi_m^{-1} S_m$. The set
\[
\{x\in Z: \dist(x,Z')\geq \delta_m\} 
\]
is compact, hence the continuous function $d(x)=\dist(\pi_m(x), S_m)$ attains its minimum. Suppose $\epsilon_m$ is less than the minimal value of this function. It follows that if $x\in Z$ and $d(x)\leq \epsilon_m$, then $x$ cannot belong to the set above. Therefore $\dist(x,Z')< \delta_m$ and the first condition is satisfied. This implies existence of $t>0$ with the required properties.
\end{proof}

\subsection{Flags of subvarieties}\label{sec:flags_subvar}
Recall that $X=X_1\times \cdots \times X_n$, each $X_k$ is covered by a finite number of open sets $U_{k,0},\ldots,U_{k,m_k}$, $S_k$ is the complement to $U_{k,0}$, $\wt{S_k} = S_k\cup \{\eta_k\}$, where $\eta_k$ is the generic point of $X_k$. For a set $L\subset\{1,\ldots,n\}$ we denote $S_L = \times_{k\in L} S_k$. Recall that $a_k(p)$ is defined as the minimal number such that $p\in U_{k, a_k(p)}$ for $p\in S_k$ and $a_k(p)=0$ for $p=\eta_k$. The number $a_k(p)$ is viewed as a vertex of the $m_k$-dimensional simplex. Also we have $a_k(\eta_k, p)$ for $p\in S_k$ defined as the pair $a_k(\eta_k), a_k(p)$, which gives an edge of the $m_k$-dimensional simplex and corresponds to the open set $U_{k, 0} \cap U_{k, a_k(p)}$.

We consider flags of subvarieties of $X$. A flag of length $m$ is a sequence of irreducible subvarieties $Z_\bullet=(Z_0\supset Z_1 \supset \cdots \supset Z_m)$. We say that a flag $Z_\bullet$ \emph{starts} with $Z_0$ and \emph{ends} with $Z_m$. We require $Z_m$ to be not empty.

\begin{defn}
Let $L\subset\{1,\ldots,n\}$, $L=\{k_1,\ldots,k_l\}$. Let $p\in S_L$. We say that a flag $Z_\bullet=Z_0\supset\cdots\supset Z_l$ is $L,p$\emph{-special} if
\begin{enumerate}
\item $Z_i$ is special for $0\leq i\leq l$,
\item $Z_i$ is an irreducible component of $Z_{i-1}\cap \pi_{k_i}^{-1} p_{k_i}$ for $1\leq i\leq l$.
\end{enumerate}
\end{defn}

\begin{defn}
We say that a flag $Z_\bullet=Z_0\supset\cdots\supset Z_l$ is \emph{strict at index $i$} if $Z_i\neq Z_{i-1}$. We say that a flag is \emph{strict} if it is strict at all indices.
\end{defn}

For a $d$-dimensional irreducible subvariety $Z$ of $X$ and a subset $L\subset\{1,\ldots,n\}$ of size $d$ the finite set of all strict $L,p$-special flags starting with $Z$ for all $p$ is denoted $\Fl_L(Z)$.

If $Z_\bullet\in \Fl_L(Z)$, $L=\{k_1, \ldots, k_d\}$, we construct a $\dim Z$-dimensional cell $a_L(Z_\bullet)$ of $\sigma$ (recall that $\sigma$ is the product of simplices $\Delta_{m_1} \times \cdots \times \Delta_{m_n}$) in the following way. For each $k\in L$ let $i$ be such that $k_i=k$ and peek the edge $a_k(\eta_k, \pi_k(Z_i))$ (note that $\pi_k(Z_i)\in S_k$). For $k\notin L$ peek the vertex $a_k(p_k(Z_\bullet))$, where $p_k(Z_\bullet)\in \wt{S_k}$ is defined as $\pi_k(Z_i)$ if this is a point or $\eta_k$ otherwise, where $i$ is the maximal number such that $k_i<k$ or $0$ if $k_1> k$. The product of these edges and vertices gives a cell in $\sigma$, which we denote by $a_L(Z_\bullet)$.

Now we prove that the refinements constructed in the last section are nice.
Suppose $\vec\veps$ is good for a sequence of numbers $\vec\delta=(\delta_1,\delta_1',\ldots,\delta_n,\delta_n')$. The \emph{$\vec\delta$-neighborhood} of an $L$-special flag $Z_\bullet$ is defined to be the set
\[
n(Z_\bullet, \vec\delta) \eq \{x\in X: \dist(x, Z_i)<\delta_{k_i},\;i\geq 1\}.
\]

\begin{prop}\label{prop:sf0}
Suppose $Z\subset X$ is special, $L\subset\{1,\ldots,n\}$, $x$ is a point for which either $x\in Z$ or $\dist(x, Z)<\delta_{\min L -1}'$ (if $\min L>1$). Suppose $p\in S_L$ and $\pi_L(x)\in U_L^p(\vec\veps)$. Then there exists an $L,p$-special flag $Z_\bullet$ starting with $Z$ such that $x$ belongs to its $\vec\delta$-neighborhood.
\end{prop}
\begin{proof}
We construct the flag step by step. If 
\[
\dist(x,Z_{i-1})<\delta_{k_{i-1}}'\leq\delta_{k_i-1}',\quad \dist(\pi_{k_{i}} x, p_{k_i})<\veps_{k_i},
\]
by the second property of good sequences we obtain 
\[
\dist(x,Z_{i-1}\cap \pi_{k_i}^{-1} p_{k_i})<\delta_{k_i}.
\]
Therefore on each step we can choose $Z_i$ as an irreducible component of $Z_{i-1}\cap \pi_{k_i}^{-1} p_{k_i}$ so that $\dist(x,Z_i)<\delta_{k_i}$.
\end{proof}

\begin{prop}\label{prop:sf}
Let $x$, $Z$, $L$ and $p$ be as in the proposition above. Take an $L,p$-special flag $Z_\bullet$ such that $x$ belongs to the $\vec\delta$-neighborhood of $Z_\bullet$. If $L'\subset L$ is such that $\pi_{L'}(x)\in R_{L'}(\vec\veps)$, then $Z_\bullet$ is strict at all indices $i$ for which $k_i\in L'$.
\end{prop}
\begin{proof}
If $k_i\in L'$, then $\dist(\pi_{k_i}, p_{k_i})>\veps_{k_i}'$. By the third property of good sequences we obtain $\dist(x, Z_{k_i})\geq \delta_{k_i}'>\delta_{k_{i-1}}'$. Since $\dist(x, Z_{k_{i-1}})<\delta_{k_{i-1}}$, $Z_{k_{i-1}}\neq Z_{k_i}$.
\end{proof}

\begin{cor}
If $\vec\veps$ is a good choice of numbers, then for any special subvariety $Z$ of dimension less than $d$ and $L\subset\{1,\ldots,n\}$ of size $d$ the intersection $Z\cap \pi_L^{-1} R_L(\vec\veps)$ is empty. Therefore the refinement $\U'(\vec\veps)$ is nice.
\end{cor}
\begin{proof}
If the intersection was not empty, by Proposition \ref{prop:sf} there would exist a strict $L,p$-special flag which starts with $Z$. Therefore the dimension of $Z$ would be at least $d$. 
\end{proof}

\subsection{Decomposition of the residue according to flags}
Let $X$, $\sigma$, $\U$, $M$ be as in the previous section. Let $Z$ be a subvariety of $X$, $Z\in M$. Let $d=\dim Z$ and $\omega$ be a meromorphic differential form on $Z$ of degree $d$ which is regular outside the special subvarieties of $Z$. Let $a\in\sigma_d$. Let $\U'=\U(\vec\veps)$ be the nice refinement of $\U$ corresponding to a good choice of numbers $\vec\veps$. We will show how to compute
the residue $\res_{a,\U'} \omega$.

Recall that the residue was defined as the integral of $\omega$ along the class $h_a(Z, \U')\in H_d(Z\cap U_a')$. Since 
\[
\sigma\eq\prod_{k=1}^n \Delta_{m_k},
\]
the cell $a$ is given by a sequence of cells $a_1\in\Delta_{m_1},\ldots,a_n\in\Delta_{m_n}$. The open set $U_a'$ is the product of sets $U_{k,a_k}(\veps_k, \veps_k')$ for $k=1,\ldots,n$. Note that for $k=1,\ldots,n$ the set $U_{k,a_k}(\veps_k, \veps_k')$ is not empty only if either $a_k$ is a point or $a_k$ is the edge joining the $0$ vertex and some other vertex. Therefore we get a non-zero residue only if there is a set $L\subset\{1,\ldots,n\}$ of size $d$ and 
\[
\dim a_k\eq \begin{cases} 1 & \text{for $k\in L$,}\\ 0 & \text{otherwise.}\end{cases}
\]
Moreover in this case we have $U_a'\subset \pi_L^{-1} R_L$. 

We have seen that any point of $Z\cap \pi_L^{-1} R_L$ belongs to the $\vec\delta$-neighborhood of a flag from $\Fl_L(Z)$. In fact we have
\begin{prop}\label{prop:unfl}
Any point $x\in Z\cap \pi_L^{-1} R_L$ belongs to the $\vec\delta$-neighborhood of a unique flag from $\Fl_L(Z)$. 
\end{prop}
\begin{proof}
Suppose $x$ belongs to the $\vec\delta$-neighborhood of two flags $Z_\bullet$ and $Z_\bullet'$. Let $i$ be the minimal index for which $Z_i\neq Z_i'$. Since $\dist(x, Z_i)<\delta_{k_i}$ and $\dist(x, Z_i')<\delta_{k_i}$, by the fourth property of good sequences we have $\dist(x, Z_i\cap Z_i')<\delta_{k_i}'$. Let $Y$ be an irreducible component of $Z_i\cap Z_i'$ such that $\dist(x, Y)<\delta_{k_i}'$. Let $K$ be the subset of elements of $L$ which are greater than $k_i$. By Proposition \ref{prop:sf} one can construct a strict $K$-special flag starting from $Y$ whose $\vec\delta$-neighborhood contains $x$. But one can see that the length of the flag equals to the dimension of $Z_i$ and is at least one more than the dimension of $Y$. Hence such flag does not exist and we obtain a contradiction.
\end{proof}

This gives us a possibility to decompose $Z\cap \pi_L^{-1} R_L$ according to the set $\Fl_L(Z)$.
\begin{cor}\label{cor:unfl}
The set $Z\cap \pi_L^{-1} R_L$ is the union of non-intersecting open sets $R_{L,Z_\bullet} = Z\cap \pi_L^{-1} R_L \cap n(Z_\bullet, \vec\delta)$, one for each $Z_\bullet\in \Fl_L(Z)$. Correspondingly we obtain the decomposition
\[
\res_{a}^{\Int}\omega \eq \sum_{Z_\bullet \in \Fl_L(Z)} \res_{a,L,Z_\bullet}^{\Int}\omega,\;\text{where}\;
\res_{a,L,Z_\bullet}^{\Int}\omega \eq \int_{h_a(Z,\U')\cap R_{L,Z_\bullet}}\omega.
\]
\end{cor}

On the other hand, the set $Z\cap \pi_L^{-1} R_L$ is a union of non-intersecting open sets $U_a'$ for $a$ running over the products of cells $\times_{k=1}^n a_k$ with $\dim a_k=1$ for $k\in L$ and $\dim a_k=0$ for $k\notin L$. It appears that each $R_{L,Z_\bullet}$ is contained in exactly one such $U_a'$ and we will see that $a=a_L(Z_\bullet)$.

The last variety in the flag, $Z_d$, is a point from $S_L$. 
\begin{prop}
Suppose $x\in R_L$. If $x\in R_{L,Z_\bullet}$ and $\pi_L(x)\in R_L^p$, then $p=Z_d$.
\end{prop}
\begin{proof}
Since $x\in Z\cap\pi_L^{-1}R_L^p$, one has an $L,p$-flag whose $\vec\delta$-neighborhood contains $x$. By Proposition \ref{prop:unfl} this flag must be $Z_\bullet$. Hence $Z_\bullet$ is $L,p$-special. Thus $Z_d=p$.
\end{proof}

If $x\in Z\cap R_L$ then $\pi_k(x)$ for any $k\in L^c$ belongs to $U_k^{p_k'}(\veps,\veps')$ for exactly one $p_k'\in \wt{S_k}$. If it was not true, then $x$ would belong to $R_{L\cup\{k\}}$, which is a contradiction. This defines a point $p_k'(x)\in\wt{S_k}$ for each $k\in L^c$.
\begin{prop}
For any $x\in Z\cap R_L$ if $x\in R_{L,Z_\bullet}$, then $p_k(Z_\bullet)=p_k'(x)$ for all $k\in L^c$.
\end{prop}
\begin{proof}
Let $K$ be the set of such $k$ that either $k\in L$ or $p_k'\in S_k$ (i.e. $p_k'$ is not the generic point). For $k\in L$ let $p_k'$ be such that $\dist(\pi_k(x), p_k')<\veps_k$. This defines a point $p'\in S_K$. Then for every $k\in K$ $\dist(\pi_k(x), p_k')<\veps_k$. Thus $\pi_K(x)\in U_L^{p'}$. 

By Proposition \ref{prop:sf0} there exists a $K, p'$-special flag $Z'$ starting from $Z$ whose $\vec\delta$-neighborhood contains $x$. This flag must be strict for all indices which correspond to elements of $L$. Since the number of this indices equals the dimension of $Z$ the flag must be not strict at all other indices. Therefore this flag defines a $L,p'$-special flag, which must be $Z$ by Proposition \ref{prop:unfl}. At the same time this shows that $p_k=p_k'$ for all $k\in K$. 

Let $k\notin K$. Then $\dist(\pi_k(x), S_k)>\veps_k'$. Suppose $i$ is the maximal index for which $k_i<k$ or $0$ if $k< k_1$. If $i\neq 0$, by the third property of good sequences 
\[
\dist(x, Z_{k_i}\cap\pi_k^{-1} S_k)\;\geq\; \delta_{k-1}'\;\geq\; \delta_{k_i}\;\geq\; \dist(x, Z_{k_i}).
\]
Therefore $Z_{k_i}\neq Z_{k_i}\cap\pi_k^{-1} S_k$ which means that $p_k=\eta_k$. The case $i=0$ is obvious.
\end{proof}

We see that 
\begin{cor}
For any cell $a$ of dimension $d$ the set $Z\cap U_a'$ is the disjoint union of open sets $R_{L,Z_\bullet}$ where $Z_\bullet$ runs over all $Z_\bullet\in \Fl_L(Z)$ with $a_L(Z_\bullet)=a$.
\end{cor}
\begin{rem}
Therefore we may omit $a$ in the notation $\res_{a,L,Z_\bullet}^{\Int}$.
\end{rem}

\subsection{Relation with iterated residues}\label{comp_res}
Note that $Z$ and $X_L$ have the same dimension. 
\begin{prop}
If the restriction of $\pi_L$ to $Z$ is not surjective, then the set $\Fl_L(Z)$ is empty.
\end{prop}
\begin{proof}
If not, then $\dim \pi_L(Z)<d$. Let $L=\{k_1,\ldots,k_d\}$. Let $Z_\bullet\in \Fl_L(Z)$. Consider the corresponding flag of irreducible subvarieties of $X_L$:
\[
\pi_L(Z)\eq\pi_L(Z_0)\supset \pi_L(Z_{1}) \supset \cdots \supset \pi_L(Z_{d}).
\]
Because of the dimension reasoning there must be an index $i$ with $\pi_L(Z_i)=\pi_L(Z_{i-1})$. This implies
\[
\pi_{k_i}(Z_{i-1}) \eq \pi_{k_i}(Z_i) \;\subset\; S_{k_i}.
\]
Therefore $Z_i=Z_{i-1}$, so the flag is not strict, which is a contradiction.
\end{proof}

We may therefore suppose without loss of generality that $\pi_L:Z\ra X_L$ is surjective. Hence the set of points on $Z$ where this map is not \'etale is a proper closed subvariety. The irreducible components of this subvariety have dimension $d-1$ or less and are special. Therefore $\pi_L^{-1} R_L\cap Z$ belongs to its complement. This means the following is true.
\begin{prop}
The projection
\[
\pi_L:\pi_L^{-1} R_L\cap Z\ra R_L
\]
is an unramified covering.
\end{prop}

Let $p=\pi_L Z_d\in S_L$.
Since for a flag $Z_\bullet\in \Fl_L(Z)$ the set $R_{L,Z_\bullet}$ is open and closed in $\pi_L^{-1} R_L\cap Z$ we get
\begin{prop}
The projection
\[
\pi_L:R_{L,Z_\bullet} \ra R_L^p
\]
is an unramified covering.
\end{prop}

\begin{rem}
By the proposition we see that since $R_L^p$ is a product of annuli, $R_{L,Z_\bullet}$ is a disjoint union of products of annuli.
\end{rem}

We are going to determine the cycle $h_a(Z,\U')\cap R_{L,Z_\bullet}\in H_d(R_{L,Z_\bullet})$. The $d$-th homology group of a product of $d$ annuli is $\Z$. Hence there is a canonical generator $h_c$ of $H_d(R_L^p)$. In fact $h_c$ can be defined as the product $c_1\times\cdots\times c_d$, where $c_i$ is the circle in $R_{k_i}^{p_{k_i}}$ going around $p_{k_i}$ counterclockwise. 

\begin{prop}
The cycle $h_a(Z,\U')\cap R_{L,Z_\bullet}$ is the pullback of $(-1)^{\frac{d(d-1)}2}h_c$ via the projection $\pi_L:R_{L,Z_\bullet} \ra R_L^p$.
\end{prop}
\begin{proof}
\begin{rem}
It is clear that we can decrease numbers $\veps_k$ and increase numbers $\veps_k'$. The sequence obtained in this way will be also good. Moreover the open subsets of the new hypercover are contained in the corresponding open subsets of the old one. Therefore the residues computed with respect to these hypercovers are equal. Therefore one can assume that the projection $\pi_{L}:R_{L,Z_\bullet} \ra R_{L}^{p}$
extends to an unramified covering for the closures of $R_{L,Z_\bullet}$ in $Z$, $R_{L}^{p}$ in $X_L$.
\end{rem}

Consider the commutative diagram:
\[
\begin{CD}
H_{2d}(X) @>>> H_{2d}(\ol{R_{L,Z_\bullet}}, \partial R_{L,Z_\bullet}) @<{\pi_L^*}<< H_{2d}(\ol{R_L^p},\partial R_L^p) \\
@V{h_a}VV @V{h_a}VV @V{h_a}VV\\
H_d(U_a') @>>> H_d(R_{L,Z_\bullet}) @<{\pi_L^*}<< H_d(R_L^p)
\end{CD}
\]
We see that it is enough to prove that the image of the fundamental class of $H_{2d}(\ol{R_L^p},\partial R_L^p)$ by the map $h_a$ in $H_d(R_L^p)$ is $(-1)^{\frac{d(d-1)}2}h_c$. There is a direct product decomposition $\ol{R_L^p}=\times_{k\in L} \ol{R_k^{p_k}}$. The hypercover on $\ol{R_L^p}$ is the product hypercover. For $k\in L$ let $\phi_k$ be the fundamental class of $H_2(\ol{R_L^p}, \partial R_L^p)$. Let us lift it to a hyperchain $\wt{\phi_k}$. 

The hypercover of $\ol{R_L^p}$ is the \v{C}ech hypercover associated to the cover with two open sets:
\begin{align*}
V_0 &\eq\ol{n(p_k, \veps_k)}\setminus \ol{n(p_k, \veps_k)}\;& \text{corresponding to the cell}\;& a(\eta_k),\\
V_1 &\eq n(p_k, \veps_k)\setminus n(p_k, \veps_k)\;& \text{corresponding to the cell}\;& a(p_k),\\
V_0\cap V_1 &\eq n(p_k, \veps_k)\setminus\ol{n(p_k, \veps_k)}\;& \text{corresponding to the cell}\;& a(\eta_k,p_k).
\end{align*}

Let $\veps'<r<\veps$. Consider topological chains
\begin{align*}
c^k_0 \eq \ol{n(p_k, \veps_k)}\setminus n(p_k,r)\in C_2(V_0),\;
c^k_1 &\eq \ol{n(p_k, r)}\setminus n(p_k,\veps_k')\in C_2(V_1),\\
c^k_{01} &\eq \partial n(p_k,r)\in C_1(V_0\cap V_1).
\end{align*}
They define a hyperchain $c^k$. We have $c_0+c_1=\phi_k$. The hyperchain is closed, therefore $h_{01}(\phi_k)=c^k_{01}$, which is the canonical generator of $H_1(R_k^{p_k})$.

Since the product of $\phi_k$ is the fundamental class of $H_{2d}(\ol{R_L^p},\partial R_L^p)$, the product of the hyperchains $c=\times_{k\in L} c^k$ lifts the fundamental class. The term of $c$ at $\times_{k\in L}(\eta_k, p_k)$ is, by the definition of the product for hyperchains, $(-1)^{\frac{d(d-1)}2} \times_{k\in L} c^k_{01}$. This is exactly $(-1)^{\frac{d(d-1)}2} h_c$.
\end{proof}

Let $Z_\bullet$ be a strict $L,p$-special flag, $L=\{k_1,\ldots,k_d\}$, $p=(p_{k_1},\ldots,p_{k_d})$. Let $Z'=Z_1$, $k=k_1$, $L'=\{k_2,\ldots,k_d\}$. Let $Z_\bullet'$ be the flag $Z_1\supset\cdots\supset Z_d$, $Z_\bullet'\in \Fl_{L'}(Z')$. Let $t_i$ be a local parameter on $X_{k_i}$ at the point $p_{k_i}=\pi_{k_i} Z_d$. Let $p'=\pi_{L'} p$. For $Z'$ we have the projection 
\[
\pi_{L'}:R_{L',Z_\bullet'} \ra R_{L'}^{p'},
\]
which is also an unramified covering. Let 
\[
\wt Z \eq Z'\times_{X_{L'}} Z \eq (X_{k_1}\times Z') \times_{X_L} Z. 
\]
We obtain the canonical diagrams
\[
\begin{CD}
\wt Z @>{\rho}>> Z  &\qquad & \wt Z @>{\rho}>> Z\\
@V{\tau}VV @VVV     @V{\wt\tau}VV @VVV     \\
Z' @>>> X_{L'}      & & (X_{k_1}\times Z') @>>> X_L
\end{CD}
\]
The diagonal embedding $Z'\ra Z'\times_{X_{L'}} Z'$ induces a morphism $\Delta:Z'\ra \wt Z$, which is a section to the natural projection $\tau:\wt Z \ra Z'$. Let 
\[
R_Z\eq Z \,\cap\, \pi_L^{-1}(U_k^{p_k}\times R_{L'}^{p'}) \,\cap\, n(Z_\bullet, \vec\delta).
\]
Consider the fiber products over $U_k^{p_k}\times R_{L'}^{p'}$:
\[
\begin{CD}
\wt R @>{\rho}>> R_Z \\
@V{\wt\tau}VV @V{\pi_L}VV     \\
U_k^{p_k}\times R_{L',Z_\bullet'} @>>> U_k^{p_k}\times R_{L'}^{p'}
\end{CD}
\]
Again, we have a section $\Delta:R_{L',Z_\bullet'}\ra \wt R$. Let $\wt R'$ be the union of connected components of $\wt R$ which intersect the image of $\Delta$. Let $\rho'$ be the restriction of $\rho$ to $\wt R'$.
\begin{rem}
The set $R_Z$ is open and closed in $Z\cap\pi_L^{-1}(U_k^{p_k}\times R_{L'}^{p'})$ because Proposition \ref{prop:unfl} and the first sentence of Corollary \ref{cor:unfl} work if we replace $R_L$ by $U_k^{p_k}\times R_{L'}^{p'}$.
\end{rem}
\begin{rem}\label{rem:RZ-nice}
All special subsets of $Z$ which intersect $R_Z$ are contained in $Z'$. Therefore the map $\pi_L$ on the diagram is an unramified covering outside $\{p_k\}\times R_{L'}^{p'}$ and $\pi_k^{-1} p_k \cap R_Z = R_{L',Z'_\bullet}$.
\end{rem}

\begin{prop}
The map $\rho'$ is an analytic isomorphism.
\end{prop}
\begin{proof}
Since $\rho'$ is a base change of the unramified covering $\pi_{L'}:R_{L',Z'_\bullet}\ra R_{L'}^{p'}$, it is an unramified covering itself. Therefore it is enough to construct a continuous section $s':R_Z\ra \wt R'$ to $\rho'$ which extends the diagonal map. This is equivalent to constructing a retraction $s:R_Z\ra R_{L',Z'_\bullet}$ which respects the projection $\pi_{L'}$.

Take a compact connected set $A\subset R_{L'}^{p'}$ such that $\pi_{L'}^{-1}(A)\cap R_{L',Z'_\bullet}=A_1\cup\cdots\cup A_m$ is a disjoint union of spaces isomorphic to $A$.

One can choose disjoint open subsets $V_1,\ldots, V_m$ in $\pi_{L'}^{-1}(A)\cap R_Z$ such that $A_i\subset V_i$. The set $C=\pi_{L'}^{-1}(A)\cap R_Z\setminus(V_1\cup\cdots\cup V_m)$ is closed. Therefore $\pi_L(C)$ is closed. Take $\alpha>0$ such that $n(p_k,\alpha)\times A$ does not intersect $\pi_{L}(C)$. This means that the open set $\pi_L^{-1}(n(p_k,\alpha)\times A) \cap R_Z$ is contained in the union of the sets $V_i$.

Let $V_i'= \pi_L^{-1}(n(p_k,\alpha)\times A) \cap V_i$. Let us prove that $V_i'$ is connected for each $i$. If not, then $V_i'=B_1\sqcup B_2$ with $B_j$ open, closed and nonempty. Since $A_i$ is connected, for some $j$ $B_j$ does not intersect $A_i$. Therefore $\pi_L B_j$ is open, closed and nonempty. Hence it must be the whole $U_k^{p_k}\times A$. This contradicts the assumption that $B_j$ does not intersect $A_i$.

One can construct a deformation retract retracting $\pi_{L'}^{-1}(A)\cap R_Z$ inside $\pi_L^{-1}(n(p_k,\alpha)\times A)\cap R_Z$. Therefore there are exactly $m$ connected components of $\pi_{L'}^{-1}(A)\cap R_Z$, each containing exactly one $A_i$. Therefore there is a unique map $s_A:\pi_{L'}^{-1}(A)\cap R_Z\ra \pi_{L'}^{-1}(A)\cap R_{L',Z'_\bullet}$ which is identity on $\pi_{L'}^{-1}(A)\cap R_{L',Z'_\bullet}$ and makes the diagram below commutative.
\[
\begin{CD}
\pi_{L'}^{-1}(A)\cap R_Z  @>>> U_k^{p_k}\times A\\
@V{s_A}VV @VVV\\
\pi_{L'}^{-1}(A)\cap R_{L',Z'_\bullet} @>>> A
\end{CD}
\]
Patching these $s_A$ together gives $s$ as required proving the first statement of the proposition.
\end{proof}

Take a meromorphic form $\omega$ on $Z$ which is holomorphic outside the special subvarieties of $Z$. The form $\omega$ can be written as
\[
\omega \eq f d t_1 \wedge d t_2 \wedge \cdots \wedge d t_n,
\]
where $f$ is a rational function on $Z$. Let $K$ be the field of fractions of $Z'$. Then $Z\times \spec K$ is a curve and $\Delta(\spec K)$ is a point. Therefore the algebraic residue $\res_{\Delta(\spec K)} \rho^* f d t_1$ is defined. We have
\begin{prop}
Put
\[
\omega' \eq (\res_{\Delta(\spec K)} \rho^* f d t_1) d t_2 \wedge \cdots \wedge d t_n.
\]
Then
\[
\res_{L,Z_\bullet}^{\Int} \omega \eq (-1)^{d-1} 2\pi\I \res_{L',Z_\bullet'}^{\Int} \omega'.
\]
\end{prop}
\begin{proof}
By the definition
\[
\res_{\U,L,Z_\bullet} \omega \eq \int_{h_a\cap R_{L,Z_\bullet}} \omega.
\]
Since $\rho'$ is an isomorphism,
\[
\int_{h_a\cap R_{L,Z_\bullet}} \omega \eq \int_{\rho'^*(h_a \cap R_{L,Z_\bullet})} \rho^* \omega.
\]
We have 
\begin{align*}
\rho'^* (h_a \cap R_{L,Z_\bullet}) 
&\eq (-1)^{\frac{d(d-1)}2}\wt R' \,\cap\, \rho^* \pi_L^* h_c  \\
&\eq (-1)^{\frac{d(d-1)}2}\wt R' \,\cap\, \wt{\tau^*}(h_{c k_1}\times (R_{L',Z'_\bullet}\cap\pi_{L'}^* h_c')),
\end{align*}
where $h_{c k_1}$ is the circle in $R_{k_1}^{p_{k_1}}$ and $h_c'$ is the product of circles in $R_{L'}^{p'}$.
By Fubini's theorem
\[
\int_{\wt R'\cap \wt{\tau^*}(h_{c k_1}\times (R_{L',Z'_\bullet}\cap\pi_{L'}^* h_c'))} \rho^* \omega \eq \int_{R_{L',Z'_\bullet}\cap\pi_{L'}^* h_c'} g(z') d t_2\wedge\cdots\wedge d t_d,
\]
where for $z'\in R_{L',Z'_\bullet}$
\[
g(z') \eq \int_{\wt R'\cap \tau^{-1} z' \cap \pi_{k_1}^{-1} h_{c k_1}} \rho^* f d t_1.
\]
The last integral is nothing else than 
\[
2\pi \I \res_{\Delta(z')} \rho^* f d t_1.
\]
Therefore
\[
g(z') d t_2\wedge\cdots\wedge d t_d \eq 2\pi\I\omega'
\]
Taking into account that
\[
R_{L',Z'_\bullet}\cap\pi_{L'}^* h_c' \eq (-1)^{\frac{(d-1)(d-2)}2}R_{L',Z'_\bullet}\cap h_{a',Z'},
\]
where $a'$ is the cell obtained from $a$ by replacing the component $a(\eta_{k_1},p_{k_1})$ with the component $a(p_{k_1})$, we obtain the statement.
\end{proof}

Let us denote by $\res_{L,Z_\bullet} \omega$ the iterated algebraic residue of $\omega$ with respect to the flag $Z_\bullet$. This is defined by induction on the dimension of $Z$ by the formula
\[
\res_{L,Z_\bullet}f d t_1\wedge\cdots\wedge d t_d \eq \res_{L',Z_\bullet'} (\res_{\Delta(\spec K)} \rho^* f d t_1) d t_2 \wedge \cdots \wedge d t_d.
\]
Thus we obtain a formula for our residue.
\begin{cor}\label{cor_3_1_22}
\[
\res_{L,Z_\bullet}^{\Int} \omega \eq (-1)^{\frac{d(d-1)}2}\,(2\pi\I)^d \res_{L,Z_\bullet} \omega.
\]
\end{cor}

For any subvariety $Z\subset X$ and any hyperform
using the decomposition of the residue according to flags we also can state a formula for the trace:
\begin{cor}\label{cor_3_1_23} For a meromorphic hyperform $\omega$ of degree $2 d$ on $Z\subset X_1\times\dotsb\times X_n$ which is holomorphic outside the special subvarieties of $Z$, $\dim Z = d$,
\[
\Tr_Z^{\Int} \omega = (-1)^{\frac{d(d-1)}2} \, (2\pi\I)^d \sum_{L\subset\{1,\dotsc,n\},\, |L|=\dim Z} \quad \sum_{Z_\bullet \in \Fl_L(Z)} \res_{L,Z_\bullet} \omega_{a_L(Z_\bullet)}.
\]
\end{cor}

Dropping the coefficient $(-1)^{\frac{d(d-1)}2} (2\pi\I)^d$ we may also define the algebraic version of the trace:
\[
\Tr_Z \omega \;:=\; \sum_{L\subset\{1,\dotsc,n\},\, |L|=\dim Z} \quad \sum_{Z_\bullet \in \Fl_L(Z)} \res_{L,Z_\bullet} \omega_{a_L(Z_\bullet)}.
\]
Then the latter corollary can be reformulated as follows:

\begin{cor} For a meromorphic hyperform $\omega$ of degree $2 d$ on $Z\subset X_1\times\dotsb\times X_n$ which is holomorphic outside the special subvarieties of $Z$, $\dim Z = d$,
\[
\Tr_Z^{\Int} \omega = (-1)^{\frac{d(d-1)}2} \, (2\pi\I)^d \Tr_Z \omega.
\]
\end{cor}

Now we make a summary of the computation of $\Tr_Z \omega$. The necessary notation is explained in Section \ref{sec:flags_subvar}.
\begin{enumerate}
 \item List all the subsets $L\subset\{1,\ldots,n\}$ of size $\dim Z$.
 \item For each $L$ list all the flags $Z_\bullet \in \Fl_L(Z)$.
 \item For each flag find the cell $a_L(Z_\bullet)$ of $\Delta_{m_1}\times\cdots\times \Delta_{m_n}$, which is a product of several vertices and edges (vertices for $k\in L^c$ and edges for $k\in L$).
 \item Compute the iterated residue corresponding to the flag $Z_\bullet$ of the form which is the component of $\omega$ on the open set $U_{a_L(Z_\bullet)}=\times_k U_k$, where $U_k$ is an open set from the cover of $X_k$ for $k\in L^c$, or an intersection of two open sets of the cover of $X_k$ for $k\in L$.
 \item Add these residues.
\end{enumerate}

\subsection{Gauss-Manin}\label{gauss-manin}
Suppose we have a morphism of algebraic varieties $X\ra S$. Let $\sigma$ be an abstract cell complex and $\U=(U_a)$ a hypercover on $X$ with respect to Zariski's topology. We assume that $S=\spec R$ is affine and all the open sets of the hypercover are affine. 

Associated to the de Rham complex on $X$ we have the complex of hypersections
\[
\Omega_X^0(\U)\ra\Omega_X^1(\U)\ra\cdots.
\]
Consider the complex of hypersections corresponding to the relative de Rham complex.
\[
\Omega_{X/S}^0(\U)\ra \Omega_{X/S}^1(\U)\ra\cdots.
\]
These are both complexes of $R$-modules.

\begin{prop}
The following natural sequence is exact for all $k\geq 0$:
\[
\Omega_S^1(S)\otimes_R\Omega_X^{k-1}(\U)\ra \Omega_X^k(\U)\ra \Omega_{X/S}^k(\U)\ra 0.
\]
\end{prop}
\begin{proof}
Since the corresponding sequence of sheaves is exact, it induces an exact sequence over any affine set.
\end{proof}

We define Gauss-Manin connection as follows. Let $\omega\in\Omega_{X/S}^k(\U)$ be closed. Lift it to $\ol\omega\in\Omega_{X}^k$. Then $d \ol\omega\in \kernel(\Omega_X^{k+1}(\U)\ra\Omega_{X/S}^{k+1}(\U))$. Choose $\ol\eta\in \Omega_S^1(S)\otimes_R\Omega_X^{k}(\U)$ which maps to $d \ol\omega$. Let $\eta$ be the projection of $\ol\eta$ in $\Omega_{X/S}^{k}(\U)\otimes_R\Omega_S^1(S)$. We say that $\eta$ is a Gauss-Manin derivative of $\omega$. Of course, this construction depends on several choices. Although $\eta$ is not well-defined, by abuse of notation we will write
\[
\eta\eq\nabla\omega
\]
if $\eta$ is a Gauss-Manin derivative of $\omega$.

Suppose we have a family of hyperchains $c_s\in C_k(\U) $, $s\in S$. This means that $c_s$ is a linear combination of simplices $c_s^i$ with each $c_s^i$ being a map from $\Delta\times S$ to $X$ which composed with $X\ra S$ gives the projection $\Delta\times S\ra S$. Here $\Delta$ denotes a standard simplex. We require the maps $c_s^i$ to be of the same type as we require for simplices, i.e. semi-algebraic. For $\omega\in\Omega_{X/S}^k(\U)$ consider the integral
\[
f(s)\eq\int_{c_s} \omega, \; s\in S.
\]
For a path $s_0 s_1$ in $S$ let $c_{s_0 s_1}\in C_{k+1}(\U)$ denote the corresponding hyperchain over this path. It provides a homotopy:
\[
\partial_h c_{s_0 s_1} \,+\, (\partial_h c)_{s_0 s_1}\eq c_{s_1} \,-\, c_{s_0}.
\]
\begin{lem}\label{lem:gauss-manin:2}
Let $\omega\in\Omega_{X/S}^k(\U)$ be a closed relative hyperform and $c=(c_s)_{s\in S}$ be a family of hyperchains, $c_s\in C_k(\U)$.
If $\eta=\nabla\omega$ with $\ol\eta$ and $\ol\omega$ as in the definition of $\nabla$, then we have
\[
d \int_{c_s} \omega \eq \int_{c_s} \nabla\omega + R(\partial_h c, \ol\omega).
\]
Here $R$ is the bilinear operator which is constructed as follows. The value of $R$ on a vector in $S$ represented by a path $s_t$ is
\[
\langle [s_t], R(\partial_h c, \ol\omega) \rangle \eq \frac{\partial}{\partial t}\bigg\vert_{t=0} \int_{(\partial_h c)_{s_0 s_t}} \ol\omega.
\]
\end{lem}
\begin{proof}
For a fixed $s_0\in S$ let $s_t$ be a path in $S$ starting from $s_0$. We have
\[
f(s_t)-f(s_0)\eq\int_{\partial_h c_{s_0 s_t}+(\partial_h c)_{s_0 s_t}} \ol\omega.
\]
The first summand can be transformed as
\[
\int_{\partial_h c_{s_0 s_t}} \ol\omega \eq \int_{c_{s_0 s_t}} \ol\eta.
\]
Over the path we are considering $\Omega_S^1$ is generated by $dt$. Therefore $\ol\eta=dt\wedge \ol\eta_0$. We obtain
\[
\int_{c_{s_0 s_t}} \ol\eta\eq\int_{c_{s_0 s_t}} dt \wedge \ol\eta_0\eq
\int_{c_{s_0 s_t}} (d(t\ol\eta_0) - t d\ol\eta_0)\eq\int_{\partial_h c_{s_0 s_t}} t\ol\eta_0 - \int_{c_{s_0 s_t}} t d\ol\eta_0.
\]
Examining the second term we see
\[
{\frac{\partial}{\partial t}}\bigg\vert_{t=0} \int_{c_{s_0 s_t}} t d\ol\eta_0 \eq0.
\]
The first one transforms to
\[
\int_{\partial_h c_{s_0 s_t}} t\ol\eta_0 \eq \int_{c_{s_t}-c_{s_0}} t\ol\eta_0 -\int_{(\partial_h c)_{s_0 s_t}} t\ol\eta_0.
\]
The second term gives
\[
{\frac{\partial}{\partial t}}\bigg\vert_{t=0} \int_{(\partial_h c)_{s_0 s_t}} t\ol\eta_0\eq0.
\]
The first one gives
\[
{\frac{\partial}{\partial t}}\bigg\vert_{t=0} \int_{c_{s_t}-c_{s_0}} t\ol\eta_0 \eq \lim_{t\ra 0} \int_{c_{s_t}} \ol\eta_0\eq\int_{c_{s_0}} \ol\eta_0.
\]
\end{proof}

For example, if $c_s$ is closed:
\begin{cor}
If $c_s$ is a closed family of hyperchains and $\omega$ is a closed relative hyperform, then
\[
d \int_{c_s} \omega \eq \int_{c_s} \nabla\omega.
\]
\end{cor}

\section{The Abel-Jacobi map for products of curves}
By a curve we mean a smooth projective curve over $\C$.
Let $X_1$, $X_2$,\dots,$X_n$ be curves. Put 
\[
X \eq X_1\times X_2\times \dots \times X_n.
\]
Let $x\in Z^k(X,1)$ be a higher cycle,
\[
x \eq \sum_i (W_i, f_i).
\]
Recall that $\dim_\C W_i = n-k+1$, $f_i\in\C(W_i)$, $\gamma_i = f_i^*[0,\infty]$, $\gamma = \sum_i \gamma_i$, $\partial\xi = \gamma$. The Abel-Jacobi map was defined as
\begin{equation}\label{defaj}
\langle AJ^{k,1} [x], [w] \rangle \eq \sum_i\int_{W_i\setminus\gamma_i}\omega\log f_i \,+\, 2\pi\I\, \int_\xi\omega
\end{equation}
for $w\in F^{n-k+1}\A^{2n-2k+2}(X)$, $dw=0$. 
\begin{rem}
Since $X$ is a product of curves its cohomology has no torsion. This implies that the class of $\gamma$ is trivial.
\end{rem}

\subsection{Triangulations}
In fact we should justify the construction of $\gamma_i$ and the integration in (\ref{defaj}). The problems are that $W_i$ are not necessarily smooth, the rational functions $f_i$ do not necessarily define maps to $\PP^1$ and $\log f_i$ is not bounded. To define all the objects we may embed $X$ as a semi-algebraic set into some $\R^N$ and  consider the semi-algebraic subsets $W_i$, $|\Div f_i|$, $f_i^{-1} [0,\infty]$. We may apply \cite{hironaka:tri} to get a triangulation of $\R^N$ which is compatible with all the sets mentioned above. This triangulation is semi-algebraic and smooth on the interiors of simplices. Therefore we can integrate smooth forms over simplices. Moreover we have the necessary bounds on the growth of $\log f_i$ restricted to any simplex which belongs to $W_i$ and does not belong to $|\Div f_i|$. In fact any such simplex intersects $|\Div f_i|$ only along the boundary and $\log f_i$ grows not faster than some multiple of the logarithm of the distance to the boundary. We will only consider simplices obtained by a linear subdivision of simplices of the constructed triangulation. Formal linear combinations of simplices of equal dimension are called chains.

For each $i$ let us consider the space $V_i$ constructed from $W_i$ by cutting out $\gamma_i$ and attaching two copies of $\gamma_i$ glued together along the boundary. Let us denote the two copies of $\gamma_i$ by $\gamma_i^+$ and $\gamma_i^-$ and suppose that they are attached in such a way that the function $\log f_i$ extends to $\gamma_i^+$ and $\gamma_i^-$, $\partial V_i = \gamma_i^+ - \gamma_i^-$, and the value of $\log f_i$ on $\gamma_i^-$ is $2\pi\I$ plus the value on $\gamma_i^+$. Denote by $\iota_i^+$, $\iota_i^-$ the natural isomorphisms $\gamma_i\To\gamma_i^+$, $\gamma_i\To\gamma_i^-$. Denote by $p_i$ the natural projection $V_i\To W_i$. The space $V_i$ is naturally endowed with the triangulation coming from the triangulation of $W_i$.

Let $l\geq 0$ and $\Delta$ be the standard simplex of dimension $l$,
\[
\Delta \eq \{(x_0, x_1, \dots, x_l)\;|\; \sum_{j=0}^l x_j = 0,\; x_j\geq 0\}.
\]
For $\epsilon>0$ denote
\[
\Delta_\epsilon \eq \{(x_0, x_1, \dots, x_l)\in\Delta\;|\; x_j\geq \epsilon\}.
\]
Let $\sigma:\Delta \To V_i$ be a simplex in $V_i$. Let $\omega$ be a smooth $l$-form on a neighborhood of $\sigma(\Delta)$. 
\begin{defn}
Suppose $\sigma(\Delta)$ is not contained in $|\Div f_i|$. Put
\[
\int_\sigma \omega \log f_i \eq \lim_{\epsilon\To 0} \int_{\Delta_\epsilon} \sigma^* \omega \log f_i.
\]
\end{defn}
\begin{defn}
A simplex $\sigma$ is called good if it is not contained in $|\Div f_i|$ for all $i$ and any simplex of its boundary is not contained in $|\Div f_i|$ for all $i$. A chain is called good if it is a linear combination of good simplices.
\end{defn}
One can check that the Stokes formula holds:
\begin{prop}
If $\sigma$ is a good simplex and $\omega$ is a smooth $l-1$-form on a neighborhood of $\sigma(\Delta)$, then
\[
\int_\sigma d\omega \log f_i \eq \int_{\partial \sigma} \omega \log f_i \,-\, \int_\sigma \frac{d f_i}{f_i}\wedge \omega.
\]
\end{prop}

\subsection{Using hypercovers}
Let us choose a Zariski affine cover on each of the curves $X_j$ this gives a hypercover $\U_j$ on $X_j$ indexed by $\sigma_j$, hence a product hypercover $\U$ on $X$ indexed by $\sigma$. This induces hypercovers on $V_i$. Let us lift chains $\Div f_i$ to hyperchains $\wt{\Div f_i}$, then lift the chains $\gamma_i$ to hyperchains $\wt\gamma_i$, then lift $V_i$ to $\wt V_i$ and $\xi$ to $\wt\xi$. Lifting $\gamma_i$ determines a lifting of $\gamma_i^+$, $\gamma_i^-$. These liftings are denoted $\wt\gamma_i^+$, $\wt\gamma_i^-$. We choose a nice analytic refinement $\U'$ of $\U$ and assume that all the hyperchains we are considering are in fact hyperchains on $\U'$. The following relations hold:
\[
\partial_h\wt\gamma_i \eq -\wt{\Div f_i},\quad \partial_h \wt\xi\eq\sum_i\wt\gamma_i,\quad \partial_h\wt V_i \eq \wt\gamma_i^+-\wt\gamma_i^-.
\]
For $\omega\in F_{n-k+1}\A^{2n-2k+2}(X)$, $d\omega=0$ we can denote by the same letter $\omega$ the corresponding hyperform. We get
\[
\langle AJ^{k,1}[x],[\omega]\rangle\eq \sum \int_{\wt V_i} \omega \log f_i \,+\, 2\pi\I\,\int_{\wt\xi}\omega.
\]
Let $L\subset\{1,\ldots,n\}$, $L'\subset\{1,\ldots,n\}$ such that $L\cap L'=\varnothing$, $|L|=n-k+1$, $|L'|=n-k+1$. Let $\omega_k$ be a holomorphic $1$-form on $X_k$ for $k\in L$ (which is considered as a hyperform) and a closed algebraic $1$-hyperform on $X_k$ for $k\in L'$. Recall that $1$-hyperform on $X_k$ is a collection of $1$-forms on the open sets of the cover of $X_k$ and functions on the pairwise intersections. Saying algebraic we require the functions and the $1$-forms to be regular. For $k\in L'$ by the Hodge theory one can choose a smooth $0$-hyperform $g_k$ such that $\omega_k'=\omega_k-d g_k$ is a smooth $1$-form. Put $\omega_k'=\omega_k$ for $k\in L$. We obtain a smooth closed form on $X$
\[
\omega'\eq\times_{k\in L\cup L'} \omega_k'\;\in\; F_{n-k+1}\A^{2n-2k+2}(X),
\]
and an algebraic closed hyperform
\[
\omega \eq \times_{k\in L\cup L'} \omega_k \;\in\; \Omega^{2n-2k+2}(\U)
\]
which is zero on $U_a$ if $\dim a>n-k+1$. Moreover they differ by a coboundary,
\[
\omega'-\omega \eq d g, \quad g\;\in\; \Omega_{\smooth}^{2n-2k+1}(\U),
\]
with $g$ zero on $U_a$ if $\dim a>n-k$.

\begin{prop}
One can compute $AJ^{k,1}$ using algebraic forms:
\[
\langle AJ^{k,1}[x],[\omega]\rangle\eq \sum_i \int_{\wt V_i} \omega \log f_i \,+\, 2\pi\I\,\int_{\wt\xi}\omega.
\]
\end{prop}
\begin{proof}
By the definition
\begin{multline*}
\langle AJ^{k,1}[x],[\omega']\rangle\eq \sum_i \int_{\wt V_i} \omega' \log f_i + 2\pi\I\,\int_{\wt\xi}\omega' \\\eq \sum_i \int_{\wt V_i} \omega \log f_i + 2\pi\I\,\int_{\wt\xi}\omega + \sum_i \int_{\wt V_i} (dg) \log f_i + 2\pi\I\,\int_{\wt\xi} dg.
\end{multline*}
The summand in the third term gives
\[
\int_{\wt V_i} (dg) \log f_i \eq \int_{\wt\gamma_i^+-\wt\gamma_i^-} g \log f_i - \int_{\wt V_i} \frac{d f_i}{f_i} \wedge g.
\]
The second integral is zero. Indeed, if $a\in \sigma$ is such that $\dim a>n-k$, then $(\frac{d f_i}{f_i} \wedge g)_a$ is zero as it was mentioned above. Otherwise $(\frac{d f_i}{f_i} \wedge g)_a\in\Omega^{2n-2k+2-\dim a}(U_a)$. Therefore the degree of this form is at least $n-k+2$, which is greater than the dimension of $V_i$. The first integral can be further transformed to
\[
\int_{\wt\gamma_i^+-\wt\gamma_i^-} g \log f_i\eq-2\pi\I\,\int_{\wt\gamma_i} g
\]
which cancels the corresponding summand in the term
\[
\int_{\wt\xi} dg\eq\sum_i\int_{\wt\gamma_i} g.
\]
\end{proof}

\subsection{Differentiating the Abel-Jacobi map}
Suppose we have families of curves $X_i\ra S$ with $S$ affine and 
\[
X \eq X_1\times_S X_2\times_S \dots \times_S X_n.
\]
Suppose we have a family of cycles $x_s\in Z^k(X_s,1)$. Let $\omega\in\Omega_{X/S}^{2n-2k+2}(\U)$. Suppose $\omega_a=0$ if $\dim a>n-k+1$. Let us compute 
\[
d \langle AJ^{k,1}[x_s],[\omega]\rangle.
\]
We first generalize it. Let for any $\omega\in \Omega_{X/S}^{2n-2k+2}(\U)$
\[
I(\omega)\eq \sum_i \int_{\wt V_i} \omega \log f_i \,+\, 2\pi\I\,\int_{\wt\xi}\omega.
\]
We also define a new operation on hyperforms. Let $\ol\omega\in\Omega_{X}^{2n-2k+2}(\U)$. For any index $i$ consider the hyperform $\frac{d f_i}{f_i}\wedge\ol\omega$. It is zero on all $U_a\cap W_i$ with $\dim a<n-k+1$ and is a form of maximal degree on $U_a\cap W_i$ with $\dim a=n-k+1$. Therefore there exists an element $\ol \phi_{f_i} \omega\in \Omega^1(S)\otimes \Omega_{X}^{2n-2k+2}(\U)$ such that the difference $\frac{d f_i}{f_i}\wedge\ol\omega-\ol \phi_{f_i} \omega$ is zero on all $U_a\cap W_i$ with $\dim a\leq n-k+1$. Its image in $\Omega^1(S)\otimes \Omega_{W_i/S}^{2n-2k+2}(\U)$ will be denoted by $\phi_{f_i} \omega$.
\begin{prop}\label{prop:diff_aj:1}
We have 
\[
d I(\omega) \eq I(\nabla\omega) \,+\, \int_{\wt W_i} \phi_{f_i} \omega.
\]
\end{prop}
\begin{proof}
Let $\eta=\nabla \omega$. Then (see Lemma \ref{lem:gauss-manin:2})
\[
d \int_{\wt\xi}\omega \eq \int_{\wt\xi}\eta \,+\, \sum_i R(\wt\gamma_i,\ol\omega).
\]
Consider the hyperform $\omega\log f_i$ on $\U\cap V_i$. The restriction of $\omega\log f_i$ to $\U'\cap V_i$ is closed simply because $\Omega_{W_i/S}^{2n-2k+3}(\U'\cap W_i)=0$. We have 
\[
d(\ol\omega \log f_i)\eq \ol\eta\log f_i \,+\, \frac{d f_i}{f_i} \wedge \ol\omega.
\]
The sum $\eta\log f_i + \phi_{f_i}\omega$ gives a Gauss-Manin derivative of $\omega\log f_i$ on $V_i$ with respect to the hypercover $\U'$. Therefore
\[
d \int_{\wt V_i} \omega \log f_i \eq R(\wt\gamma_i^+ - \wt\gamma_i^-,\ol\omega\log f_i) \,+\, \int_{\wt V_i} \eta\log f_i \,+\, \int_{\wt V_i} \phi_{f_i}\omega.
\]
The first summand equals to 
\[
R(\wt\gamma_i^+-\wt\gamma_i^-,\ol\omega\log f_i) \eq -2\pi\I R(\wt\gamma_i,\ol\omega).
\]
The third summand equals to
\[
\int_{\wt V_i} \phi_{f_i}\omega\eq\int_{\wt W_i} \phi_{f_i}\omega.
\]
\end{proof}

We also can compute $I(\omega)$ for an exact hyperform.
\begin{prop}\label{prop:diff_aj:2}
Let $\omega=d\eta$ with $\eta\in\Omega_{X/S}^{2n-2k+1}(\U)$. Then
\[
I(\omega)\eq- \sum_i\int_{\ol W_i}\frac{d f_i}{f_i}\wedge \eta
\]
\end{prop}
\begin{proof}
By the definition
\[
I(\omega)\eq\sum_i \int_{\wt V_i} \omega \log f_i \,+\, 2\pi\I\,\int_{\wt\xi}\omega.
\]
We have $\omega\log f_i=d(\eta \log f_i)-\frac{d f_i}{f_i}\wedge\eta$. Applying the Stokes formula we obtain the result.
\end{proof}

Let us denote 
\[
J_i(\omega)\eq(2\pi\I)^{-n+k-1} \int_{\wt V_i} \omega.
\]
for any $\omega\in \Omega^{2n-2k+2}_{W_i/S}(\U)$. We see that the following relations hold:
\[
d I(\omega) \eq I(\nabla\omega) + (2\pi\I)^{n-k+1} \sum_i J_i(\phi_{f_i}\omega),\; I(d\eta)\eq-(2\pi\I)^{n-k+1} \sum_i J_i(\frac{d f_i}{f_i}\wedge\eta).
\]
Note that the value of $J_i$ can be expressed as a certain sum of iterated residues. Therefore it is possible to compute $J_i$.

\subsection{Extensions of $\D$-modules} \label{subs:ext_dmod}
Recall that $S$ is an affine variety. Suppose $S=\spec R$ for a commutative ring $R$. Suppose $S$ is smooth and $\Omega(R)$ is free. Let $\D$ be the ring of differential operators on $S$. Consider $R$-modules $\Omega_X^i(\U)$. We have two filtrations on $\Omega_X^i(\U)$. The first is the Hodge filtration. Elements of $F^j\Omega_X^i(\U)$ are those hyperforms which have as components only forms of rank at least $j$. We have another filtration defined as $G^j\Omega_X^i(\U):=\Omega^j(R)\wedge\Omega_X^{i-j}(\U)$. The exterior derivative respects these filtrations. We have 
\[
\Omega_{X/S}^i(\U)\eq\Omega_X^i(\U)/G^1\Omega_X^i(\U).
\]

To understand the results of the previous section we introduce two operations on hyperforms. The first one is $\Psi_0:\Omega_X^{2n-2k+1}\ra R$ defined as
\[
\Psi_0(\eta)\eq(2\pi\I)^{-n+k-1}\sum_i\int_{\wt W_i}\frac{d f_i}{f_i}\wedge \eta \quad (\eta\in\Omega_X^{2n-2k+1}(\U)).
\]
The second one is $\Psi_1:\Omega_X^{2n-2k+2}\ra \Omega(R)$ defined as
\[
\Psi_1(\omega)\eq(2\pi\I)^{-n+k-1}\sum_i\int_{\wt W_i}\phi_{f_i}\omega \quad 
(\omega\in\Omega_X^{2n-2k+2}(\U)),
\]
the operation $\phi$ being defined in the previous section.

\begin{rem}
The integral with respect to $\wt W_i$ can be expressed as a sum of iterated residues as is proved in the previous chapter. Therefore the operations $\Psi_0$ and $\Psi_1$ can be defined purely algebraically
\end{rem}

\begin{prop}\label{prop:ext-dmod:1}
The first operation satisfies the following properties:
\begin{enumerate}
\item $\Psi_0$ is $R$-linear.
\item $\Psi_0|_{G^1\Omega_X^{2n-2k+1}(\U)}=0$.
\item $\Psi_0|_{F^{n-k+1}\Omega_X^{2n-2k+1}(\U)}=0$.
\item $\Psi_0(\eta)=-(2\pi\I)^{-n+k-1} I(d\eta)$ for $\eta\in\Omega_X^{2n-2k+1}(\U)$.
\item If $d\eta\in G^1\Omega_X^{2n-2k+2}(\U)$, then $\Psi_0(\eta)=0$.
\end{enumerate}
\end{prop}

\begin{prop}\label{prop:ext-dmod:2}
The second operation satisfies the following properties:
\begin{enumerate}
\item $\Psi_1$ is $R$-linear.
\item $\Psi_1|_{G^2\Omega_X^{2n-2k+2}(\U)}=0$.
\item $\Psi_1|_{F^{n-k+2}\Omega_X^{2n-2k+2}(\U)}=0$.
\item $\Psi_1(u\wedge\eta)=-u\Psi_0(\eta)$ for $u\in\Omega(R),\eta\in\Omega_X^{2n-2k+1}(\U)$.
\item $\Psi_1(\omega)=(2\pi\I)^{-n+k-1} (d I(\omega)-I(\nabla\omega))$ for $\omega\in\Omega_X^{2n-2k+2}(\U)$ such that $d\omega\in G^1\Omega_X^{2n-2k+3}(\U)$, and $\nabla\omega$ is defined as the class of $d\omega$ in $\Omega(R)\otimes\Omega_X^{2n-2k+2}(\U)$.
\item $\Psi_1(d \eta)=d \Psi_0(\eta)$ for $\eta\in\Omega_X^{2n-2k+1}(\U)$.
\end{enumerate}
\end{prop}

Let $B^{2n-2k+2}\subset Z^{2n-2k+2}\subset \Omega_{X/S}^{2n-2k+2}(\U)$ be defined as 
\[
B^{2n-2k+2}\eq\image d_{X/S},\; Z^{2n-2k+2}\eq\kernel d_{X/S}.
\]
Let 
\[
H^{2n-2k+2} \eq Z^{2n-2k+2}/B^{2n-2k+2}.
\]
Note that the Hodge theory implies 
\[
F^i \Omega_{X/S}^{2n-2k+2} \cap B^{2n-2k+2}\eq d_{X/S}(F^i\Omega_{X/S}^{2n-2k+1}).
\]
Let $M$ be defined as 
\[
M\eq(R\oplus Z^{2n-2k+2})/\image(\Psi_0,d_{X/S}).
\]
We then have the following commutative diagram with exact rows:
\[
\begin{CD}
0 @>>>B^{2n-2k+2} @>>> Z^{2n-2k+2} @>>> H^{2n-2k+2} @>>> 0\\
& & @V{\Psi_0}VV @V{i_2}VV @| \\
0 @>>>R @>{i_1}>> M @>>> H^{2n-2k+2} @>>> 0
\end{CD}
\]
Note that there is a canonical section $s:F^{n-k+1} H^{2n-2k+2}\ra M$ which is defined by sending a class $[\omega]\in F^{n-k+1} H^{2n-2k+2}$ ($\omega\in Z^{2n-2k+2}$) to $(0,\omega)$.

Both $R$ and $H^{2n-2k+2}$ are modules over $\D$. On $M$ we have the following $\D$-module structure:
\[
v(r, \omega)\eq(v(r)-\langle v, \Psi_1(\ol\omega)\rangle , \langle v, \nabla \ol\omega\rangle),
\]
where $r\in R$, $v\in\Der(R)$, $\omega\in Z^{2n-2k+2}$, $\ol\omega\in\Omega_X^{2n-2k+2}(\U)$ represents $\omega$.
One can check that this indeed gives a $\D$-module structure. The section $s$ is compatible with the $\D$-module structure so that the following diagram commutes:
\[
\begin{CD}
F^{n-k+2} H^{2n-2k+2} @>{\nabla}>> \Omega(S)\otimes F^{n-k+1} H^{2n-2k+2} \\
@V{s}VV @V{\id_S\otimes s}VV\\
M@>>> \Omega(S)\otimes M.
\end{CD}
\]
Let us choose $S$ smaller so that $H^{2n-2k+2}$ is free and $F^{n-k+1} H^{2n-2k+2}$ is a direct summand. Extend the section $s$ to a section $\wt s:H^{2n-2k+2}\ra M$. This provides an isomorphism of $R$-modules $M\cong R\oplus H^{2n-2k+2}$. We define a homomorphism $\Psi\in \Hom_R(H^{2n-2k+2},\Omega(R))$ in the following way. For any $[\omega]\in H^{2n-2k+2}$ put
\[
\Psi([\omega])\eq\nabla \wt s([\omega]) - (\id\otimes \wt s)(\nabla[\omega]) \;\in\; \Omega(R).
\]
One can see that $\Psi|_{F^{n-k+2}H^{2n-2k+2}}=0$ and $\Psi$ is correctly defined up to the differential of an element of $\Hom_R(H^{2n-2k+2},\Omega(R))$. The structure of $\D$-module on $R\oplus H^{2n-2k+2}$ induced from $M$ can be recovered as follows:
\[
v(r,h)\eq(v(r)+\Psi(h),\nabla_v h) \qquad (h\in H^{2n-2k+2}, r\in R, v\in\Der(R)).
\]

Let $N$ be the kernel of the multiplication homomorphism
\[
N \;:=\; \kernel(\D\otimes_R H^{2n-2k+2} \ra H^{2n-2k+2}).
\]
It inherits the filtration $F^\bullet$ from the one on $H^{2n-2k+2}$. We have the canonical homomorphism of $\D$-modules $\Psi'_{\alg}:F^{n-k+1} N\ra R$ defined as follows.  
\[
\Psi'_{\alg} (\sum_j \alpha_j \otimes h_j)=\sum_j \alpha_j s(h_j)\in R\quad (\sum_j \alpha_j \otimes h_j\in N,h_j\in F^{n-k+1} N,\alpha_j\in\D).
\]

On the other hand we could have defined $\Psi'_{\analyt}:F^{n-k+1} N\ra R$
\[
\Psi'_{\analyt} (\sum_j \alpha_j \otimes h_j)\eq\sum_j \alpha_j \langle AJ^{k,1}[x], h_j\rangle,
\]
but the properties of $\Psi_0$, $\Psi_1$ imply:
\begin{cor}
We have $\Psi'_{\analyt}=(2\pi\I)^{n-k+1}\Psi'_{\alg}$. Therefore $\Psi'_{\analyt}$ can be computed algebraically.
\end{cor}

\subsection{Products of elliptic curves}\label{prod-ellc}
Suppose we have a family of elliptic curves $\pi:E\ra S$ and $X_i=E$. We have $S=\spec R$. We suppose that $R$ is a $1$-dimensional domain and its field of fractions is denoted as $R_0$. Let 
\[
V \eq R^1\pi_*\Omega_{E/S}^\bullet.
\]
This is a locally free module over $R$ of rank $2$. The Hodge filtration has two pieces, each a rank $1$ locally free $R$-module. Denote
\[
M^1 \eq F^1 V.
\]
This is a line bundle of modular forms of weight $1$. We denote $M^j=(M^1)^{\otimes j}$ for $j\in\Z$. We have the canonical pairing
\[
F^1 V \otimes_R (V/F^1) \ra R.
\]
With the help of this pairing we identify
\[
V/F^1 \cong M^{-1}.
\]

Consider the Kodaira-Spencer map 
\[
KS: M^1\ra \Omega(R)\otimes M^{-1}.
\]
Suppose that the family $E$ is not constant and $S$ is small enough, so that $KS$ is an isomorphism. Therefore we have $\Omega(R)\cong M^2$. Let $H^1=R^1\pi_*\Omega_{E/S}^\bullet$, which is a $\D$-module of rank $2$ over $R$.

Suppose $n=2k-2$. Then, using the notations of the previous section ($H^n=R^n\pi_*\Omega_{E^n/S}^\bullet$)
\[
M^n \eq F^n H^{2n-2k+2} \eq F^n H^n.
\]
Let $H^n_\sym$ denote the direct summand of $H^n$ which corresponds to 
\[
H^n_\sym\eq\Sym^n H^1.
\]

We have the following fact:
\begin{prop}
There is a unique element $B_n\in M^{n+2}\otimes_R F^n N$ which contains only differential operators of degree at most $n+1$ and the corresponding symbol in $M^{n+2}\otimes_R \Der(R)^{\otimes (n+1)}\otimes_R M^n\cong R$ is $1$. 
\end{prop}
\begin{proof}
Let us prove uniqueness. Let $\D^n$ denote the differential operators of order at most $n$. The following map is a monomorphism: 
\[
m:\D^n\otimes_R F^n H^n\ra H^n.
\]
To prove this look at the filtration by the order of differential operator on $\D^n\otimes_R F^n H^n$ and the Hodge filtration on $H^n$. The graded pieces of the map are the Kodaira-Spencer maps which are injective.

Existence follows from the fact that the image of $m$ is $\Sym^n H^1$. So we can pick any element of $M^{n+2}\otimes_R \D^{n+1}\otimes_R M^n$ with symbol $1$ and then subtract an element of $M^{n+2}\otimes_R \D^n\otimes_R M^n$ which maps to the same element in $M^{n+2} \otimes_R \Sym^n H^1$.
\end{proof}

One can apply $\Psi'$ to $B=B_n$ and get an element 
\[
(2\pi\I)^{-n+k-1}\Psi'_{\analyt}(B)\eq\Psi'_{\alg}(B) \;\in\; M^{n+2}.
\]
Therefore we have a canonical modular form constructed from a family of cycles $\{x_s\in Z^k(X_s,1)\}_{s\in S}$.

\subsection{Analytic computations}\label{subs:anal_comp}
We still have $n=2k-2$. Suppose we have a family of elliptic curves $\pi:E\ra S$ over an affine smooth irreducible curve. We now translate the notions from Section \ref{prod-ellc} to the analytic language. The elliptic curve over a point $s\in S$ is denoted $E_s$. Let $U$ be an analytic subset in $S$ homeomorphic to a disk. Choose families of $1$-cycles $c_1$, $c_2$ over $U$ such that $\cc_1(s)$ and $\cc_2(s)$ generate $H_1(E_s, \Z)$ and the intersection number is $\cc_1\cdot \cc_2=1$.

Any other choice $\cc_1'$, $\cc_2'$ can be obtained from the choice $\cc_1$, $\cc_2$ by a transformation 
\[
\gamma\eq\begin{pmatrix}a & b\\ c & d\end{pmatrix}\in SL_2(\Z):\qquad \cc_2' \eq a \cc_2 + b \cc_1,\; \cc_1' \eq c \cc_2 + d \cc_1.
\]

Let $\omega$ be a closed relative differential $1$-form on $E$. We denote
\[
\Omega_1(\omega)\eq\int_{\cc_1} \omega, \qquad \Omega_2(\omega)\eq\int_{\cc_2} \omega.
\]
The cup product provides a pairing:
\[
(\omega_1,\omega_2) \eq \int_{E_s} \omega_1\wedge \omega_2.
\]
Let us integrate $\omega_1=df$ over the universal cover of $E_s$. Then 
\[
(\omega_1,\omega_2) \eq \int_{\partial \wt E_s} f \omega_2 \eq \Omega_1(\omega_1)\Omega_2(\omega_2)-\Omega_1(\omega_2) \Omega_2(\omega_1),
\]
where $\wt E_s$ denotes a fundamental domain of the universal cover of $E_s$.

If $\omega$ is holomorphic we put
\[
z\eq\frac{\Omega_2(\omega)}{\Omega_1(\omega)}.
\]
This locally defines an isomorphism between $S$ and the upper half plane. Indeed,
\[
\Im z \eq \frac{\Omega_2(\omega)\ol\Omega_1(\omega)-\Omega_1(\omega) \ol\Omega_2(\omega)}{2\I \Omega_1(\omega)\ol\Omega_1(\omega)} \eq -\frac{\Im \int_{E_s} \omega \wedge\ol\omega}{2 \Omega_1(\omega)\ol\Omega_1(\omega)}.
\]
If we represent $E_s$ as a quotient $\C/\Lambda$ with $\omega=dx+\I dy$, then
\[
\omega\wedge\ol\omega \eq -2 \I dx\wedge dy,
\]
therefore
\[
\Im z \eq \frac{\vol_\omega E_s}{|\Omega_1(\omega)|^2} \;>\; 0,
\]
where $\vol_\omega E_s$ is the volume of $E_s$ defined with the help of the form $\omega$. For another choice $\cc_1'$, $\cc_2'$ we obtain
\[
z' \eq \frac{a z + b}{c z + d}.
\]
We define a canonical isomorphism of the analytic version of the sheaf $M^1$ of modular forms of weight $1$, as defined in the previous section, and the pullback via $z$ of the usual sheaf of modular forms of weight $1$ on the upper half plane:
\[
\omega \ra \Omega_1(\omega).
\]

Let $X$ be a formal variable. We identify $H^1(E_s, \C)$ with the space of polynomials of degree not greater than $1$ in $X$ in the following way:
\[
\langle 1, \cc_2\rangle\eq-1,\; \langle 1, \cc_1\rangle\eq0,\; \langle X, \cc_1\rangle\eq1,\; \langle X, \cc_2\rangle\eq0.
\]

Let $\omega$ be a closed differential $1$-form. Then the corresponding polynomial is 
\[
[\omega]_\cc\eq\Omega_1(\omega) X-\Omega_2(\omega).
\]
In particular, if $\omega$ is holomorphic,
\[
[\omega]_\cc\eq\Omega_1(\omega) (X-z).
\]
If $\cc'=\gamma \cc$, then one can check that
\[
[\omega]_{\cc'} \eq \gamma(\Omega_1(\omega) X - \Omega_2(\omega)),
\]
where the action on polynomials is defined as
\[
\gamma(p) \eq p|_{-1} \gamma^{-1}\qquad (p=p_1 X + p_0).
\]
Therefore the map
\[
\omega \ra f_\omega\eq\frac{\Omega_1(\omega) X-\Omega_2(\omega)}{X-z} \eq \Omega_1(\omega) + \frac{\Omega_2(\omega)-z \Omega_1(\omega)}{z-X}
\]
defines an isomorphism between the sheaf $V=H^1(E,\C)$ and the pullback of the sheaf of quasi-modular forms of weight $1$ and depth $1$. \emph{Quasi-modular forms} of weight $w$ and depth $d$ are functions of the form
\[
f(z, X) \eq \sum_{i=0}^d \frac{f_i(z)}{i!\, (z-X)^i}
\]
which transform like modular forms of weight $w$ in $z$ and weight $0$ in X:
\[
f(\gamma(z), \gamma(X))\, (c z+d)^{-w} \eq f(z, X).
\]

One can see that the pairing can be written as
\[
(aX+b, a' X + b') \eq - a b' +a' b.
\]
Therefore if $a(X-z)\in F^1 V$ and $a'X+b'\in V$, then
\[
(a(X-z), a'X+b') \eq -a (a' z + b').
\]
If $\omega$ is a holomorphic differential and $\eta$ an arbitrary closed differential, then
\[
(\omega, \eta) \eq f_\omega (f_{\eta})_1,
\]
where $(f_{\eta})_1$ denotes the coefficient at $(z-X)^{-1}$ of $f_{\eta}$ and is a modular form of weight $-1$. Therefore the isomorphism $V/F^1 V\ra M^{-1}$ is given by sending $f$ to $f_1$.

\emph{The Gauss-Manin derivative} with respect to the parameter $z$ sends the cohomology class with periods $\Omega_1$, $\Omega_2$ to the cohomology class with periods $\frac{\partial \Omega_1}{\partial z}$, $\frac{\partial \Omega_2}{\partial z}$. Therefore on the level of quasi-modular forms it can be written as 
\[
\frac{\partial}{\partial z} + \frac{1}{z-X}.
\]
If we take a modular form $f$ of weight $1$, the Gauss-Manin derivative of the corresponding family of cohomology classes will be given by
\[
\left(\frac{\partial f}{\partial z} + \frac{f}{z-X}\right) dz.
\]
Therefore the \emph{Kodaira-Spencer} map sends 
\[
f\ra f dz.
\]
If $t$ is a local parameter on $S$ then the Kodaira-Spencer map acts as
\[
f \ra \frac{f}{\left(\frac{dt}{dz}\right)} dt,
\]
therefore the isomorphism $\Omega(S)\ra M^2$ acts by sending $dt$ for a function $t$ to the modular form $\frac{dt}{dz}$.

Next we construct the canonical differential operator from the previous section. Let $U$ be an open subset in $S$ such that there exist modular forms $f$ of weight $n$ and $g$ of weight $n+2$ with non-zero values on $U$. Let $D(f,g)$ be the operator
\[
\phi\ra \frac{1}{g} \left(\frac{\partial}{\partial z}\right)^{n+1} \frac{\phi}{f}.
\]
This is a differential operator which sends functions to functions, therefore $D(f,g)\in \D^{n+1}(U)$. Moreover $D(f,g) f = 0$, therefore 
\[
B_n(U) \eq g \otimes D(f,g)\otimes f
\]
defines a section of $M^{n+2}\otimes_R F^n N$. Its symbol is
\[
g \frac{1}{f g dz^{n+1}} f,
\]
which goes to $1$ under the isomorphism $M^{n+2}\otimes_R \Der(R)^{\otimes (n+1)}\otimes_R M^n\cong R$.

Let us consider the natural map
\[
\frac{H^{2k-2}}{F^k H^{2k-2} + 2\pi\I\,H^{2k-2}_\Z} \ra \frac{H^{2k-2}}{F^1 H^{2k-2} + 2\pi\I\,H^{2k-2}_\Z} \cong M^{2-2k}/2\pi\I\,M^{2-2k}_\Z,
\]
where $H^{2k-2}_\Z\subset H^{2k-2}$ is the subsheaf of integral cohomology classes and $M^{-j}_\Z$ for any $j\geq 0$ is the subsheaf of $M^{-j}$ generated by $1,z,\ldots,z^j$. The image of the section $AJ^{k,1}[x]$ under this map will be denoted by $A_x$. More explicitly, one can obtain $A_x$ by choosing a modular form $f$ of weight $2k-2$, and then dividing by $f$ the pairing of $A_x$ with the holomorphic differential $2k-2$~-form with periods given by $f$. It is clear that the operator 
\[
\left(\frac{\partial}{\partial z}\right)^{2k-1}: M^{2-2k}\ra M^{2k}
\]
vanishes on $M^{2-2k}_\Z$. Therefore it is defined on $M^{2-2k}/M^{2-2k}_\Z$ and it is clear that 
\[
\Psi_{\analyt}'(B_{2k-2}) \eq \left(\frac{\partial}{\partial z}\right)^{2k-1} A_x \;\in\; M^{2k}.
\]

Also we have the canonical (non-holomorphic) section of $F^{k-1}H^{2k-2}$. This is defined by the polynomial-valued function
\[
Q_z^{k-1}\eq\left(\frac{(X-z)(X-\ol z)}{z-\ol z}\right)^{k-1} \;\in\; \Sym^{k-1} V \;=\; H^{2k-2}.
\]

\begin{thm} \label{thm:anal-comp:9}
Suppose $S$ is a smooth affine curve.
Suppose we have a family of elliptic curves $\{E_s\}_{s\in S}$ and an algebraic family of higher cycles $\{x_s\}_{s\in S}$, $x_s\in Z^k(E_s^{2k-2}, 1)$. Suppose there is a map $\varphi$ from $S$ to $\HH/\Gamma$ which lifts the canonical map from $S$ to $\HH/SL_2(\Z)$. Suppose the symmetric part of $AJ^{k,1}[x]$ is constant along the fibers of $\phi$. Take the corresponding modular form $A_x$ of weight $2-2k$ on the image of $\phi$ defined locally up to polynomials in $z$ of degree not greater than $2k-2$ as above. Suppose 
\[
\left(\frac{\partial}{\partial z}\right)^{2k-1} A_x = (-1)^{k-1} D_0^{\frac{k-1}2} \alpha g_{k,z_0}^{\HH/\Gamma}(z)
\]
for a CM point $z_0$ of discriminant $D_0$ and a non-zero rational number $\alpha$. Then for any point $z$ in the image of $\phi$, $\phi(s)=z$,
\[
G_k^{\HH/\Gamma}(z, z_0) = 2 \alpha^{-1} \frac{(2k-2)!}{(k-1)!}\, \Re\left(D_0^{\frac{1-k}2} ( AJ^{k,1}[x], Q_z^{k-1} )\right).
\]
Take $N_A$, $N_B$, $N$ as in Corollary \ref{thm:global_st:6}. Suppose $z$ is a CM point of discriminant $D$ and make $N_A$ larger if necessary to satisfy $N_A (k-1)!\, \alpha^{-1}\in \Z$. Then the algebraicity conjecture is true and one has
\[
(D D_0)^{\frac{k-1}2} \,\widehat{G}^{k, \HH/\Gamma}(z, z_0) \;\equiv\; \alpha^{-1} \log (x_s \cdot Z_{z}) \mod \frac{2\pi\I}{N} \Z,
\]
where $Z_z$ is a subvariety of $E_s^{2k-2}$ which intersects $x_s$ properly and has cohomology class
\[
\cl Z_z = \frac{(2k-2)!}{(k-1)!}\, D^{\frac{k-1}2} Q_{z}^{k-1},
\]
and $\cdot$ denotes the intersection number as in Theorem \ref{thm:spec_val1}, which is an algebraic number if $x_s$ is defined over $\ol{\Q}$.
\end{thm}
\begin{proof}
Note that the function $A_x$ is holomorphic of weight $2-2k$. Hence it is (locally) of type $F_{k,2-2k}$ (see Section \ref{eigenvalues}). Therefore the following function is of type $F_{k,0}$:
\[
f \; := (-1)^{k-1}\, \frac{(k-1)!}{(2k-2)!} \; \delta^{k-1} A_x.
\]
We obtain the function $\wt f$ with values in $V_{2k-2}$. The function $\wt f$ satisfies (by Theorem \ref{thm:eigenfunc:5})
\[
(\wt f, (X-z)^{2k-2}) = (-1)^{k-1}\, \frac{(2k-2)!}{(k-1)!} (\wt f, \delta^{1-k} Q_z^{k-1}) = A_x.
\]
There is another formula for $\wt f$ (see Proposition \ref{prop:eigenfunc:6}):
\[
\wt f = \sum_{i=0}^{2k-2} \frac{(X-z)^i}{i!}\, \left(\frac{\partial}{\partial z}\right)^i A_x.
\]
This formula shows that when we add to $A_x$ a polynomial $p(z)$, the function $\wt f$ changes by $p(X)$. Therefore $\wt f$ is defined up to elements of $2\pi\I\,V_{2k-2}^\Z$. On the other hand, $AJ^{k,1}[x]$ is a function with values in 
\[
V_{2k-2}/(F^k V_{2k-2}+2\pi\I\,V_{2k-2}^\Z)
\]
($F^k$ corresponds to the polynomials divisible by $(X-z)^k$). So the difference satisfies
\[
AJ^{k,1}[x] - \wt f \;\in\; F^1V_{2k-2}/(F^k V_{2k-2}+2\pi\I\,V_{2k-2}^\Z).
\]

Consider the local sections of $F^{k-1} N$ given by 
\[
\xi_i = \frac{\partial}{\partial z} \otimes (X-z)^{i} \,+\, i\, (X-z)^{i-1}\qquad (i\geq k).
\]
It is easy to check that $\Psi'_{\alg} \xi_i=0$ using the property (iii) of the function $\Psi_1$. Therefore $\Psi'_{\analyt} \xi_i=0$. This means
\[
\frac{\partial}{\partial z} (AJ^{k,1}[x], (X-z)^{i}) \,+\, i\,(AJ^{k,1}[x], (X-z)^{i-1}) \eq 0\qquad (i\geq k).
\]
The function $\wt f$ satisfies similar property. Therefore their difference also does. Since
\[
(AJ^{k,1}[x] - \wt f, (X-z)^{2k-2}) \eq 0,
\]
we prove by induction that
\[
(AJ^{k,1}[x] - \wt f, (X-z)^i) \eq 0 \qquad (i\geq k-1).
\]

Next we note that (see Theorem \ref{thm:eigenfunc:5})
\[
d \wt f \eq \frac{(X-z)^{2k-2}}{(2k-2)!}\, \left(\left(\frac{\partial}{\partial z}\right)^{2k-1} A_x\right) dz.
\]
So if the hypothesis is true,
\[
d \wt f \eq (-1)^{k-1} \frac{\alpha D_0^{\frac{k-1}2}}{(2k-2)!}\, (X-z)^{2k-2}\, g_{k, z_0}^{\HH/\Gamma}(z)\, dz.
\]
Therefore one can choose (see Theorem \ref{int_pairing} and Corollary \ref{thm:global_st:6})
\[
I^{A,\Gamma}_{N_A (X-z)^{2k-2} g_{k, z_0}^{\HH/\Gamma}(z) dz} \eq (-1)^{k-1} N_A \alpha^{-1} (2k-2)!\, D_0^{\frac{1-k}2}\, \wt f,
\]
recall that $A = 2\pi\I \frac{(2k-2)!}{(k-1)!}\, D_0^{\frac{1-k}2} V_{2k-2}^\Z$. This shows that
\[
\widehat{G}_k^{\HH/\Gamma} (z, z_0) \;\equiv\; \alpha^{-1} \frac{(2k-2)!}{(k-1)!}\, D_0^{\frac{1-k}2} (\wt f, Q_{z}^{k-1}) \mod 2\pi \I (D D_0)^{\frac{1-k}2} \frac{1}{N} \Z,
\]
so the statement follows from Theorem \ref{thm:spec_val1}.
\end{proof}

\begin{rem}
A cycle with cohomology class $\frac{(2k-2)!}{(k-1)!} D^{\frac{k-1}2} Q_z^{k-1}$ was constructed in the introduction by taking the graphs of complex multiplication by $a z$ and $a \bar z$ ($a$ is the leading coefficient of the minimal quadratic equation of $z$), denoted $Y_{a z}$ and $Y_{a \bar z}$ respectively, and adding the products of $k-1$ copies of $Y_{a z} - Y_{a \bar z}$ for all possible splittings of the $2k-2$ elliptic curves of $E_z^{2k-2}$ into pairs. 
\end{rem}

\begin{rem}
As we will show later (Section \ref{sec:torsion}) for the case $\Gamma=PSL_2(\Z)$, $k=2$ Corollary \ref{thm:global_st:6} gives $N_B=2$, $N_A=1$. Since $N=(k-1)! N_A N_B$ this gives $N=2$. But if the numerator of $\alpha$ is greater then $1$, we take $N_A$ to be the numerator of $\alpha$ to satisfy the conditions of the theorem, so the statement will hold for $N=2 N_A$.
\end{rem}

\chapter{Cohomology of elliptic curves}
This chapter studies the Weierstrass family of elliptic curves 
\[
y^2\eq x^3+ax+b.
\]
We first study expansions at infinity of various functions. As a coordinate we use the formal integral of the holomorphic differential form $\frac{dx}{2y}$. Also we note that the base ring $\C[a,b]$ is isomorphic to the ring of modular forms for $SL_2(\Z)$ and we choose an isomorphism $\mu$.

In Section \ref{periods_dif} we state the precise relation between periods of differential forms of second kind and values of quasi-modular forms. 

In Section \ref{deriv-mf} we choose lifts of vector fields on the base to vector fields on the total space of the family. It happens that particularly nice formulae can be obtained for lifts of the Euler vector field and the Serre vector field. Therefore it is natural to choose these vector fields as a basis. We represent cohomology of elliptic curves by two differential forms of second kind $\frac{dx}{2y}$ and $\frac{xdx}{2y}$.

In Section \ref{repr-coh} we choose representatives of two cohomology classes, corresponding to the forms of second kind $\frac{dx}{2y}$ and $\frac{xdx}{2y}$, as hyperforms on the total space. This choice satisfies an important property. Whenever we apply the Gauss-Manin derivative (see Section \ref{gauss-manin}) to these representatives, the result is again a linear combination of these representatives. We express the hyperforms at infinity. Indeed, for computation of residues later it will be enough to know only these expressions, the result does not depend on the global information.

\section{Certain power series}

Let $R=k[a,b]$ be the ring of polynomials in two variables $a$, $b$. Denote by $K$ the field of fractions of $R$. Let $G_m$ be the multiplicative group. Let $G_m$ act on $R$ by the law
\[
a \ra \lambda^4 a,\; b\ra \lambda^6 b\qquad (\lambda\in G_m).
\]

We consider the family over $R$ given by the equation
\[
y^2\eq x^3+ax+b.
\]
This can be 'compactified' to the projective variety $E$ over $R$ given by the homogeneous equation in $\wt x$, $\wt y$, $\wt z$:
\[
\wt y^2 \wt z \eq \wt x^3+a \wt x \wt z^2 + b \wt z^3.
\]

The action of $G_m$ extends to the action on $E$ in the following way:
\[
\wt x\ra \lambda^2 \wt x,\; \wt y\ra \lambda^3 \wt y,\; \wt z\ra \wt z\qquad (\lambda\in G_m).
\]
Therefore the affine chart $\wt z=1$ is stable under the action. We denote this chart by $U_0$. In fact $E$ is an elliptic curve outside the zero locus of the discriminant
\[
\Delta\eq-16(4 a^3+27 b^2).
\]

If a rational function $\phi$ on $E$ transforms according to 
\[
\phi\ra \lambda^k \phi\qquad (\lambda\in G_m),
\]
then we say that $\phi$ is of weight $k$. Let us denote the space of rational functions of weight $k$ by $F_k$. The action of $G_m$ gives rise to the vector field whose derivation is the Euler operator, $\delta_e$. This operator acts on homogeneous rational functions as follows:
\[
\delta_e f \eq k f \qquad(f\in F_k).
\]

We have the zero section $s_0:\spec R \ra E$ given by sending 
\[
\wt x \ra 0, \wt y\ra 1, \wt z \ra 0.
\]
Let $t=-x/y = -\wt x/ \wt y \in F_{-1}$. This is a local parameter at $s_0$. We can express $x$ and $y$ as Laurent series in $t$:
\begin{align*}
x &\eq t^{-2}-a t^2-b t^4-a^2t^6 - 3 a b t^8+O(t^{10}),\\
y\eq-t^{-1}x &\eq -t^{-3}+a t + b t^3 + a^2 t^5 + 3 a b t^7+O(t^9).
\end{align*}

The invariant differential form $\omega=\frac{dx}{2y}$ has expansion
\[
\omega\eq\frac{dx}{2y}\eq(1+2a t^4+3b t^6+6 a^2 t^8+20 a b t^{10}+O(t^{12})) dt.
\]
Consider the formal integral of $\omega$:
\[
z \eq \int\omega\eq t+ \frac{2a}5 t^5+\frac{3b}7 t^7+\frac{2 a^2}3 t^9+\frac{20 a b}{11} t^{11}+O(t^{13}).
\]
In fact $z$ is the logarithm for the formal group law of the elliptic curve. We can now take $z$ as a new local parameter and express $x$ and $y$ in terms of $z$:
\begin{align*}
x &\eq z^{-2}-\frac{a}5 z^2-\frac{b}7 z^4+\frac{a^2}{75}z^6 + \frac{3 ab}{385}z^8+O(z^{10}),\\
y\eq\frac{\partial}{2\partial z}x &\eq -z^{-3}-\frac{a}5 z -\frac{2b}7 z^3 + \frac{a^2}{25} z^5 + \frac{12 a b}{385} z^7+O(z^9).
\end{align*}

Let us fix an isomorphism between $R$ and the ring of modular forms in the following way:
\[
\mu(a) \eq -\frac{E_4}{2^4 3},\; \mu(b) \eq  \frac{E_6}{2^5 3^3}.
\]
Then the integral of $-x dz$ can be expressed as follows:
\begin{align*}
v_0&\eq-\int x dz \eq z^{-1} +\frac{a}{15}z^3 +\frac{b}{35} z^5 - \frac{a^2}{525}z^7 - \frac{ab}{1155} z^9\\
&\eq z^{-1}-\frac{E_4}{720} z^3+\frac{E_6}{30240}z^5-\frac{E_4^2}{1209600} z^7+\frac{E_4 E_6}{47900160} z^9+O(z^{11})\\
&\eq z^{-1}+\sum_{k\geq 2} \frac{B_{2k}E_{2k}}{(2k)!} z^{2k-1}.
\end{align*}
In fact, this follows from the corresponding identity over the complex numbers which can be proved using the Taylor expansion of the Weierstrass $\wp$-function. We define
\[
v\eq v_0+\frac{E_2}{12} z \eq z^{-1}+\sum_{k\geq 1} \frac{B_{2k}E_{2k}}{(2k)!} z^{2k-1}\in R[E_2]((z)).
\]

Note that for $a=-\frac{1}{2^4 3}$, $b=\frac{1}{2^5 3^3}$ (this corresponds to $E_4=1$, $E_6=1$ and the curve is degenerate) we can find expansions of $v_0$, $x$ and $y$ explicitly:
\[
v_0 \eq \frac{1}{e^z-1} + \frac{1}{2} - \frac{z}{12},\qquad
x \eq \frac{e^z}{(e^z-1)^2}+\frac{1}{12},\qquad
y \eq -\frac{e^{2z}+e^z}{2(e^z-1)^3}.
\]
This corresponds to the fact that the Fourier expansion of $E_{2k}$ starts with $1$.

\section{Periods of differentials of second kind}\label{periods_dif}
Let us view $E$ as an elliptic curve over $K$. We will consider odd differential forms on $U_0$. Each such form has an expansion of the type
\[
\sum_{k\in \Z} a_{2k} z^{2k} dz\qquad (a_{2k}\in K).
\]
In fact such a form is determined by its coefficients $a_0, a_{-2}, a_{-4},\ldots$. Moreover, given a polynomial $P$ there is a unique form which has Laurent expansion starting with $P(z^{-2})dz$. The space of odd differential forms on $U_0$ has basis 
\[
\omega_k\eq\frac{x^k dx}{y} \qquad (k\geq 0).
\]

An odd function on $U_0$ is a function of the form $Q(x)y$. It has expansion
\[
\sum_{k\in\Z} a_{2k-1} z^{2k-1}\qquad (a_{2k-1}\in K).
\]
Such a function is determined by its coefficients $a_{-3}, a_{-5},\ldots$. Conversely, for each polynomial $P$ there exists a form which has Laurent expansion starting with $z^{-3}P(z^{-2})$. 

It follows that the space of odd forms modulo the space of exact odd forms is $2$-dimensional with basis 
\[
\omega\eq\frac{dx}{2y},\qquad \eta\eq\frac{x dx}{2y}.
\]
This space is canonically isomorphic to the first de Rham cohomology group of $E$. We denote it by $H^1_K$.

\begin{prop} 
Consider the map which sends an odd differential form $\kappa$ to the following element of $K[E_2]$:
\[
\mu(\kappa) \;:=\; \res(\kappa v).
\]
This map vanishes on exact forms. Therefore it defines a map from $H^1_K$ to $K[E_2]$.
\end{prop}
\begin{proof}
Indeed, if $f=Q(x)y$, then
\begin{align*}
\res(v df) 
&\eq -\res(f dv)
\eq\res(Q(x) y (x-\frac{E_2}{12}) dz)\\
&\eq\frac 12 \res(Q(x)(x-\frac{E_2}{12}) dx) 
\eq 0.
\end{align*}
\end{proof}

One has
\[
\mu(\omega)\eq1,\qquad \mu(\eta)\eq\frac{E_2}{12}.
\]
Therefore using $\mu$ one can build an isomorphism of algebras over $K$:
\[
\Sym H^1_K/(\omega - 1) \xrightarrow[\sim]{} K[E_2].
\]
The symbol $\Sym H^1_K$ denotes the algebra of symmetric tensors of $H^1_K$. 

For any $G_m$-module $M$ which is $G_m$-equivariant over some $G_m$-field we denote by $M(1)$ the same module but with the twisted action. If $m^M$ is the action on $M$ then the action on $M(1)$ is defined as follows:
\[
m^{M(1)}_\lambda a \eq \lambda m^M_\lambda a\qquad(a\in M,\lambda\in G_m.)
\]

Using this notation the isomorphism constructed above can be made into a $G_m$-equivariant isomorphism:
\[
\mu:\Sym H^1_K(1)/(\omega - 1) \xrightarrow[\sim]{} K[E_2].
\]

In the following two propositions the value of a quasi-modular form $f$ of weight $k$ on a pair of numbers $\omega_1, \omega_2\in\C$ with $\tau=\frac{\omega_2}{\omega_1}\in\HH$ is defined as follows:
\[
f(\omega_1, \omega_2) \;:=\; (2\pi \I)^k f(\tau) \omega_1^{-k}.
\]

\begin{prop}\label{prop:periods2}
Let $k=\C$. Let $a_0, b_0\in\C$ and $f\in K$ a rational function which is defined at the point $(a_0, b_0)$. Suppose $f$ has weight $k$. Choose a basis $\cc_1, \cc_2$ of the first homology for the curve $y^2=x^3+a_0 x+ b_0$. Let $\omega_i=\int_{\cc_i} \omega$ with $\tau=\frac{\omega_2}{\omega_1}\in\HH$. Then
\[
f(a_0,b_0) \eq \mu(f)(\omega_1, \omega_2).
\]
\end{prop}

\begin{prop}
Let $k=\C$. Let $a_0, b_0\in\C$ and $[\kappa]\in H^1_K$ represented by an odd differential form $\kappa$ which is defined at the point $(a_0, b_0)$. Suppose $\kappa$ has weight $k$. Choose a basis $\cc_1, \cc_2$ of the first homology for the curve $y^2=x^3+a_0 x+ b_0$. Let $\omega_i=\int_{\cc_i} \omega$ with $\tau=\frac{\omega_2}{\omega_1}\in\HH$. Then
\[
\int_{\cc_1} \kappa \eq \omega_1 \mu(\kappa)(\omega_1, \omega_2).
\]
\end{prop}

\section{Derivations on modular forms}\label{deriv-mf}
Consider the Gauss-Manin connection on the module $H^1_K$. It is a map
\[
\nabla: H^1_K \ra \Omega^1(K/k)\otimes H^1_K.
\]
This map is equivariant with respect to the $G_m$-action. Consider $\nabla \omega$. This is an element of $\Omega^1(K/k)\otimes H^1_K$. Let us view this element as a map 
\[
	\nabla[\omega]:\Der(K/k)\ra H^1_K.
\]
Since both spaces have dimension $2$ and one can check that the map is surjective, it is an isomorphism. For example, one has
\[
\nabla[\omega](\delta_e) \eq -[\omega].
\]

We would like to compute $\nabla[\omega]$ for a general vector field $\partial$. Let $\partial\in \Der(K/k)$ be a derivation of the field $K$. By $\partial^*$ we denote a lift of $\partial$ to the affine set $U_0$. We assume that $\partial^*$ is even, so it is given by an even function $\partial^* x$ and an odd function $\partial^* y$. There is no canonical choice of this lift. On the other hand we let $\partial$ act on formal Laurent series in $z$ simply by setting $\partial z=0$. We obtain
\[
2 y (\partial y -\partial^* y) \eq (3 x^2 + a) (\partial x -\partial^* x).
\]
Therefore there exists a Laurent series $\alpha$ with the property
\[
\partial y \eq \partial^* y + \alpha y',\qquad \partial x \eq \partial^* x+ \alpha x',
\]
where $'$ denotes the derivative of a Laurent series with respect to $z$. Since $\partial$ commutes with $'$, we have
\[
2\partial y \eq \partial x' \eq (\partial^* x)' + \alpha x'' + \alpha' x'\eq (\partial^* x)' + 2\partial y - 2 \partial^* y + \alpha' x'.
\]
Therefore
\[
\alpha'x'\eq2\partial^* y - (\partial^* x)'.
\]
It is easy to see that the right hand side is a regular odd function on $U_0$, hence it is a product of $y$ and a polynomial of $x$. Since $x'=2y$, we obtain that $\alpha'$ is a polynomial of $x$. Consider the form
\[
\alpha' dz.
\]
It is a linear combination over $K$ of an exact form, $\omega$, and $\eta$. Therefore $\alpha$ is a linear combination of a regular odd function on $U_0$, $z$, and $v_0$. We obtain
\[
\partial x \eq P(x) + (A_\partial z + B_\partial v_0) x'\qquad (P\in K[x],\; A_\partial\in K,\; B_\partial\in K).
\]
Note that 
\[
\partial x \eq -\frac{\partial a}5 z^2-\frac{\partial b}7 z^4+O(z^{6}).
\]
Therefore $\deg P\leq 2$. Moreover, looking at the expansions
\[
\begin{array}{cccccccc}
x^2 & \eq & z^{-4} && -\frac{2a}{5} & - \frac{2b}{7} z^2 & +
    \frac{a^2}{15} z^4 & +\, O(z^6),\\ 
x & \eq & & z^{-2} & & -\frac{a}5 z^2 & -\frac{b}7 z^4 & +\,
    O(z^{6}),\\
v_0 x' & \eq & -2 z^{-4} && -\frac{8a}{15} & - \frac{22b}{35} z^2 & + \frac{2 a^2}{35} z^4 & +\,O(z^6),\\
z x' & \eq & & -2 z^{-2} & & -\frac{2 a}{5} z^2 & -\frac{4 b}{7} z^4 & +\,O(z^6),
\end{array}
\]
we see that 
\[
\partial x \eq A_\partial( z x'+ 2x) + B_\partial(v_0 x' + 2 x^2 + \frac {4a}3).
\]
The elements $A_\partial$, $B_\partial$ can be found from the equations:
\[
-\frac{\partial a}{5} \eq -\frac{4 a}{5} A_\partial - \frac{6 b}{5} B_\partial ,\;
-\frac{\partial b}{7} \eq -\frac{6 b}{7} A_\partial + \frac{4 a^2}{21} B_\partial,
\]
or
\[
\partial a \eq 4 a A_\partial + 6 b B_\partial,\qquad
\partial b \eq 6 b A_\partial - \frac{4 a^2}{3} B_\partial.
\]

We see that it is convenient to choose the following vector fields as a basis:
\begin{align*}
\delta_e: & \;\delta_e a \eq 4a,\; \delta_e b \eq 6b,\\
\delta_s: & \;\delta_s a \eq 6 b,\; \delta_s b \eq -\frac{4 a^2}{3}.
\end{align*}
The first one, $\delta_e$ is the Euler derivative, which was already defined. The second one is the Serre derivative. One can check that our definition coincides with the standard one, i.e.
\[
\delta_s E_4 \eq -\frac{E_6}{3},\qquad 
\delta_s E_6 \eq -\frac{E_4^2}{2}.
\]
In particular, we see that $\delta_e$ and $\delta_s$ define vector fields on $\spec R$. Choosing $\alpha_e=z$, $\alpha_s=v_0$ we also obtain liftings of these vector fields to regular vector fields on $U_0$:
\[
\delta_e^* \eq \delta_e - \alpha_e \frac{d}{dz},\qquad
\delta_s^* \eq \delta_s - \alpha_s \frac{d}{dz}.
\]
For reference we summarise the values of the derivations:
\begin{align*}
\delta_e^* z& \eq -z,& \delta_e^* v_0& \eq v_0,& 
\delta_e^* x& \eq 2x,& \delta_e^* y& \eq 3y, \\
\delta_s^* z& \eq -v_0,& \delta_s^* v_0& \eq -y - \frac{a}{3}z,&
\delta_s^* x& \eq 2 x^2 + \frac{4a}3,& \delta_s^* y& \eq 3xy, \\
z'& \eq 1,& v_0'& \eq -x,& x'& \eq 2 y, & y'& \eq 3 x^2 + a.
\end{align*}

It is also convenient to introduce differential forms $d_e$, $d_s$ on the base as the dual basis to the basis $\delta_e$, $\delta_s$.  Then we have
\begin{align}\label{dxdy}
dx &\eq 2y (dz + z d_e + v_0 d_s) + 2 x d_e + (2 x^2 + \frac{4a}3) d_s,\\
dy &\eq (3 x^2 + a) (dz + z d_e + v_0 d_s) + 3 y d_e + 3 x y d_s.
\end{align}

It is now easy to compute $\nabla_\partial [\omega]$:
\[
\nabla_\partial [\omega] \eq [d \partial^* z] \eq -[\alpha_\partial' dz],
\]
so we obtain:
\[
\nabla_{\delta_e}[\omega] \eq -[\omega], \qquad
\nabla_{\delta_s}[\omega] \eq [\eta].
\]
The derivatives of $[\eta]$ are given by
\[
\nabla_{\delta_e}[\eta] \eq [\eta], \qquad
\nabla_{\delta_s}[\eta] \eq \frac{a}{3}[\omega].
\]

\section{Representing cohomology classes}\label{repr-coh}
In this section we consider in details how to represent cohomology classes on elliptic curves using the language of hypercovers (Section \ref{hypercovers}). We would like to represent cohomology classes for families of elliptic curves by differential forms on the total space.  For an elliptic curve $E$ we consider two affine open sets $U_0$ and $U_1$. The set $U_0$ was mentioned before, it is the complement of the neutral element $[\infty]$ of $E$. The equation of $U_0$ is $y^2=x^3+ax+b$. The set $U_1$ can be chosen to be any affine open set which contains $[\infty]$. We will not fix $U_1$ since sometimes we need $U_1$ to be ``small enough''. The intersection is denoted $U_{int} = U_0\cap U_1$. The triple $U_0, U_1, U_{int}$ defines a hypercover of $E$ which is a particular case of a \v{C}ech hypercover, the corresponding abstract chain complex is the segment, $\Delta_1$.

\subsection{Hyperforms}
Next, $0$-hyperforms are triples $(f_0, f_1, f_{int})$ where $f_0$ and $f_1$ are functions on $U_0$, $U_1$ correspondingly, and $f_{int}$ is forced to be $0$.  $1$-hyperforms are triples $(\theta_0, \theta_1, \theta_{int})$ where $\theta_0$ and $\theta_1$ are differential $1$-forms on $U_0$, $U_1$ correspondingly, and $\theta_{int}$ is a function on $U_{int}$. $2$-hyperforms are triples $(\iota_0, \iota_1, \iota_{int})$ where $\iota_0$ and $\iota_1$ are differential $2$-forms on $U_0$, $U_1$ correspondingly, and $\iota_{int}$ is a differential $1$-form on $U_{int}$. Note that in the case of a single elliptic curve or relative forms for a family of elliptic curves $\iota_0$ and $\iota_1$ must be zero. However if we consider absolute differential forms on a family of elliptic curves this is not the case.

According to the formulae of Section \ref{hypercovers} the differentials are given as follows:
\begin{align*}
d(f_0, f_1, 0) &\eq (d f_0, d f_1, f_1-f_0),\\ 
d(\theta_0, \theta_1, \theta_{int}) &\eq (d\theta_0, d\theta_1, d\theta_{int}-\theta_1+\theta_0).
\end{align*}

\subsection{The class $[\omega]$.}
We are going to write down some representatives for $H^1(E, \C)$. Consider the differential form $\frac{dx}{2y}$. We want to write it down as a regular form on the Weierstrass family. For this we recall that we lifted two vector fields $\delta_e$ and $\delta_s$ on the base to regular vector fields $\delta_e^*$, $\delta_e^*$ on $U_0$. Consider an arbitrary vector field $\partial$ on the base. Expressing it as a linear combination of $\delta_e$ and $\delta_s$ and using the lifts $\delta_e^*$ and $\delta_s^*$ we construct a lift of $\partial$, denoted $\partial^*$. Then we have
\[
2 y \partial^* y \eq (3 x^2 + a) \partial^* x + x\partial a + \partial b.
\]
Therefore
\[
i_{\partial^*} (dx\wedge dy) \eq \partial^* x dy - \partial^* y dx \eq (-x\partial a - \partial b) \frac{dx}{2y} \; \mod da, db.
\]
We see that if we choose $\partial=-\frac\partial{\partial b}$ we obtain a differential form which is regular on $U_0$ and represents the same relative form as $\frac{dx}{2y}$.

More explicitly, one can check that
\[
\frac\partial{\partial b} \eq 12\Delta^{-1}(4 a \delta_s - 6 b \delta_e),\quad
\frac\partial{\partial a} \eq -12\Delta^{-1}(6 b \delta_s + \frac{4 a^2}{3} \delta_e).
\]
Then we compute
\[
i_{-\left(\frac\partial{\partial b}\right)^*} (dx\wedge dy) \eq 8 \Delta^{-1}((-12 a x^2 +18 bx- 8 a^2) dy + (18 a x y - 27 b y) dx).
\]

Let us compute the Laurent expansion of the form above. More generally, for any $\partial$ we have
\begin{multline*}
\partial^* x dy - \partial^* y dx \eq (\partial^* x y' - \partial^* y x') (dz + z d_e + v_0 d_s) + (\partial^* x \delta_e^*y - \partial^* y \delta_e^* x) d_e \\+ 
(\partial^* x \delta_s^*y - \partial^* y \delta_s^* x) d_s.
\end{multline*}
We have already seen that
\[
\partial^* x y' - \partial^* y x' \eq -x \partial a - \partial b.
\]
For the remaining part it is enough to compute
\[
\delta_s^* x \delta_e^* y - \delta_e^*x \delta_s^* y \eq 4 a y.
\]
Hence
\[
i_{-\left(\frac\partial{\partial b}\right)^*} (dx\wedge dy) \eq 
dz + z d_e + v_0 d_s - 12\Delta^{-1}(16 a^2 y d_e + 24 a b y d_s).
\]
Noting that $d a = 4 a d_e + 6 b d_s$ we obtain
\[
i_{-\left(\frac\partial{\partial b}\right)^*} (dx\wedge dy) \eq 
dz + z d_e + v_0 d_s - 48\Delta^{-1}a y da.
\]

We would like to construct forms which have as small order of pole at infinity as possible. Therefore we take as $\omega_0$ the form
\[
\omega_0 \eq i_{-\left(\frac\partial{\partial b}\right)^*} (dx\wedge dy) + 48\Delta^{-1}a y da \eq dz + z d_e + v_0 d_s.
\]

As the form $\omega_1$, since it is always possible to approximate power series in $z$ by rational functions, one can choose any form such that
\[
\omega_1 \eq dz + z d_e + O(z^N) d_e + O(z^N) d_s \;\text{for some $N>>0$.}
\]
Indeed, one may start with $t=x/y$. This is a function on a neighborhood of infinity. Then express $dt = \alpha \frac{dx}{2y} \;\mod da, db$, where $\alpha$ is a rational function which equals to $-1$ at infinity. Therefore $\alpha^{-1} dt = dz \mod da, db$. One can add correction terms of the form $f da$, $f db$ to construct $\omega_1$ as needed.

Next since we want $(\omega_0, \omega_1, \omega_{int})$ to be in the $F^1$ of the Hodge filtration we put $\omega_{int}=0$.

In this way we have constructed a hyperform which represents the cohomology class $[\omega]$:
\[
\omega \eq (\omega_0, \omega_1, 0) \eq (dz + z d_e + v_0 d_s, dz + z d_e + O(z^N) d_e + O(z^N) d_s, 0).
\]

\subsection{The class $[\eta]$}
To construct a representative of $[\eta]$ we consider $x \omega_0$. Its Laurent expansion is
\[
x\omega_0 \eq x dz + x z d_e + x v_0 d_s.
\]
Here it is possible to make the pole of the coefficient at $d_s$ smaller by adding functions regular on $U_0$. Indeed, it is easy to see that the Laurent series $x v_0 + y$ has no pole. Therefore we put
\[
\eta_0 \eq x \omega_0 + y d_s \eq x dz + x z d_e + (x v_0 + y) d_s.
\]

Next we choose $\eta_{int}$. Since we need $d \eta_{int} + \eta_0 - \eta_1 = 0 \mod da, db$ it is natural to take as $\eta_{int}$ some kind of a formal integral of $-\eta_0$. Therefore we choose $\eta_{int}$ as any function which satisfies
\[
\eta_{int} \eq v_0 + O(z^N).
\]

We finally choose $\eta_1$ as follows:
\[
\eta_1 \eq -\frac{a}3 z d_s + O(z^N) d_z + O(z^N) d_e + O(z^N) d_s.
\]

Then the hyperform $\eta$ is defined as
\begin{multline*}
\eta \eq (\eta_0, \eta_1, \eta_{int}) \eq (x dz + x z d_e + (x v_0 + y) d_s, \\
-\frac{a}3 z d_s + O(z^N) dz+ O(z^N) d_e + O(z^N) d_s, v_0 + O(z^N)).
\end{multline*}

\subsection{Gauss-Manin derivatives}
Let us now compute the Gauss-Manin derivatives for the hyperforms constructed above. First we compute $d\omega$:
\[
d\omega_0 \eq -d_e\wedge dz + x d_s \wedge dz \eq -d_e\wedge \omega_0 + d_s \wedge \eta_0 \;\mod d_e\wedge d_s,
\]
\begin{multline*}
d\omega_1 \eq -d_e\wedge dz + O(z^{N-1}) d_e\wedge dz + O(z^{N-1}) d_s\wedge dz \\
\eq -d_e\wedge \omega_1 + d_s \wedge \eta_1 + O(z^{N-1}) d_e\wedge dz + O(z^{N-1}) d_s\wedge dz \;\mod d_e\wedge d_s.
\end{multline*}
\[
d\omega_{int} - \omega_1 + \omega_0 \eq v_0 d_s + O(z^N) d_e + O(z^N) d_s.
\]
This shows that
\begin{multline*}
d\omega \eq -d_e \wedge \omega + d_s \wedge \eta 
\\+ (0, O(z^N) d_e\wedge dz + O(z^N) d_s\wedge dz, O(z^N) d_e + O(z^N) d_s) \;\mod d_e\wedge d_s.
\end{multline*}

Similarly
\[
d\eta_0 \eq x d_e\wedge dz + \frac{a}3 d_s\wedge dz \eq d_e \wedge \eta_0 + \frac{a}3 d_s \wedge \omega_0 \;\mod d_e \wedge d_s.
\]
\[
d\eta_1 \eq \frac{a}{3} d_s\wedge dz + O(z^{N-1})d_e\wedge dz + O(z^{N-1})d_s\wedge dz  \;\mod d_e\wedge d_s.
\]
\[
d\eta_{int} - \eta_1 + \eta_0 \eq v_0 d_e + O(z^{N-1}) d_e + O(z^{N-1}) d_s.
\]
Therefore
\begin{multline*}
d\eta \eq d_e \wedge \eta + \frac{a}{3} d_s \wedge \eta 
\\+ (0, O(z^{N-1}) d_e\wedge dz + O(z^{N-1}) d_s\wedge dz, O(z^{N-1}) d_e + O(z^{N-1}) d_s) \;\mod d_e\wedge d_s.
\end{multline*}

Thus we have proved:
\begin{prop}
For hyperforms $\omega$, $\eta$ defined above the following identities hold up to hyperforms of type $(0, O(z^{N-1}), O(z^{N-1}) d_e + O(z^{N-1}) d_s)$ and hyperforms from the piece $G^1$ of the filtration $G^\bullet$:
\begin{align*}
\nabla_{\delta_e} \omega &\eq -\omega, &\nabla_{\delta_s} \omega &\eq \eta,\\
\nabla_{\delta_e} \eta &\eq \eta, &\nabla_{\delta_s} \eta &\eq \frac{a}3 \omega.
\end{align*}
\end{prop}

\subsection{Computation of the Poincar\'e pairing} \label{subs:comp_poincare}
We would like to compute the Poincar\'e pairing of $\omega$ and $\eta$ to illustrate Corollary \ref{cor_3_1_23}. By the definition
\[
\langle \omega, \eta \rangle \eq \int_E [\omega \wedge \eta] \eq \int_{\Delta_E\subset E\times E} [\omega \times \eta].
\]
Therefore we need to compute the integral on the right hand side. The hyperform $\omega\times\eta$ has components indexed by pairs $\alpha, \beta$ where $\alpha, \beta \in \{0, 1, int\}$. According to Section \ref{products}) 
\[
\omega\times\eta \eq \begin{pmatrix} 
\omega_0 \times \eta_0 & \omega_0 \times \eta_1 & \omega_0 \times \eta_{int}\\
\omega_1 \times \eta_0 & \omega_1 \times \eta_1 & \omega_1 \times \eta_{int}\\
-\omega_{int}\times \eta_0 & -\omega_{int} \times \eta_1 & \omega_{int} \times \eta_{int}
\end{pmatrix}.
\]
Recall that $\omega_{int}=0$. Therefore the last row is zero.

We have
\[
\int_{\Delta_E} [\omega \times \eta] 
\eq \Tr^{\Int}_{\Delta_E} \omega \times \eta
\eq 2\pi\I\, \Tr_{\Delta_E} \omega \times \eta,
\]
where
\[
\Tr_{\Delta_E} \omega \times \eta \eq
\sum_{L\subset\{1,2\},|L|=1} \qquad \sum_{Z_\bullet\in\Fl_L(\Delta_E)} \res_{L, Z_\bullet} (\omega\times\eta)_{a_L(Z_\bullet)}.
\]
Clearly, we have only one flag in $\Fl_{L}(\Delta_E)$ for each $L\subset\{1,2\}$, $|L|=1$. The flag is $\Delta_E, [\infty\times\infty]$. Let us denote this flag $\fl$. 
We find $a_{\{1\}}(\fl) = [0,1]\times[1]$ and $a_{\{2\}}(\fl) = [0]\times [0,1]$. Therefore
\[
\Tr_{\Delta_E} \omega\times\eta = 2\pi\I (-\res_{\{1\}, \fl} \omega_{int} \eta_1 + \res_{\{2\}, \fl} \omega_0\eta_{int}).
\]
The first residue is zero since $\omega_{int}$ is zero. The second one equals to
\[
\res_{[\infty]\times[\infty]} \eta_{int}\omega_0.
\]
The residue in the last expression can be computed using Laurent series expansion:
\[
\res_{[\infty]\times[\infty]} \eta_{int}\omega_0 \eq \res v_0 dz \eq 1.
\]
Therefore we have proved that $\Tr_{\Delta_E} \omega\times\eta = 1$ and
\[
\langle \omega, \eta \rangle \eq 2 \pi \I.
\]

\chapter{Examples}
In the first section of this chapter we give examples of higher cycles on the square of the Weierstrass family. We construct two higher cycles $\cycle_1$, $\cycle_2$ and then show that their difference is annihilated by $12$ in the Chow group. The goal here is simply to give example of showing that two cycles are equivalent. The construction of these cycles was inspired by other constructions of higher cycles on products of elliptic curves in \cite{GL1}, \cite{GL2}. The author came up first with the cycle $\cycle_1$ and the main theorem was proved using the cycle $\cycle_1$. This cycle may seem better than $\cycle_2$ since it is composed of smooth varieties. Later the author found the more simple and natural construction of $\cycle_2$, but since it is composed of a singular variety (and using the singularity is crucial since we have a function which has at the same time zero and pole there), we needed to develop somewhat more involved machinery to deal with general varieties in the third chapter, though the computations with $\cycle_2$ given here became simpler than it was before with $\cycle_1$.

In the second section we compute the corresponding $D$-module and the invariant $\Psi'_{\alg}$ for the second cycle (we could use the first one, but due to the first section we know that the result would be the same up to torsion). We verify that the constructed cycle satisfies requirements for Theorem \ref{thm:anal-comp:9} for $\Gamma=PSL_2(\Z)$, $k=2$ and $z_0=\I$. This proves the algebraicity conjecture in this case for any second CM point $z$ which is not equivalent to $\I$ (Theorem \ref{thm:grfunc:1}). Then we give example of computing the value of the Green function which proves that
\[
\sqrt{28}\, \widehat{G}_2^{\HH/PSL_2(\Z)}\left(\frac{-1+\sqrt{-7}}{2}, \I\right) \;\equiv\;  8 \log(8-3\sqrt{7}) \mod \pi\I \Z
\]
and
\[
G_2^{\HH/PSL_2(\Z)}\left(\frac{-1+\sqrt{-7}}{2}, \I\right) =  \frac{8}{\sqrt{7}} \log(8-3\sqrt{7}).
\]

In the last section we obtain values of the constants $N_A=1$, $N_B=2$ and $N=2$ for the case $\Gamma=PSL_2(\Z)$, $k=2$.

\section{Examples of higher cycles}
In this section we construct some higher cycles in $Z^2(E\times E, 1)$ where $E$ is an elliptic curve.

\subsection{Notations}
We consider the Weierstrass family $y^2 = x^3 + ax +b$. Denote the total space by $\E$. A fiber, i.e., a single elliptic curve, is denoted by $E$.  We cover $E$ by two charts. The first one is 
\[
U_0 \eq E\setminus \{[\infty]\},
\]
which was mentioned before. The second one is 
\[
U_1 \eq E\setminus \{(y=0)\}.
\]
The coordinates on $U_1$ are $t, u$ and the equation is
\[
u = t^3 + a t u^2 + b u^3.
\]
The gluing maps are given as follows:
\[
U_0\ra U_1: \; (x,y) \ra (x y^{-1}, y^{-1}), \;\; U_1\ra U_0: \; (t,u) \ra (t u^{-1}, u^{-1}).
\]

\subsection{The first cycle}
The first cycle will be on $E\times E$ for $b\neq 0$. To denote the coordinates on the first $E$ in $E\times E$ we will use index $1$, for the second one we use $2$. So $U_0\times U_0$ is given by two equations in $4$ variables:
\[
y_1^2 = x_1^3 + a x_1 + b,\;\; y_2^2 = x_2^3 + a x_2 + b.
\]

Consider the subvariety $W$ which is the closure in $E\times E$ of the subvariety given by the equation $x_2 y_1 + \I x_1 y_2=0$. To compute the closure first consider $U_0\times U_1$. We have equations:
\[
y_1^2 = x_1^3 + a x_1 + b,\;\; u_2 = t_2^3 + a t_2 u_2^2 + b u_2^3,\;\; t_2 u_2^{-1} y_1 + \I x_1 u_2^{-1}=0.
\]
This gives 
\[
x_1 = \I y_1 t_2,\;\; y_1^2 = -\I t_2^3 y_1^3 + \I a t_2 y_1 + b.
\]
when $(t_2, u_2) = 0$ we obtain $(x_1, y_1) = (0, \pm \sqrt{b})$. So we see that $W$ contains two more points on $U_0\times U_1$. In fact $t_2$ is a local parameter on $W$ at these points, $u_2$ has order $3$ and $x_1$ has order $1$.

Analogously $W$ contains two more points on $U_1\times U_0$, namely $(x_2, y_2) = (\pm \sqrt{b},0)$. 

It remains to look at $U_1\times U_1$. The equations are
\[
u_1 = t_1^3 + a t_1 u_1^2 + b u_1^3,\;\;u_2 = t_2^3 + a t_2 u_2^2 + b u_2^3,\;\; t_2 u_1^{-1} u_2^{-1} + \I t_1 u_1^{-1} u_2^{-1}=0.
\]
Therefore we have
\[
t_1 = \I t_2.
\]
One can check that the point $[\infty]\times[\infty]$ also belongs to $W$ and the local parameter there can be chosen as $t_1$ or $t_2$.

Let $f$ be the rational function on $W$ given as $y_1-\I y_2$. We first compute the divisor of $f$.  Since $y$ has a triple pole at $[\infty]$ and the projections $W\to E$ are unramified at $[\infty]$ we see that $f$ has triple pole at the points $[\infty]\times(0, \pm\sqrt{b})$ and $(0, \pm\sqrt{b})\times[\infty]$. To study the behaviour of $f$ at $[\infty]\times[\infty]$ we write $f$ as
\[
f = u_1^{-1} - \I u_2^{-1} = (u_2 - \I u_1) u_1^{-1} u_2^{-1}.
\]
Using the fact that $u = t^3 + a t^7 + b t^9 + \cdots$ on $U_1$ at $[\infty]$ we see that
\[
f\sim 2 b \I t_2^{3} \;\; \text{at $[\infty]\times[\infty]$,}
\]
so $f$ has a triple zero.

Finally we look for zeroes of $f$ on the set $U_0\times U_0$. We need to find all common solutions of the equations
\[
y_1^2 = x_1^3 + a x_1 + b,\;\; y_2^2 = x_2^3 + a x_2 + b,\;\;x_2 y_1 + \I x_1 y_2=0,\;\;y_1-\I y_2=0.
\]
We get
\[
y_1=\I y_2,\;\; (x_1+x_2) y_1 y_2 = 0.
\]
In the case $y_1=y_2=0$ we see that for $\lambda_1, \lambda_2, \lambda_3$~--- the distinct roots of the polynomial $x^3+ax+b$ the $9$ points $(\lambda_i,0)\times(\lambda_j,0)$ are solutions. The case $x_1+x_2=0$ does not give any solution unless $b=0$. We will not check that the $9$ zeroes are indeed simple since we already know that $f$ has total multiplicity of poles $12$ and triple zero was already found. Therefore
\begin{multline*}
\Div f = \sum_{i,j} (\lambda_i,0)\times(\lambda_j,0) + 3 [\infty]\times[\infty] - 3 [\infty]\times(0, \sqrt{b}) \\- 3 [\infty]\times(0, -\sqrt{b}) - 3 (0, \sqrt{b})\times[\infty] - 3 (0, -\sqrt{b})\times [\infty].
\end{multline*}

It is not difficult to correct $(W,f)$ and obtain a cycle. The following combination is a cycle in $Z^2(E\times E, 1)$:
\[
(W,f) - \sum_i((\lambda_i,0)\times E, y_2) - 3(E\times[\infty], y_1) + 3([\infty]\times E, x_2) + 3(E\times[\infty], x_1).
\]

\subsection{The second cycle}\label{second-cycle}
The second cycle is also on $E\times E$, but it is given by a single summand. We take $W$ as the closure of the subvariety defined by equation $x_1+x_2=0$. It is clear from the equation that when $x_1$ is infinite, $x_2$ must be infinite as well and vice versa. Therefore the closure contains only one additional point $[\infty]\times[\infty]$.

The problem with this $W$ is that it is singular. Namely on $U_1\times U_1$ the equations are
\[
u_1 = t_1^3 + a t_1 u_1^2 + b u_1^3,\;\;u_2 = t_2^3 + a t_2 u_2^2 + b u_2^3,\;\; t_1 u_1^{-1} + t_2 u_2^{-1}=0.
\]
So we obtain the following equations for $W$ inside $E\times E$:
\[
t_1 u_2 = -u_1 t_2,\;\; t_1^2+t_2^2 = 2 b t_2^2 u_1^2.
\]
One can see that the function $t_1 t_2^{-1}$ belongs to the integral closure of the structure ring, but does not belong to the structure ring itself. It is easy to check that if we pass to the normalization of $W$ we obtain $2$ points over $[\infty]\times[\infty]$, namely the one with $t_1 t_2^{-1} = \I$ and the one with $t_1 t_2^{-1} = -I$.

The function $f$ remains the same,
\[
f = y_1 - \I y_2.
\]

It is clear that $f$ does not have zeroes or poles on $W\cap U_0\times U_0$. The only remaining point is $[\infty]\times[\infty]$. We use expansion $u = t^3 + a t^7 + b t^9 + \cdots$.
The second equation for $W$ gives
\[
t_1^2 + t_2^2 = 2 b t_2^2 (t_1^6 + 2 a t_1^{10} + 2 b t_1^{12}+\cdots).
\]
This implies
\[
t_2^2 = -t_1^2 - 2 b t_1^8 - 4 ab t_1^{12} - 8 b^2 t_1^{14}+\cdots.
\]
There are $2$ solutions of this equation:
\[
t_2 = \pm \I (t_1 + b t_1^7 + 2 a b t_1^{11} + \frac{7}{2} b^2 t_1^{13}+\cdots).
\]
The corresponding value of $u_2$ can be computed:
\[
u_2 = \mp \I (t_1^3 + a t_1^7 + 2 b t_1^9+\cdots).
\]

We can now compute the expansion of $f$. Along the branch $t_2\sim \I t_1$ we have
\[
f = \frac{u_2 - \I u_1}{u_1 u_2} \sim 2 t_1^{-3}.
\]
Along the branch $t_2\sim -\I t_1$ we have
\[
f = \frac{u_2 - \I u_1}{u_1 u_2} \sim b t_1^3.
\]
Therefore $\Div f = 0$ and $(W,f) \in Z^2(E\times E, 1)$.

\subsection{Equivalence of the first and the second cycles}
Denote
\[
f_1 := x_2 y_1 + \I x_1 y_2,\;\;
f_2 := x_1 + x_2,\;\;
f_3 := y_1-\I y_2.
\]
Denote some algebraic sets
\begin{multline*}
D_1 := Z(f_1), \;\; 
D_2 := Z(f_2), \;\;
D_3 := Z(f_3), \;\;
Y_1 := Z(y_1), \;\;
Y_2 := Z(y_2), \\
X_1 := Z(x_1), \;\;
X_2 := Z(x_2), \;\;
Z_1 := [\infty]\times E, \;\;
Z_2 := E\times [\infty],
\end{multline*}
where $Z(f)$ denotes the closure of the zero locus of a function $f$.

Then we have
\begin{multline*}
\Div f_1 = D_1 - 3 Z_1 - 3 Z_2, \;\;
\Div f_2 = D_2 - 2 Z_1 - 2 Z_2, \;\;
\Div f_3 = D_3 - 3 Z_1 - 3 Z_2, \\
\Div y_1 = Y_1 - 3 Z_1, \;\;
\Div y_2 = Y_2 - 3 Z_2.
\end{multline*}

Using this notation the first cycle can be expressed as
\[
\cycle_1 = (D_1, f_3) - (Y_1, y_2) + 3(Z_1, x_2) + 3(Z_2, x_1 y_1^{-1}).
\]
The second cycle is
\[
\cycle_2 = (D_2, f_3).
\]

Let us write
\[
 \{f, g\} := (\Div f, g|_{\Div f}) - (\Div g, f|_{\Div g})\;\; \text{for functions $f$, $g$.}
\]
Note that all $\{f, g\}$ are trivial in the Chow group by definition.

We consider the following element:
\[
\Psi := \{f_1^3 f_2^{-3} y_1^{-1} y_2^{-1}, f_3\}.
\]
It can be expressed as follows:
\begin{multline*}
\Psi = 3 (D_1, f_3) - 3 (D_2, f_3) - (Y_1, f_3) - (Y_2, f_3) 
- (D_3, f_1^3 f_2^{-3} y_1^{-1} y_2^{-1}) \\
+ 3 (Z_1, f_1^3 f_2^{-3} y_1^{-1} y_2^{-1}) 
+ 3(Z_2, f_1^3 f_2^{-3} y_1^{-1} y_2^{-1}).
\end{multline*}
Now each term will be treated separately. The second and the third ones give
\[
- (Y_1, y_1 - \I y_2) - (Y_2, y_1 - \I y_2) = - (Y_1, y_2) - (Y_2, y_1) + (Y_1, \I).
\]
Using the fact that $y_1 = \I y_2$ on $D_3$ the fourth term equals to 
\[
- \left(D_3, \frac{(x_2 y_1 + \I x_1 y_2)^3}{y_1 y_2 (x_1+x_2)^3}\right) = - (D_3, -y_2) = -(D_3, y_2) + (D_3, -1).
\]
Computing asymptotics at infinity we get
\[
\left(Z_1, \frac{(x_2 y_1 + \I x_1 y_2)^3}{y_1 y_2 (x_1+x_2)^3}\right) 
= \left(Z_1, \frac{x_2^3}{y_2}\right),\;\;
\left(Z_2, \frac{(x_2 y_1 + \I x_1 y_2)^3}{y_1 y_2 (x_1+x_2)^3}\right) 
= \left(Z_2, -\I\frac{x_1^3}{y_1}\right).
\]
Therefore
\begin{multline*}
\Psi = 3 \cycle_1 - 3 \cycle_2 -(D_3, y_2) + 2 (Y_1, y_2) - (Y_2, y_1) - 3 (Z_1, y_2) + 6 (Z_2, y_1) \\+ (Y_1 + 2 D_3  + Z_2, \I).
\end{multline*}
The last term is torsion so we could neglect it but for curiosity we will track it further.

Next put
\[
\Psi' = \{(y_1 - \I y_2)^2 x_2^{-3}, y_2\}.
\]
We have
\[
\Psi' = 2 (D_3, y_2) - 6 (Z_1, y_2) - 3 (X_2, y_2) - 2(Y_2, y_1) + 3(Y_2, x_2) + 3(Z_2, -1).
\]
Therefore
\begin{multline*}
2\Psi + \Psi' = 6 (\cycle_1 - \cycle_2) + 4 (Y_1, y_2) - 4 (Y_2, y_1) - 12(Z_1, y_2)+ 12 (Z_2, y_1) \\
+ 3 (Y_2, x_2) - 3(X_2, y_2) + (Y_1, -1).
\end{multline*}

Next we put
\[
\Psi'' = \{y_1, y_2\} = (Y_1, y_2) - 3(Z_1, y_2) - (Y_2, y_1) + 3(Z_2, y_1).
\]
This gives
\[
2\Psi + \Psi' - 4\Psi'' = 6 (\cycle_1 - \cycle_2) + 3(Y_2,x_2) - 3(X_2, y_2) + (Y_1, -1).
\]

To kill the remaining terms (which are already {\em decomposable}) we take
\[
\sum_{i=1}^3 \{x_2-\lambda_i, \lambda_i\} 
= 2 (Y_2, x_2) - 2 (Z_2, \lambda_1 \lambda_2 \lambda_3) = 2(Y_2, x_2) - 2 (Z_2, -b),
\]
\[
\{x_2, b\} = 2 (X_2, y_2) - 2 (Z_2, b).
\]

Concluding,
\[
12 (\cycle_1 - \cycle_2) = 4\Psi + 2\Psi' - 8\Psi''.
\]
Therefore
\begin{prop}
The difference between the two cycles constructed above is torsion in the higher Chow group $CH^2(E\times E, 1)$. More precisely, it is annihilated by $12$.
\end{prop}

\section{Computation of the Abel-Jacobi map}

\subsection{General remarks}
Let us compute the associated invariants of the Abel-Jacobi map for the cycles constructed above. We will use the second cycle $(W, f)$. We need to compute $\Psi_1$ for certain hyperforms (see Section \ref{subs:ext_dmod}). The map $\Psi_1$ is defined as an integral over $W$. The situation is similar to the one in Section \ref{subs:comp_poincare}.

Since the only intersection of $W$ and $Z_1$ or $Z_2$ is at $[\infty]\times[\infty]$ we have only one flag $\fl = (W, [\infty]\times[\infty])$ to consider. For a hyperform $\theta$ with components $\theta_{0,0}$, $\theta_{0,1}$, \ldots, $\theta_{\Int, \Int}$ we obtain
\[
\int_W \theta = 2\pi\I (\res_{[\infty]\times[\infty]} \theta_{\Int, 1} + \res_{[\infty]\times[\infty]} \theta_{0,\Int}).
\]

Let us compute $\Psi_1$ for the following hyperforms:
\[
\theta^0 := \omega\times\omega,\;\; 
\theta^1 := \eta\times\omega+\omega\times\eta,\;\; 
\theta^2 := \eta\times\eta.
\]

\subsection{Some Laurent series expansions}
Let us use as the coordinate $z = z_1$, which is the formal integral of $\frac{d x_1}{y_1}$. The equation of $W$ becomes
\begin{multline*}
z^{-2} - \frac{a}5 z^2 - \frac{b}{7} z^4 + \frac{a^2}{75} z^6 + \frac{3ab}{385} z^8 + \cdots \\=
-z_2^{-2} + \frac{a}5 z_2^2 + \frac{b}{7} z_2^4 - \frac{a^2}{75} z_2^6 - \frac{3ab}{385} z_2^8 + \cdots
\end{multline*}
We can solve this equation inverting both sides and considering the right hand side as a power series in $z_2^2$:
\[
z_2^2 = -z^2 - \frac{2b}{7} z^8 + \frac{4ab}{55} z^{12}+\cdots.
\]
Then we find $z_2$ by taking the square root:
\[
z_2 = \pm\I(z + \frac{b}{7} z^7 - \frac{2ab}{55} z^{11}+\cdots).
\]
Let us call the case with plus ``case 1'' and the case with minus ``case 2''. For the function $f$ we obtain
\begin{align}\label{eq:laurent_ser4}
f &= -2 z^{-3} - \frac{2a}{5} z + \frac{3b}{7} z^3 + \frac{2a^2}{25} z^5 - \frac{53 ab}{385} z^7 + \cdots \qquad &\text{in case 1,}\\ \label{eq:laurent_ser5}
f &= -b z^3 + \frac{ab}{5} z^7 + \cdots \qquad &\text{in case 2.}
\end{align}
In fact we know that the product of the two series above must be $2b$ since $(y_1-\I y_2)(y_1+\I y_2)=2 b$ on $W$. Therefore we can get a better approximation for case 2 using case 1:
\[
f = -b z^3 + \frac{ab}{5} z^7 - \frac{3b^2}{14} z^9 - \frac{2a^2 b}{25} z^{11} + \frac{17 a b^2}{110} z^{13}+ \cdots \qquad \text{in case 2.}
\]
For computation of $\phi_f$ we need the logarithmic derivative of $f$. 
\begin{multline*}
\frac{df}{f} = (-3 z^{-1} + \frac{4a}{5} z^3 - \frac{9b}{7} z^5 - \frac{12a^2}{25} z^7 + \frac{86ab}{77} z^9 + \cdots) dz \\
+ (\frac{4a}{5} z^4 - \frac{9b}{7} z^6 - \frac{12a^2}{25} z^8 + \frac{86 a b}{77} z^{10}+\cdots) d_e\\
+ (\frac{6b}{5} z^4 + \frac{2a^2}{7} z^6 - \frac{18ab}{25} z^8 + \frac{-172 a^3 + 774 b^2}{1155} z^{10}+\cdots) d_s \qquad \text{in case 1,}
\end{multline*}
\begin{multline*}
\frac{df}{f} = (3 z^{-1} - \frac{4a}{5} z^3 + \frac{9b}{7} z^5 + \frac{12a^2}{25} z^7 - \frac{86ab}{77} z^9 + \cdots) dz \\
+ (6 - \frac{4a}{5} z^4 + \frac{9b}{7} z^6 + \frac{12a^2}{25} z^8 - \frac{86 a b}{77} z^{10}+\cdots) d_e\\
+ (-\frac{4 a^2}{3 b} - \frac{6b}{5} z^4 + \frac{2a^2}{7} z^6 + \frac{18ab}{25} z^8 - \frac{172 a^3 + 774 b^2}{1155} z^{10}+\cdots) d_s \qquad \text{in case 2.}
\end{multline*}

It is clear that $\Psi_1(\omega\times \omega) = 0$. We turn to the hyperform 
\[
\theta^1:=\eta\times\omega + \omega\times\eta.
\]

\subsection{The hyperform $\theta_1$}
There are two components which we need, namely $\theta^1_{0,\Int}$ and $\theta^1_{\Int, 1}$. Denote by $\omega^{(i)}$, $\eta^{(i)}$, $v_0^{(i)}$ the corresponding objects coming from the $i$-th curve ($i=1,2$). We can compute (up to $O(z^N)$):
\begin{align*}
\theta^1_{0,\Int} &= \eta_{\Int}^{(2)} \omega_0^{(1)} = v_0^{(2)} (d z_1 + z_1 d_e + v_0^{(1)} d_s),\\
\theta^1_{\Int,1} &= - \eta_{\Int}^{(1)} \omega_1^{(2)} = - v_0^{(1)} (d z_2 + z_2 d_e).
\end{align*}

Since in the end we will take the sum of residues on a curve at the same point, we can restrict our considerations to the sum 
\[
\theta^1_s:=\theta^1_{0,\Int} + \theta^1_{\Int,1}.
\]

We compute the corresponding expansion in case 1.
\begin{multline*}
\theta^1_s = \I(-2 z^{-1} -\frac{2a}{15} z^3 - \frac{6b}{7} z^5 + \frac{2 a^2}{525} z^7 + \frac{62 a b}{231} z^9 + \cdots) dz \\
+\I(-2-\frac{2 a}{15} z^4-\frac{6 b}{7} z^6+\frac{2 a^2}{525} z^8 + \frac{62 a b}{231} z^{10} + \cdots) d_e \\
+ \I(-z^{-2} - \frac{2a}{15} z^2 + \frac{b}{7} z^4 + \frac{299 a^2}{1575} z^6 - \frac{64 a b}{1155} z^{8} + \cdots) d_s.
\end{multline*}
It is easy to see that for case 2 the expansion is precisely the negative of this one.

Finally we take
\begin{multline*}
\frac{df}{f}\wedge \theta^1_s = (-6 \I z^{-1} + \cdots) d_e\wedge dz + (-3 \I z^{-3} + \frac{2 \I a}{5} z + \cdots) d_s \wedge dz \\+ (\cdots) d_s\wedge d_e \qquad \text{in case 1,}
\end{multline*}
\begin{multline*}
\frac{df}{f}\wedge \theta^1_s = (6 \I z^{-1} + \cdots) d_e \wedge dz + (-3 \I z^{-3} -\frac{8 \I a^2}{3b} z^{-1} + \frac{2 \I a}{5} z + \cdots) d_s \wedge dz \\+ (\cdots) d_s\wedge d_e \qquad \text{in case 2.}
\end{multline*}
This means that
\begin{align*}
\phi_f \theta^1_s &= d_e \otimes (-6 \I z^{-1} + \cdots) dz + d_s \otimes (-3 \I z^{-3} + \frac{2 \I a}{5} z + \cdots) dz \qquad \text{in case 1,}\\
\begin{split}
\phi_f \theta^1_s &= d_e\otimes (6 \I z^{-1} + \cdots) dz + d_s\otimes (-3 \I z^{-3} -\frac{8 \I a^2}{3b} z^{-1} + \cdots) dz \qquad \text{in case 2.}
\end{split}
\end{align*}

Hence, taking residues and adding,
\[
\Psi_1(\theta^1) = -\frac{8 \I a^2}{3 b} d_s.
\]

\subsection{The hyperform $\theta_2$}
Analogously to the previous case we first write down
\begin{align*}
\theta^2_{0,\Int} &= \eta_{\Int}^{(2)} \eta_0^{(1)} = v_0^{(2)} (x_1 d z_1 + x_1 z_1 d_e + (x_1 v_0^{(1)}+y_1) d_s),\\
\theta^2_{\Int,1} &= - \eta_{\Int}^{(1)} \eta_1^{(2)} = v_0^{(1)} ( \frac{a}3 z_2 d_s).
\end{align*}

Adding we obtain the following expansion (in case 1):
\begin{multline*}
\theta^2_s = \I(-z^{-3} + \frac{2a}{15} z + \frac{11 b}{35} z^3 + \cdots) dz
+ \I(-z^{-2} + \frac{2a}{15} z^2 + \frac{11 b}{35} z^4 + \cdots) d_e\\
+ \I(\frac{2a}{3} + \frac{2b}{5} z^2  + \frac{2 a^2}{315} z^4 + \cdots) d_s.
\end{multline*}

This gives
\[
\phi_f\theta^2_s = d_s\otimes (2 \I a z^{-1}+\cdots) dz + d_e \otimes (-3 \I z^{-3} + \frac{2 \I a}{5} z + \cdots) dz \qquad \text{in case 1,}
\]
\begin{multline*}
\phi_f\theta^2_s = d_s\otimes (-\frac{4 \I a^2}{3 b} z^{-3} + 2 \I a z^{-1}+\cdots) dz + d_e \otimes (3 \I z^{-3} - \frac{2 \I a}{5} z + \cdots) dz \\ \text{in case 2.}
\end{multline*}

Again we add residues to see that
\[
\Psi_1(\theta^2) = 4 \I a d_s.
\]

\subsection{The extension of $\D$-modules}
Using the computation of $\Psi_1$ above and since the basis of hyperforms was chosen so that their Gauss-Manin derivatives can be again expressed in this basis, we can explicitly describe the extension of $\D$-modules from Section \ref{subs:ext_dmod}. The $\D$-module $M$ is generated over $O_S$ by $1$, $\theta^0$, $\theta^1$, $\theta^2$. The action of the derivations $\delta_s$ and $\delta_e$ (denoted $\delta_s'$, $\delta_e'$) is defined as follows:
\begin{gather*}
\delta_e' 1 = 0,\;\;
\delta_e' \theta_0 = -2\theta_0,\;\;
\delta_e' \theta_1 = 0,\;\;
\delta_e' \theta_2 = 2\theta_2,\\
\delta_s' 1 = 0,\;\;
\delta_s' \theta_0 = \theta_1,\;\;
\delta_s' \theta_1 = 2\theta_2 + \frac{2a}{3} \theta_0 + \frac{8 \I a^2}{3 b} ,\;\;
\delta_s' \theta_2 = \frac{a}{3} \theta_1 - 4\I a.
\end{gather*}

Consider a curve of the following type:
\[
\Delta = -16(4 a^3 + 27 b^2) = \const.
\]
Note that $\delta_s$ is a derivation on such curve. For this curve the element $B\in M^4\otimes\D\otimes M^2$ must be of the form
\[
m \otimes (\delta_s^3- \frac{4a}3 \delta_s - 4b) \otimes [\omega]^2,\qquad m\in M^4.
\]
Here we note that $[\omega]$, the fiberwise cohomology class of $\omega$, is a modular form of weight $1$. The canonical pairing identifies the image of $[\eta]$ in $H^1/F^1 H^1$ with $2\pi\I [\omega]^{-1}$. Therefore the Kodaira-Spencer map acts as
\[
KS: [\omega] \ra d_s \otimes 2\pi\I [\omega]^{-1}.
\]
Therefore the symbol of $\delta_s$ is $2\pi\I [\omega]^{-2}$. Hence
\[
B = (2\pi\I)^{-3} [\omega]^4 \otimes (\delta_s^3- \frac{4a}3 \delta_s - 4b) \otimes [\omega]^2.
\]

We evaluate $\Psi'_\alg$ on $B$. We have
\begin{multline*}
[\omega]^2 = \theta^0,\;\; 
\delta_s'[\omega]^2 = \theta^1,\;\;
\delta_s'^2[\omega]^2 = 2\theta_2 + \frac{2a}{3} \theta_0 + \frac{8 \I a^2}{3 b},\\
\delta_s'^3[\omega]^2 = \frac{4a}3 \theta_1 + 4b \theta_0 + 24 \I a + \frac{32\I a^4}{9 b^2}.
\end{multline*}
Therefore
\[
\Psi'_\alg(B) = (2\pi\I)^{-3} [\omega]^4 \left(24\I a + \frac{32\I a^4}{9 b^2}\right).
\]

The modular form as a function on $\HH$ which computes the periods of the cohomology class above can be found using Proposition \ref{prop:periods2}. It is the function
\[
2\pi\I \mu\left(24\I a + \frac{32\I a^4}{9 b^2}\right) = -\frac{\pi E_4(E_4^3-E_6^2)}{E_6^2} = -1728 \pi \frac{E_4(\tau)}{j(\tau)-j(\I)}.
\]

\subsection{Green's function}
The Abel-Jacobi map for the family of cycles we study is a function of $a$ and $b$. As we have proved, the derivative of $AJ_s^{2,1}[\cycle_2]$ (the symmetric part) is $0$ along the vector field $\delta_e$. Therefore $AJ^{2,1}[\cycle_2]$ is invariant under the action of the multiplicative group $G_m$, which sends $a,b$ to $\lambda^4 a, \lambda^6 b$. Restricted to the curve $\Delta=\const$ this means that
\[
AJ_s^{2,1}[\cycle_2](a,b) = AJ_s^{2,1}[\cycle_2](\veps a, b) = AJ_s^{2,1}[\cycle_2] \qquad (\veps^2+\veps+1=0).
\]
In particular $AJ^{2,1}[\cycle_2]$ depends only on the $j$-invariant
\[
AJ_s^{2,1}[\cycle_2](a,b) = AJ_s^{2,1}[\cycle_2](j),\qquad j=\frac{-2^{12} 3^3 a^3}{\Delta}.
\]
Therefore $AJ_s^{2,1}[\cycle_2]$ is constant along the fibers of $j$. We have proved that (see Section \ref{subs:anal_comp}, keep in mind that $k=2$, $n=2$)
\[
\Psi'_{\analyt}(B) = 2\pi\I\, \Psi'_{\alg}(B) = -\frac{2 \pi^2\I\, E_4(E_4^3-E_6^2)}{E_6^2}.
\]
We aim to prove that $g:=\Psi'_{\analyt}(B)$ is proportional to $g^{\HH/PSL_2(\Z)}_{2, \I}$.  Since $g$ has only poles at the orbit of $\I$ and is zero at the cusp, it is enough to study the expansion of $g$ at $I$. We have (note that $E_6(\I)=0$ and $E_2(\I) = \frac{3}{\pi}$ and use the Ramanujan's formulae for the derivatives of $E_2$, $E_4$, $E_6$)
\[
E_4(\tau) = E_4(\I) \left(1 + 2\I (\tau-\I) + \cdots\right),
\]
\[
E_6(\tau) = -\pi\I E_4(\I)^2(\tau-\I)\left(1 + \frac{7 \I(\tau-\I)}{2} +\cdots\right).
\]
Therefore
\[
\frac{E_4(E_4^3-E_6^2)}{E_6^2} = - \pi^{-2} ( (\tau-\I)^{-2} + \I (\tau-\I)^{-1} + O(1)),
\]
whereas
\[
Q_{\I}(\tau)^{-2} = \frac{-4}{(\tau-\I)^2(\tau+\I)^2} = (\tau-\I)^{-2} + \I(\tau-\I)^{-1} + O(1).
\]
This gives
\[
g = 2\I\, Q_{\I}(\tau)^{-2} + O(1),
\]
thereby proving that (see Theorem \ref{thm:global_study3} and note that $m=2$ for the point $\I$)
\[
\Psi'_{\analyt}(B) = -\I g^{\HH/PSL_2(\Z)}_{2, \I}(\tau).
\]

Therefore one may apply Theorem \ref{thm:anal-comp:9} ($\alpha=\frac12$).
This implies our first main theorem:
\begin{thm} \label{thm:grfunc:1}
For any $z\in\HH$ which is not equivalent to $\I$ let $a_0$, $b_0$ be such that the $j$-invariant of the curve $E_{a_0, b_0}$ (given by equation $y^2=x^3+a_0 x+b_0$) is $j(z)$. Then
\[
G_2^{\HH/PSL_2(\Z)}(z, \I) \eq 8 \Re (\tfrac{1}{\sqrt{-4}} AJ^{2,1}[\cycle_2(a_0, b_0)], Q_z),
\]
where $\cycle_2(a_0, b_0)$ is the higher cycle on $E_{a_0, b_0}\times E_{a_0, b_0}$ constructed in Section \ref{second-cycle} and $Q_z$ is the polynomial
\[
\frac{(X-z)(X-\zc)}{z-\zc}\in \Sym^2 H^1(E_{a_0, b_0}, \C)\subset H^2(E_{a_0, b_0}\times E_{a_0, b_0}, \C)
\]
(recall that we identify $H^1(E_{a_0, b_0}, \C)$ with the space of polynomials of degree $1$ of variable $X$).
For any CM point $z$ not equivalent to $\I$ the algebraicity conjecture is true and 
\[
\sqrt{-4D}\, \widehat{G}_2^{\HH/PSL_2(\Z)}(z, \I) \;\equiv\; 2 \log(\cycle_2(a_0, b_0)\cdot Z_z) \mod \pi\I \Z
\]
for any algebraic cycle $Z_z$ on $E_{a_0, b_0}\times E_{a_0, b_0}$ whose cohomology class is $2\sqrt{D}Q_z$.
\end{thm}

\subsection{Special values of the Green function}
The identity proved in the last section allows us to compute the values of the Green's function $G_k^{\HH/PSL_2(\Z)}(\tau, I)$ in points of complex multiplication. We consider in this section an example with $\tau = \tau_7=\frac{-1+\sqrt{-7}}2$ (in this section we denote the variable on the upper half plane by $\tau$).

The following elliptic curve corresponds to the point $\tau_7$:
\[
y^2 = x^3 - 35 x - 98.
\]
Let us denote $\mu = \tau_7$. It is a root of the polynomial $x^2+x+2$. Note that we have a decomposition
\[
x^3- 35 x - 98 = (x-7)(x+\mu+4)(x-\mu+3).
\]
The curve has the following endomorphism (we denote $\mu = \tau_7$, $\mu^2+\mu+2=0$):
\[
\phi(x, y) = \left(\mu^{-2}\left(x + \frac{(3\mu+5)^2}{x+\mu+4}\right), \mu^{-3} y\left(1 - \frac{(3\mu+5)^2}{(x+\mu+4)^2}\right)\right).
\]

It's graph $\Gamma_\phi$ has a class in $H^2(E\times E, \Z)$. This class can be represented as 
\[
[\Gamma_\phi] = c_1 [Z_1] + c_2 [Z_2] + c_3 [\Delta_E] + c_4 \frac{(X-\tau)(X-\ol\tau)}{\tau-\ol\tau} \qquad c_i\in \C.
\]
Here recall that we identify $\Sym^2 H^1(E,\C)$ with the space of polynomials of degree $2$ in variable $X$, in fact (see Section \ref{subs:anal_comp})
\[
\frac{(X-\tau)(X-\ol\tau)}{\tau-\ol\tau} = \frac{[\omega]\times[\ol\omega] + [\ol\omega]\times[\omega]}{4\I \vol_\omega E}.
\]

To find some of the constants above we compute intersections
\[
[\Gamma_\phi] \cdot [Z_1] = 1,\;\;
[\Gamma_\phi] \cdot [Z_2] = 2,\;\;
[\Gamma_\phi] \cdot [\Delta_E] = 4.
\]
Since we know that 
\[
[Z_1]\cdot [Z_2] = 1,\;\;
[Z_i]\cdot[Z_i]=0,\;\;
[\Delta_E]\cdot[Z_i] = 1,\;\;
[\Delta_E]\cdot[\Delta_E] = 0,
\]
we obtain
\[
c_1= \frac{5}{2},\;\; c_2 = \frac{3}{2},\;\; c_3 = -\frac{1}{2}.
\]
To find the last coefficient consider the following integral:
\[
\int_{\Gamma_\phi} (\omega \times \ol\omega + \ol\omega\times\omega) = \int_E (\omega \wedge \phi^* \ol\omega + \ol\omega\wedge\phi^*\omega).
\]
Since $\phi^* \omega = \mu \omega$, we obtain that the integral in question equals to the following one:
\[
(\ol\mu-\mu) \int_E \omega\wedge\ol\omega = 2\I (\mu-\ol\mu) \vol_{\omega} E.
\]
Therefore
\[
0=[\Gamma_\phi]\cdot[\Gamma_\phi] \eq
c_1 [Z_1]\cdot[\Gamma_\phi] 
+ c_2 [Z_2]\cdot[\Gamma_\phi]
+ c_3 [\Delta_E]\cdot[\Gamma_\phi]
+ c_4 \frac{\mu-\ol\mu}2,
\]
which implies $c_4 = \sqrt{-7}$.

We use these computations to express
\[
\frac{(X-\tau)(X-\ol\tau)}{\tau-\ol\tau} = \frac{1}{\sqrt{-7}} \left([\Gamma_\phi] - \frac{5}{2} [Z_1] - \frac{3}{2}[Z_2] + \frac{1}{2}[\Delta_E]\right).
\]
Therefore we can choose
\[
Z_\tau = 2[\Gamma_\phi] - 5 [Z_1] - 3[Z_2] + [\Delta_E],
\]
or $Z_\tau = [\Gamma_\phi] - [\Gamma_{\bar\phi}]$, as in the introduction. We will use the first expression.

Our next step is to employ Theorem \ref{thm:grfunc:1}. So we need to compute $Z_\tau\cdot \cycle_2$. For any curve $Z\subset E\times E$ which does not pass through $[\infty]\times[\infty]$ we have defined the intersection number $\cycle_2\cdot Z$ as the following element in $\C^\times$:
\[
\cycle_2\cdot Z = \prod_{p\in W\cap Z} f(p)^{\ord_p(W\cdot Z)}.
\]
Since this number depends only on the cohomology class of $Z$, the definition can be extended to curves which pass through $[\infty]\times[\infty]$. In particular, $Z_\tau \cdot \cycle_2$ also makes sense.

To compute $\cycle_2\cdot Z_1$ we cannot directly take the value of $f$ at the intersection $W\cdot Z$ since the only intersection of $W\cdot$ is at $[\infty]\times[\infty]$ and $f$ is not defined as this point. Instead we consider the deformation $Z_1(z) = \{z\}\times E$, where $z$ is a formal parameter and by the point $z$ we mean the point with $z_1=z$. Then there are two points of intersection with values of $f$ given by the series (\ref{eq:laurent_ser4}), (\ref{eq:laurent_ser5}). Therefore
\[
\cycle_2\cdot Z_1 = 2 b = -14^2.
\]

For the computation of $\cycle_2\cdot Z_2$ we take deformation $Z_2(z) = E\times \{z\}.$ The expansions of the branches of $f$ begin with $2\I z^{-3}$ and $b\I z^3$, so we obtain $-2b$ as the product,
\[
\cycle_2\cdot Z_2 = -2 b = 14^2.
\]

Let us turn to the computation of $\cycle_2\cdot \Delta_E$. There are now $3$ points of intersection, namely $[\infty]\times[\infty]$, $(0,\sqrt{b})\times (0,\sqrt{b})$, and $(0,-\sqrt{b})\times (0,-\sqrt{b})$. The last two points are simple. To find the value at infinity we consider the deformation $\Delta_E(z)=\Delta_E+(z,0)$. We obtain two points of intersection. The coordinates of the first one satisfy $z_2 \sim \I z_1$, $z_1-z_2\sim z$. In this case $z_1\sim \frac {1+\I}{2} z$, hence $f\sim (4+4 \I)z^{-3}$. The second point satisfies $z_2 \sim -\I z_1$, $z_1-z_2\sim z$. In this case $z_1\sim \frac{1-\I}{2} z$ and $f\sim \frac{1+\I}{4} b z^3$. Multiplication of all the principal parts gives
\[
\cycle_2\cdot \Delta_E = (4+4 \I) \frac{1+\I}{4}\, b \, (\sqrt{b}-\I\sqrt{b})\,
(-\sqrt{b}+\I\sqrt{b}) = -4 b^2 = -14^4.
\]

The computation of $\cycle_2\cdot\Gamma_\phi$ is more difficult. Consider the equation
\[
\mu^{-2}\left(x_1 + \frac{(3\mu+5)^2}{x_1+\mu+4}\right) = -x_1.
\]
It is equivalent to $x_1$ being a root of the quadratic equation
\[
t^2+t (\mu+4) -7 \mu - 21 = 0.
\]
Let us denote the solutions of the quadratic equation above by $t_1$, $t_2$. For each $t=t_i$ we put $x_1=t$ and get an equation on $y_1$ as follows:
\[
y_1^2 = t^3 - 35 t - 98.
\]
In this way we obtain $4$ points of intersection $W\cap \Gamma_\phi$. Using the equation on $t$ we transform
\[
t^3 - 35 t - 98 = 14 \mu(t+6\mu+3).
\]
Also it is easy to see that
\[
\frac{3\mu+5}{t+\mu+4} = t \frac{-\mu^2-1}{3\mu+5} = \frac{t}{3-\mu}.
\]
Therefore for each point of intersection we have
\begin{multline*}
y_2 = y_1 \mu^{-3}\left (1-\frac{t^2}{(3-\mu)^2}\right) = y_1 \mu^{-3} \left(1+\frac{t (\mu+4) - 7(\mu+3)}{7-7\mu}\right) \\
= - \frac{y_1(t - 3\mu - 5)}{8\mu+4}.
\end{multline*}

It is not hard to see that at these four points the intersection is transversal. The product of values of $f$ over these four points is equal to
\begin{multline*}
(t_1^3-35 t_1-98) (t_2^3-35 t_2-98) \left(1+ \frac{\I}{8\mu+4} (t_1 - 3 \mu - 5)\right)^2\\
\times \left(1+ \frac{\I}{8\mu+4} (t_2 - 3\mu - 5)\right)^2
\end{multline*}
The product is in fact a product of norms from $\Q(\mu, t, \I)$ to $\Q(\mu, \I)$ (we denote this norm by $N$):
\[
=\frac{\mu^2}{64} N(t+6\mu+3) N(t-3\mu-5 - \I(8\mu+4))^2.
\]
We have
\[
N(t+6\mu+3) = (6\mu+3)^2 - (6\mu+3)(\mu+4)-7\mu-21 = -28(\mu+3),
\]
\begin{multline*}
N(t-3\mu-5 - \I(8\mu+4)) = 
(3\mu+5 + \I(8\mu+4))^2 \\ + (3\mu + 5 + \I(8\mu+4))(\mu+4)
-7\mu-21 =-28\mu(2\mu+1)(\I\mu+1).
\end{multline*}
So the product above equals
\[
=2^3 7^4 \mu (I\mu+1)^2 = 14^4 (-1+\mu - 2\I - \I\mu)  = 14^4 u, \;\; \text{where}
\]
\[
u = -1+\mu - 2\I - \I\mu \;\; \text{($u$ is a unit.)}
\]

It remains to consider the intersection point at infinity.  This is done as in the previous cases. We consider the deformation $\Gamma_\phi(z) = \Gamma_\phi + (z,0)$. The intersection points satisfy $z_2\sim \mu (z_1 - z)$ and $z_2\sim \pm \I z_1$. The expansions of $f$ start with 
\[
f \eq -2 \left(\frac{\mu-\I}{\mu}\right)^3 z^{-3} + \cdots,\;\;
f \eq -b \left(\frac{\mu}{\mu+\I}\right)^3 z^3 + \cdots.
\]
Multiplication of the principal terms gives
\[
2b \left(\frac{\mu-\I}{\mu+\I}\right)^3 \eq - 2 b u^3.
\]

Therefore we have proved that
\[
\cycle_2 \cdot \Gamma_\phi \eq -2 b 14^4 u^4 \eq 14^6 u^4.
\]

Finally we obtain
\[
\cycle_2\cdot Z_\tau \eq \cycle_2 \cdot (2\Gamma_\phi - 5 Z_1 - 3 Z_2 + \Delta_E) \eq u^8.
\]
Theorem \ref{thm:grfunc:1} says in this case (note that $u^2 = \I (8-3\sqrt{7})$):
\begin{prop}\label{prop:spec-val-green:1}
\[
\sqrt{28} \widehat{G}_2^{\HH/PSL_2(\Z)}\left(\frac{-1+\sqrt{-7}}{2}, \I\right) \;\equiv\;  8 \log(8-3\sqrt{7}) \mod \pi\I \Z.
\]
Therefore
\[
G_2^{\HH/PSL_2(\Z)}\left(\frac{-1+\sqrt{-7}}{2}, \I\right) \eq  \frac{8}{\sqrt{7}} \log(8-3\sqrt{7}).
\]
\end{prop}

\section{The torsion}\label{sec:torsion}
We compute the constants $N_A$ and $N_B$ for the group $PSL_2(\Z)$, $k=2$. This corresponds to the representation $V_2$. To compute $N_B$ we need to find the group $H_0(PSL_2(\Z), V_2^\Z)$. The group $PSL_2(Z)$ is generated by two elements:
\[
S\eq \begin{pmatrix} 0 & -1\\ 1 & 0\end{pmatrix},\qquad
T\eq \begin{pmatrix} 1 & 1\\ 0 & 1\end{pmatrix}.
\]
For $v\in V_2^\Z$, $v = v_2 X^2 + v_1 X + v_0$ we have
\[
S v - v \eq (v_0 - v_2) (X^2-1)  - 2 v_1 X,\qquad T v - v \eq - 2 v_2 X + v_2 - v_1.
\]
This shows that
\[
H_0(PSL_2(\Z), V_2^\Z) \;\cong\; \Z/2\Z,
\]
so we can choose $N_B=2$.

To find $N_A$ we need to compute the parabolic cohomology group 
\[
H^1_{\parab}(PSL_2(\Z), V_2/V_2^\Z).
\]
So let $c$ be a parabolic cocycle with values in $V_2/V_2^\Z$, which is given by its values on $S$ and $T$. Since we know that the cocycle is parabolic, we can modify it by a coboundary to ensure that $c_T=0$. Let 
\[
v \eq c(S) \eq v_2 X^2 + v_1 X + v_0\qquad (v_0, v_1, v_2 \in \C/\Z).
\]
Then we have $c(ST) = c(S)$ and 
\[
v+ Sv \eq (v_0+v_2)(X^2+1),\quad v+ ST v + (ST)^2 v \eq (2 v_0 - v_1 + 2 v_2) (1+X+X^2).
\]
This shows that 
\[
v_0 \;\equiv\; - v_2 \mod \Z,\qquad v_1 \;\equiv\; 0 \mod \Z.
\]
Moreover, we are free to add to $c$ a coboundary which is zero on $T$. Therefore we can add 
\[
S (v_0) - v_0 = (X^1-1) v_0,
\]
which makes the cocycle zero. Therefore 
\[
H^1_{\parab}(PSL_2(\Z), V_2/V_2^\Z) \eq 0,
\]
and we can choose $N_A=1$. This gives $N=2$.

\bibliography{commons/refs}
\end{document}